\documentclass[final]{elsarticle}
\usepackage{amsmath,graphicx,latexsym,amssymb, amscd, psfrag,verbatim, amsfonts,amsthm}
\usepackage{amsfonts,bm,mathrsfs,subfigure}
\usepackage{algorithm,amsxtra,color}
\usepackage[algo2e,ruled,vlined]{algorithm2e}
\usepackage{graphicx}
\usepackage{setspace}
\usepackage{tabularx}
\usepackage{enumerate}
\usepackage{cases}
\usepackage{tikz}
\usepackage{pgfplots}
\pgfplotsset{compat=1.11}
\usepackage{float}
\usepackage{dcolumn}
\usepackage{booktabs}
\usepackage{algorithmicx,algpseudocode}
\usepackage[top=1in, bottom=1in, left=0.7in, right=1.25in]{geometry}
\usepackage{multirow}
\usepackage[pdftex, unicode=true,
            colorlinks,
            linkcolor=red,
            anchorcolor=blue,
            citecolor=green
            ]{hyperref}

\newtheorem{thm}{Theorem}[section]

\newtheorem{prop}[thm]{Proposition}
\newtheorem{remark}[thm]{Remark}

\newtheorem{assumption}{Assumption}

\newtheorem{alg}{Algorithm}[section]
\newcommand{\be}{\begin{equation}}
\newcommand{\ee}{\end{equation}}
\newcommand{\lb}{\left (}
\newcommand{\rb}{\right )}

\newcommand{\bff}{\bm f}
\newcommand{\bn}{\bm n}

\newcommand{\bx}{\bm x}
\newcommand{\bX}{\bm X}

\newcommand{\bu}{\bm u}
\newcommand{\bU}{\bm U}
\newcommand{\bv}{\bm v}

\newcommand{\bG}{\bm G}
\newcommand{\bV}{\bm V}

\newcommand{\bg}{\bm g}
\newcommand{\br}{\bm r}

\newcommand{\btau}{\bm \tau}
\newcommand{\bS}{\bm \sigma}

\newcommand{\bD}{\bm D}

\newcommand{\bQ}{\bm Q}
\newcommand{\bT}{\bm T}

\newcommand{\bR}{\bm R}
\newcommand{\bF}{\bm F}
\newcommand{\bI}{\bm I}

\newcommand{\bphi}{\bm \phi}

\newcommand{\bomega}{\omega}

\newcommand{\totalnorm}[1]{{\left\vert\kern-0.25ex\left\vert\kern-0.25ex\left\vert #1
    \right\vert\kern-0.25ex\right\vert\kern-0.25ex\right\vert}}

\begin{document}
\begin{frontmatter}

\title{A Unified-Field Monolithic Fictitious Domain-Finite Element Method
for Fluid-Structure-Contact Interactions and Applications to Deterministic Lateral Displacement Problems}
\author[Tongji]{Cheng Wang}
\ead{wangcheng@tongji.edu.cn}
\address[Tongji]{School of Mathematical Sciences, Tongji University, 1239 Siping Road, Shanghai 200092, China}

\author[unlv]{Pengtao Sun\corref{cor}}
\ead{pengtao.sun@unlv.edu}
\address[unlv]{Department of Mathematical Sciences,
University of Nevada Las Vegas, 4505 Maryland Parkway, Las Vegas, NV
89154, USA}

\author[ruixin]{Yumiao Zhang}
\ead{yumiao.zhang@raymind.com}
\author[PSU,KAUST]{Jinchao Xu}
\ead{xu@math.psu.edu}
\ead{jinchao.xu@kaust.edu.sa}
\address[PSU]{Department of Mathematics, Pennsylvania State University, University Park, PA 16802, USA}
\address[KAUST]{Computer, Electrical and Mathematical Science and Engineering Division, King Abdullah University of Science and Technology, Thuwal 23955, Saudi Arabia}

\author[ruixin]{Yan Chen}
\ead{yanc@raymind.com}

\author[ruixin]{Jiarui Han}
\ead{han@raymind.com}
\address[ruixin]{Shenzhen Raymind Biotechnology Co., Ltd, 4018 Jintian Road, Shenzhen, China}

\cortext[cor]{Corresponding author}

\begin{abstract}
Based upon two overlapped, body-unfitted meshes, a type of unified-field monolithic
fictitious domain-finite element method (UFMFD-FEM) is developed in this
paper for moving interface problems of dynamic fluid-structure
interactions (FSI) accompanying with high-contrast physical coefficients across the
interface and contacting collisions between the structure and
fluidic channel wall when the structure is immersed in the fluid.
In particular, the proposed novel numerical method consists of
a monolithic, stabilized mixed finite element method within the frame of
fictitious domain/immersed boundary method (IBM) for generic fluid-structure-contact
interaction (FSCI) problems in the Eulerian--updated Lagrangian description, while involving the
no-slip type of interface conditions on the fluid-structure
interface, and the repulsive contact force on the
structural surface when the immersed structure contacts the fluidic channel wall. The developed UFMFD-FEM for FSI or FSCI problems can deal with the structural motion
with large rotational and translational displacements and/or large deformation
in an accurate and efficient fashion, which are first validated by two benchmark
FSI problems and one FSCI model problem, then by experimental results
of a realistic FSCI scenario -- the microfluidic deterministic lateral displacement
(DLD) problem that is applied to isolate circulating tumor cells (CTCs) from blood cells
in the blood fluid through a cascaded filter DLD microchip in practice,
where a particulate fluid with the pillar obstacles effect in the fluidic channel,
i.e., the effects of fluid-structure interaction and structure collision,
play significant roles to sort particles (cells) of different sizes
with tilted pillar arrays. The developed unified-field, monolithic fictitious domain-based
mixed finite element method can be seamlessly extended to more sophisticated, high
dimensional FSCI problems with contacting collisions between the moving
elastic structure and fluidic channel wall.
\end{abstract}

\begin{keyword}
Fluid-structure-contact interactions (FSCI) \sep interface conditions
\sep repulsive contact force \sep unified-field monolithic fictitious domain-finite element method (UFMFD-FEM)
\sep mixed finite element
\sep deterministic lateral displacement (DLD).

\vspace{.3cm} \MSC[2020] 65M22 \sep 65M60 \sep 65M85 \sep 65Z05 \sep
65D17 \sep 70F35 \sep 70F40 \sep 74S05 \sep 74F10 \sep 76M10 \sep
76M30 \sep 76D05 \sep 76D09
\end{keyword}

\end{frontmatter}

\section{Introduction}
\begin{figure}[htb]
  \centering
  \includegraphics[height=1.3in,width=1.6in]{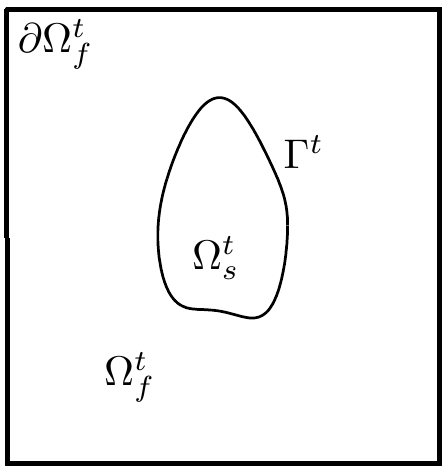}
  \hspace{0.5in}
  \includegraphics[height=1.3in,width=1.6in]{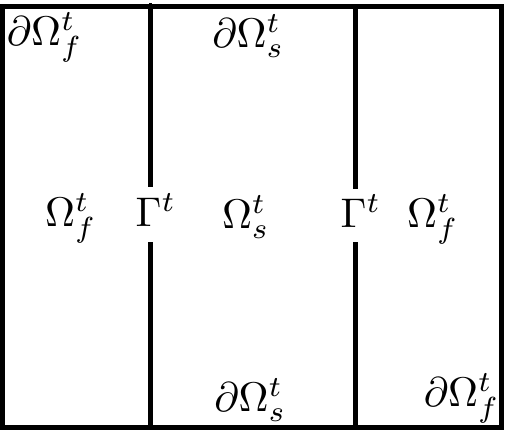}
  \hspace{0.5in}
  \includegraphics[height=1.3in,width=1.6in]{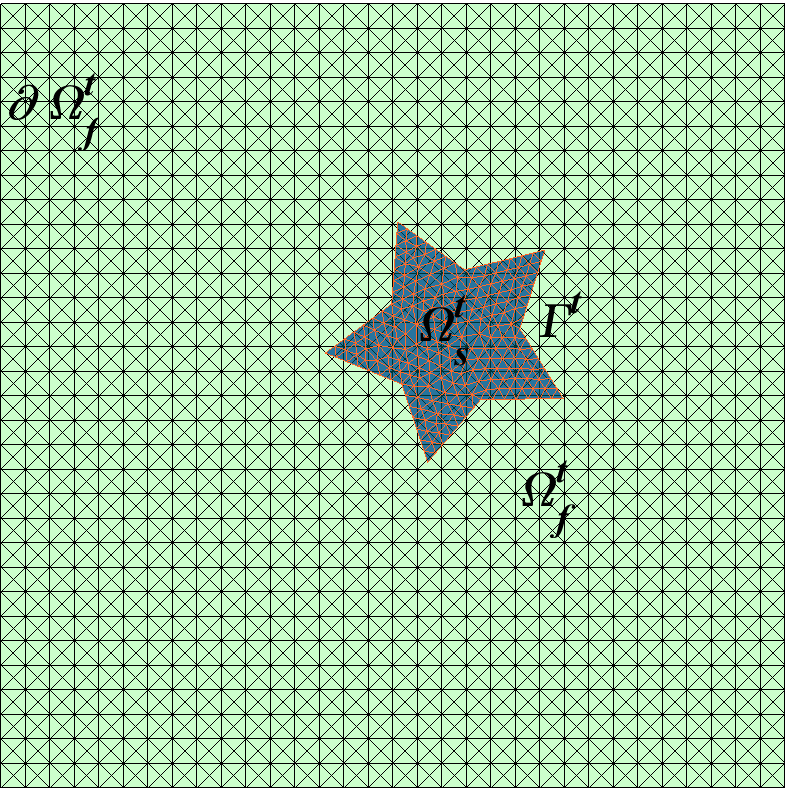}
  \caption{Illustrations of ({\bf Left:}) the immersed case of FSI;
  ({\bf Middle:}) the back-to-back case of FSI;
  ({\bf Right:}) two overlapped, body-unfitted meshes.
}
\label{fig:domain}
\end{figure}

The interaction of a flexible structure with a flowing fluid it is
submersed in or surrounded by gives rise to a rich variety of
physical phenomena with applications in many fields of engineering
and biology, e.g., the vibration of turbine blades impacted by the
fluid flow, the response of bridges and tall
buildings to winds, the floating parachute wafted by the air current,
the rotating mechanical parts of percussive drill tool driven by the pressurized liquid,
the flow of blood through the heart and arteries,
the artificial heart pump and the intravascular stent,
and etc. These examples comprise many applications of
fluid-structure interaction (FSI) problems in hydrodynamics,
aerodynamics and hemodynamics
\cite{Wang;Sun2016,Takizawa2011b,Chakrabarti2005,Dowell2001,Morand1995,Seo2013,Loon2005}.
Figure \ref{fig:domain} shows schematic domains of FSI in two cases.
To understand the phenomena of FSI problems, it is necessary to find
an effective way to model and simulate both the fluid motion and
structural motion in a monolithic fashion by investigating the interactional
mechanism between them.
In general, FSI problems require the fluid and structure fields at
the common interface to share not only the same velocity but also
the common normal stress. There are currently two major numerical
approaches classified by how the treatment of interface
conditions of FSI and the fluid mesh are handled along the moving
interface.

The first major approach is called the arbitrary Lagrangian Eulerian
(ALE) technique
\cite{Wang;Sun2016,Hirt1974,Hughes1981,Huerta1988,Nitikitpaiboon1993,Souli2010,Hu1996,Johnson;Tezduyar1997},
which adapts the fluid mesh to accommodate the deformations of the
structure on the interface
and produces a body-fitted conforming mesh. Thus, interface
conditions can be naturally satisfied within an interface
location-dependent finite element space and its variational
formulation. However, if the structure motion involves large
deformations/displacements, then fluid elements tend to become
ill-shaped, and accuracy loss will be imminent. This is commonly
fixed with re-meshing techniques. However, the re-meshing process
could be complicated, time consuming, and inaccurate. The transfer
of solutions from the degenerated mesh to the new mesh may also
introduce artificial diffusions, causing further inaccuracy. Moreover,
the worst thing of the re-meshing process is that it breaks the mesh
connectivity property of ALE mapping along the time which is crucial
to carry out an optimal convergence analysis for ALE-based finite element methods.

Therefore, when the most advanced and best cultured ALE-based scheme
inevitably reaches the point where only re-meshing helps,
one might be tempted to turn over to the
body-unfitted mesh method that leads to the second major approach --
the fictitious domain method, which adopts the non-conforming mesh by extending the fluid from its domain
$\Omega_f^t$ to the interior structural domain $\Omega_s^t$, simultaneously,
adding force-equivalent terms to the fluid equation
as a representation of the interactional mechanism of the fluid and structure.
Thus, re-meshing the fluid domain is avoided even for the case of large structural
displacements/deformations. As shown in
the right of Figure \ref{fig:domain}, the fictitious domain method requires
two overlapped, body-unfitted meshes, in which one is the background fixed (Eulerian)
fluid mesh, and the other one is the foreground moving (Lagrangian)
structural mesh.
A noticeable contribution to this approach is the immersed boundary
method (IBM) \cite{Peskin1972,peskin2002immersed} for the structures
that do not occupy volumes \cite{LeVeque1994,ZLi2001,Tan2009}. This
introduces a pseudo delta function for communications between
Eulerian and Lagrangian quantities.
To accurately represent the interaction between a fluid and a bulk
structure described by detailed constitutive laws, the extended IBM
(EIBM)
\cite{LucyZhang2004,Liu2006a,Liu2006b,LZhang2007,LZhang2008,wang2006,wang2007}
was developed. In this method, a fictitious fluid is introduced to
cover the structure domain and thus the fluid is extended to the
entire domain. To enforce the interface conditions, the FSI force is
imposed not only on the fluid-structure interface as the IBM does,
but also on every grid node in the fictitious fluid through a pseudo
delta function. The fluid equations are then solved to yield the
velocity field throughout the entire domain, afterwards, the structural
(velocity/displacement) variables are gained through an interpolation
process from the fluid mesh to the structural mesh. The obtained structural
variables can then be substituted into the suitable structural
constitutive law to update the FSI force,
which in turn can be used by the fluid equation to find the new velocity of the fluid nodes, and so forth.
Note that the IBM/EIBM is a semi-explicit finite difference scheme
within a partitioned approach due to the explicit coupling between the fluid
and the structure, therefore this method can be error-prone and may
have stability problems, especially when the time step size is
relatively large. Moreover, the involvement of a global Dirac delta
distribution makes the analyses of well-posedness, stability and
convergence properties even more difficult for IBM/EIBM.

To conquer those numerical deficiencies of IBM/EIBM, a fully
implicit monolithic immersed domain approach has been developed with
the help of the Lagrange multiplier, named the distributed Lagrange
multiplier/fictitious domain (DLM/FD) method
\cite{Glowinski1999,Glowinski2001,Yu2004,Yu2005,XShi2005}, which has
recently gained considerable popularity in the simulation of
interface problems including FSI
\cite{Boffi2014,Boffi2015,Boffi2016,Wang2016}. Same with the EIBM,
the DLM/FD method also needs to extend the fluid to the structure
domain as the fictitious fluid.
However, instead of a $d$-dimensional Dirac delta function
used in the EIBM, the DLM/FD method introduces the Lagrange
multiplier (a pseudo body force) to enforce the fictitious fluid
inside the structure to deform as the structure by constraining
the fluid and structural velocities equal to each other
in the structural domain. This
results in a monolithic system which is essentially a saddle-point
problem in regard to the Lagrange multiplier and principle unknowns.
Simultaneously, this also nests a subproblem of saddle-point type
for fluid equations. Thus a nested saddle-point problem is formed by
the DLM/FD method for FSI, for which numerical analysis has long
been missing due to its complexity. Until recently, some efforts
including the authors' achievements have been made on the analyses
of the well-posedness, stability and convergence of the nested
saddle-point problem associated with the DLM/FD method
\cite{Boffi2016,Wang2016,Lundberg;Sun2016,Sun2016,Sun2017,Wang;Sun2018,WangSun2021},
which shows that the accuracy of DLM/FD finite element method is of first-order
in $L^2$ norm and half order in $H^1$ norm for linear finite element
approximation to interface problems with non-smooth solutions across
the interface, at most.
In summary, the aforementioned immersed/body-unfitted mesh methods
have become increasingly popular in FSI simulations due to their
efficiency, flexibility and capability on handling FSI problems with
complex structure motion that involves a large
deformation/displacement for which the ALE method usually fails, although
they usually lack the resolution near the interface. In addition,
these body-unfitted mesh methods solely rely on the helps of either
Dirac delta function or Lagrange multiplier, so their limited convergence
results also depend on the discretization accuracies of Dirac delta function/Lagrange multipliers
to which a poor approximation will result in a bad convergence of fictitious domain method.

Recently, Wang et al. in \cite{WANG20171146} develop a one-field
fictitious domain method without introducing the Lagrange multiplier
to weakly constrain the fictitious fluid velocity equals the structural
velocity in the structural domain, which are instead reinforced
to equal each other in a strong sense. Inspired by their work, in this paper
we propose a unified-field monolithic fictitious domain-finite element method (UFMFD-FEM)
based upon the weak form of the presented FSCI model in Eulerian-updated Lagrangian description,
while considering the structural collision effect between the structure and fluidic channel wall
in the FSI scenario that has been not only a challenging but also an open problem so far
in terms of the definition of contact conditions and the method of finding
and exerting the repulsive contact force on the immersed structure, and etc.
We demonstrate that the weak form of the developed fictitious domain method can
derive an equivalent strong form by means of a $d$-dimensional
global Dirac delta distribution function that can be discretized by the EIBM, directly,
which illustrates that our proposed fictitious domain method is
actually in a close relation with the classical EIBM but more beyond
since we develop a monolithic mixed finite element method for the presented FSCI model
within the frame of unified-field fictitious domain method,
in contrast with the classical EIBM that is mostly discretized by finite
difference scheme for the strong form of FSI model via the
Dirac delta function. In our proposed UFMFD-FEM, neither the Dirac delta function
nor the Lagrange multiplier is involved,
instead, a unified variable pair of the velocity and the pressure are defined
in the entire domain on the fixed fluid mesh nodes only, and are shared by
both the fluid and structural equations,
besides that all advantages of the original fictitious domain/IBM are inherited,
such as the interface-unfitted and time-independent background fluid mesh, on the top of which
a moving structural mesh is updated all the time through the material/Lagrangian mapping.
In the meanwhile, our proposed method can also effectively deal with
structural collisions with the fluidic channel wall while coexisting with
fluid-structure interactions by treating the repulsive contact force,
which is solved through a nonlinear iteration process, as a traction force
on the contacting surface of structure.

The proposed UFMFD-FEM will be validated by two benchmark FSI problems,
i.e., the motion of neutrally buoyant circular cylinders (i) in simple shear flow; and
(ii) in plane Poiseuille flow, where convergence properties of
the developed UFMFD-FEM are also investigated.
In addition, we will also apply the developed UFMFD-FEM to a self-defined
FSCI model problem for an investigation of structural collision influence
on the particle's migration, and eventually, to a realistic microfluidic deterministic lateral
displacement (DLD) problem for numerically checking the DLD critical
diameter by releasing particles of different size in the fluid flow
through a cascaded filter DLD microchip. If the particle size is smaller
than the critical diameter, then the particle is isolated from the main stream
and moves down to the lower row of micro-posts in the zigzag mode,
otherwise it travels in the bumping mode over the micro-posts array
(the inner wall of fluidic channel),
where besides interactions between the particle and the
fluid, the collisions amongst the particle and the micro-posts array
are also significant to the particle's migration, and further, to
the particle's isolation. Numerical results of our method have good
agreements with physical experiments on a realistic DLD microchip.

The structure of this paper is organized as follows. In Section 2,
we present a general FSCI model with fluid-structure interface conditions
and structural contact conditions.
Then, weak formulations of the presented FSCI model
in both Eulerian-updated Lagrangian and Eulerian-total Lagrangian forms
are defined in Section 3, where an IBM-based strong form is also derived
via Dirac delta function. We propose the UFMFD-FEM in Section 4,
and corresponding numerical algorithms and implementations in Section 5.
Numerical experiments are carried out to validate the developed numerical
method in Section 6 through two benchmark FSI problems, one self-defined
FSCI model problem, and one realistic FSCI problem occurring
in a cascaded filter DLD microchip. Finally, the concluding remarks
and future work are given in Section 7.

\section{Model descriptions}\label{sec:FSImodel}
Let $\Omega$ denote an open bounded domain in $\mathbb{R}^d$ with a
polygonal boundary $\partial\Omega$. For all $t\in[0,T]$ with
$T$ as the terminal time, $\Omega$ always consists of two
subdomains, the fluid domain $\Omega_f^t$ and the structure domain
$\Omega_s^t$ with a polygonal boundary $\partial\Omega_s^t$,
which are separated by the fluid-structure interface $\Gamma^t$ and
satisfy ${\Omega}={\Omega_f^t}\cup{\Omega_s^t}$,
$\Omega_f^t\cap\Omega_s^t=\emptyset$,
$\Gamma^t=\partial\Omega_f^t\cap\partial\Omega_s^t$.
\begin{figure}[hbt]
    \centering
    \includegraphics[height=4cm]{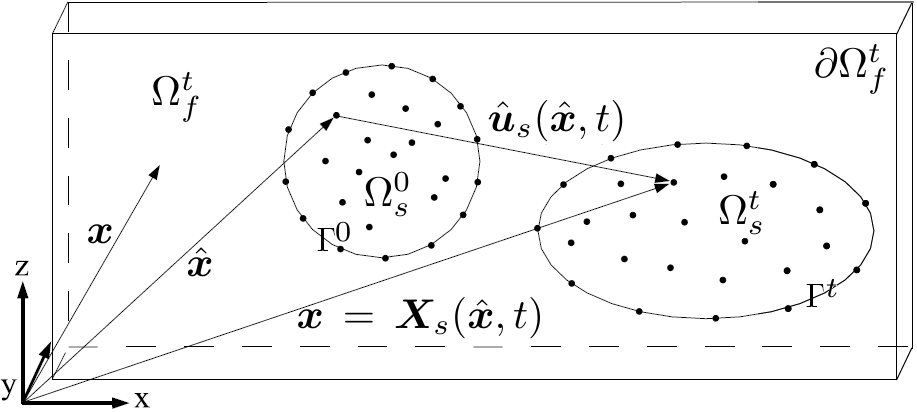}
    \caption{A schematic domain of FSI.}\label{fig:domainsketch}
\end{figure}
As shown in Figure \ref{fig:domainsketch}, the structure domain
$\Omega_s^t$ can be considered as an image of a reference
(initial/Lagrangian) domain $\Omega_s^0$ via a family of material/Lagrangian
mapping $\bX_s(t):\Omega_s^0\rightarrow\Omega_s^t$ defined as $\bx =
\bX_s(\hat\bx, t)$ for all $t\in[0, T]$, which is associated with
the material (Lagrangian) coordinate of points
$\hat\bx\in\Omega_s^0$ and corresponding spatial (Eulerian)
coordinate $\bx\in\Omega_s^t$. Using $\hat\bu_s(\hat\bx,t)$ to
denote the structural displacement defined in $\Omega_s^0$, we attain a flow map,
$\bx=\hat\bx+\hat\bu_s$. Correspondingly, we introduce the
deformation gradient tensor of structure,
$\bF=\nabla_{\hat\bx}\bx=\nabla_{\hat\bx}\bX_s(\hat\bx,t)=\bI+\nabla_{\hat\bx}\hat\bu_s$,
and the Jacobian $J=\operatorname{det}(\bF)$. In what follows, we
may use `` $\hat\cdot$ '' to denote a quantity `` $\cdot$ '' that is associated with the
reference domain of structure, $\Omega_s^0$, without further
explanation, and, we use
$(\phi,\tilde\phi)_\Psi=\int_\Psi\phi\tilde\phi d\bx$ and
$\langle\phi,\tilde\phi\rangle_{\partial\Psi}=\int_{\partial\Psi}\phi\tilde\phi
ds$ to denote a $L^2$ inner product inside a region
$\Psi\subset\mathbb{R}^d\ (d=2,3)$ and on a $(d-1)$-dimensional
region $\partial\Psi$, respectively.


\subsection{Fluid motion}\label{sec:fluideqn}
The fluid motion is described by the following incompressible
Navier-Stokes equations with respect to the fluid velocity $\bv_f$
and the fluid pressure $p_f$, 
\begin{equation}\label{fluid-equation}
\begin{array}{rcll}
\rho_f\frac{D\bv_f}{Dt}&=&\nabla\cdot\bS_f+\rho_f\bg, & \text{in } \Omega_f^t\times(0,T],\\
\nabla\cdot\bv_f&=&0, & \text{in } \Omega_f^t\times(0,T],
\end{array}
\end{equation}
where the stress tensor $\bS_f$ takes the form, $\bS_f=-p_f\bI+2\mu_f\bD(\bv_f)$,
the rate-of-strain tensor $\bD(\bv_f)=\frac{\nabla\bv_f+(\nabla\bv_f)^T}{2}$,
$\rho_f,\,\mu_f$ denote the density and dynamic
viscosity of the fluid, respectively,
$\frac{D}{Dt}$ represents the material derivative as
$\frac{D\bv_f}{Dt}=\frac{\partial\bv_f}{\partial
t}+\bv_f\cdot\nabla\bv_f$, and $\bg$ is the acceleration due to
gravity. (\ref{fluid-equation}) subjects to the following Dirichlet
boundary condition on the partial boundary $\partial\Omega_{f,D}^t$,
the possible Neumann boundary condition on the rest boundary
$\partial\Omega_{f,N}^t$, as well as the initial condition,
\begin{equation}\label{fluid-BICs}
\begin{array}{rcll}
\bv_f&=&\bv_{f,D}, & \text{on } \partial\Omega_{f,D}^t\times[0,T],\\
\bS_f\bn_{f,N}&=&\bff_{f,N}, & \text{on } \partial\Omega_{f,N}^t\times[0,T],\\
\bv_f(\bx,0)&=&\bv_{f}^0, & \text{in } \Omega_f^0,
\end{array}
\end{equation}
where $\bn_{f,N}$ is the outward normal unit vector on $\partial\Omega_{f,N}^t$.
Additionally, in the scenario of FSI problems, other than
(\ref{fluid-BICs}), the fluid equation (\ref{fluid-equation}) shall
also subject to the interface conditions defined in
(\ref{interface}) across the interface of fluid and structure,
$\Gamma^t$. In this case,
$\partial\Omega_f^t:=\partial\Omega_{f,D}^t\cup\partial\Omega_{f,N}^t\cup\Gamma^t$,
for all $t\in[0,T]$.

\subsection{Structure motion}\label{sec:solideqn}
We consider an incompressible
structure material which is, for instance, always adopted to model
the artery wall or blood cells in hemodynamic applications, and its
dynamics is defined below in terms of the structural velocity
$\bv_s$ and the
structural (hydrostatic) pressure $p_s$, 
\begin{equation}\label{solid-equation}
\begin{array}{rcll}
\rho_s\frac{D\bv_s}{Dt}&=&\nabla\cdot\bS_s+\rho_s\bg, & \text{in } \Omega_s^t\times(0,T],\\
\nabla\cdot\bv_s&=&0,& \text{in } \Omega_s^t\times(0,T],
\end{array}
\end{equation}
where $\rho_s$ is the structural density, $\bS_s$ denotes the Cauchy stress tensor of an incompressible
hyperelastic structure. For instance, if an incompressible
neo-Hookean (INH) material is adopted for the structure, then
$\bS_s=-p_s\bI+\btau_s$, and $\btau_s=\mu_s(\bF\bF^T-\bI)$ is the deviatoric stress,
here $\mu_s$ is the shear
modulus. Further, the structural hydrostatic pressure, $p_s$, plays a role of
Lagrange multiplier to reinforce the incompressibility condition of
structure, $J=1$ or $\nabla\cdot\bv_s=0$, mathematically
\cite{amabili_2018,HaoSun2021}. And, the structural velocity
$\bv_s(\bx,t)=\hat\bv_s(\hat\bx,t)=\frac{\partial\bX_s}{\partial
t}(\hat\bx,t)=\frac{\partial\hat\bu_s}{\partial t}$. Thus,
$\frac{D\bv_s}{Dt}=\frac{\partial^2\bX_s}{\partial t^2}(\hat\bx,t)$
in Lagrangian description.

The following Dirichlet boundary condition and initial condition can
be defined for (\ref{solid-equation}),
\begin{equation}\label{solid-BICs}
\begin{array}{rcll}
\bv_s&=&\bv_{s,D},& \text{on } \partial\Omega_s^t\backslash\Gamma^t\times[0,T],\\
\bv_s(\bx,0)&=&\bv_{s}^0,& \text{in } \Omega_s^0,
\end{array}
\end{equation}
which apply to some scenarios of FSI, e.g., when the structure is
not fully immersed in the fluid but owns a fixed outer boundary, as
shown in the middle of Figure \ref{fig:domain}. Otherwise
$\Gamma^t=\partial\Omega_s^t$, where if no collision occurs between
the immersed structure and the fluidic channel wall (see the left of Figure \ref{fig:domain}), then
the fluid-structure interface conditions (\ref{interface}) are
applied instead.
\begin{remark}
For the aforementioned material/Lagrangian mapping
$\bX_s(t):\Omega_s^0\rightarrow\Omega_s^t$, we further assume that
it fulfills the following conditions, $\bX_s(t)\in
W^{1,\infty}(\Omega_s^0)$, $\bX_s(t)$ is one to one, and there
exists a constant $\alpha$ such that $\|\bX_s(\hat\bx_1,
t)-\bX_s(\hat\bx_2, t)\|\geq\alpha\|\hat\bx_1-\hat\bx_2\|$ for all
$\hat\bx_1,\ \hat\bx_2\in\Omega_s^0$ and $t\in[0, T]$. Note that
this requirements imply that $\bX_s(t)$ is invertible with Lipschitz
inverse, which in particular implies that $\hat\bv_s(\hat\bx,t)\in
H^1(\Omega_s^0)^d$ if and only if $\bv_s =
\hat\bv_s\circ\bX_s^{-1}(t)\in H^1(\Omega_s^t)^d$
\cite{Gastaldi2001,Nobile;Formaggia1999}.
\end{remark}

\subsection{Fluid-structure interface conditions}\label{sec:interfacecond}
We consider the following
no-slip type interface conditions that can be applied to most cases
of FSI problems,
\begin{equation}\label{interface}
\begin{array}{rcll}
\bv_f&=&\bv_s, & \text{on } \Gamma^t\times[0,T],\\
\bS_f\bn&=&\bS_s\bn, & \text{on } \Gamma^t\times[0,T],
\end{array}
\end{equation}
which are called the kinematic and the dynamic lateral interface
condition describing the continuity of velocity, and the continuity
of normal stress, respectively,
where and in what follows, $\bn$ is the outward normal unit vector
on $\partial\Omega_s^t\cap\Gamma^t$ pointing into the fluid domain from the structural domain.

\subsection{Structural collision and contact conditions}\label{sec:contactcond}
\begin{figure}[hbt]
    \centering
    \includegraphics[height=5.2cm]{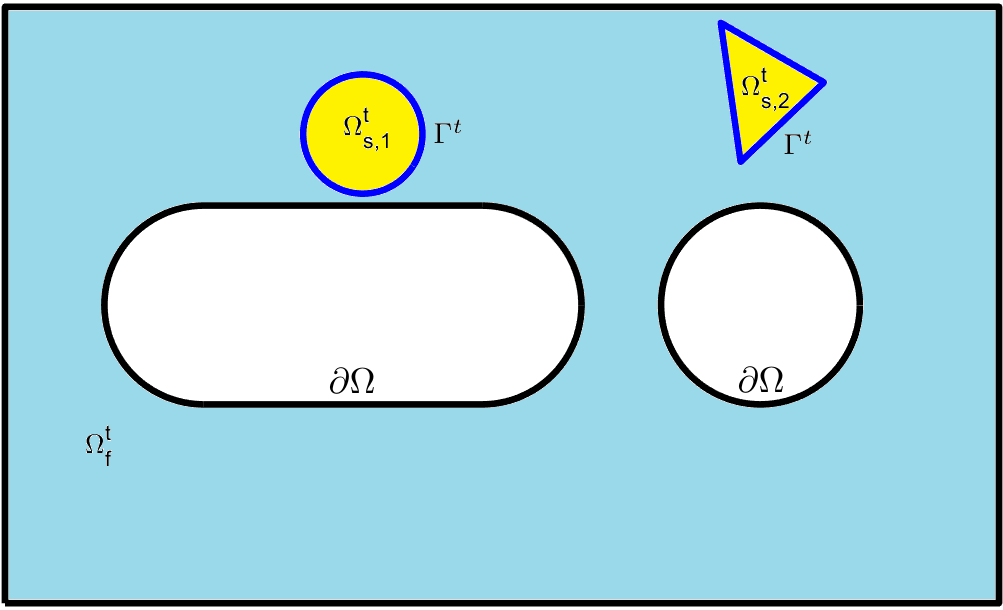}
    \includegraphics[height=5.2cm]{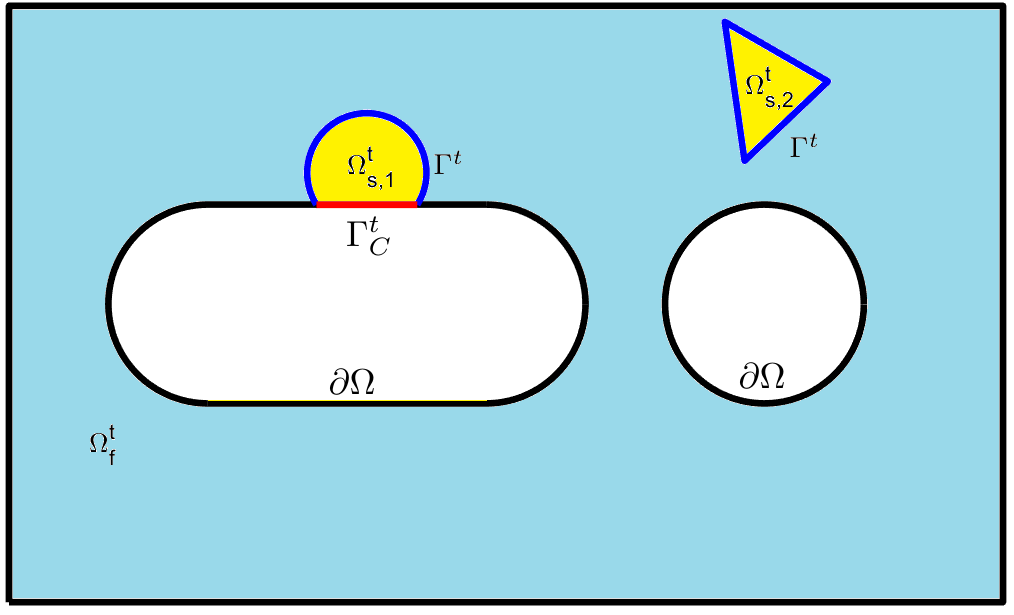}\\
    (a) \hspace{7cm} (b)
    \caption{Schematic illustrations of two types of structural boundary
        conditions: (a) no collision occurs and thus interface conditions are
        applied to the (blue) fluid-structure interface $\Gamma^t$;
        (b) the collision occurs and generates a repulsive contact force
        on the (red) contacting surface, $\Gamma_C^t$, between the structure and
        the fluidic channel wall, where the light blue area represents the fluid domain $\Omega_f^t$,
        the yellow areas are the immersed structures of different shapes, $\Omega_s^t$,
        and, the solid black lines denote the fluidic channel wall $\partial\Omega$.
    }\label{fig:contactfig_I}
\end{figure}

When a moving structure, which is immersed inside the fluid flow,
collides with the fixed wall of the fluidic channel $\partial\Omega$, the collision phenomenon
happens and generates a repulsive contact force acting on the
contacting surface $\Gamma_C^t$ that is actually a part of the
structural surface $\partial \Omega_s^t$.
To describe the mathematical model of structural collision, we
divide the structural boundary $\partial \Omega_s^t:=\Gamma^t$ into two parts:
\begin{align}
    \Gamma_{C}^t := &\left\{\bx\in\partial \Omega_s^t: \bx \not\in \mathring{\Omega}\right\},\label{eq:I_7}\\
    \Gamma_I^t := &\partial \Omega_s^t\backslash\Gamma_{C}^t,
\end{align}
where $\mathring{\Omega}$  denotes the interior of $\Omega$.
We note that $\Gamma_{C}^t$ can also be defined as $\Gamma_{C}^t:=\partial
\Omega_s^t\cap \partial \Omega$.
Thus, once $\Gamma_C^t\neq
\emptyset$ when the collision happens, we exert the contact force on
this part of structural boundary. Hence, two types of structural boundary
conditions can be raised as follows with respect to
two parts of $\partial \Omega_s^t:=\Gamma_I^t\cup\Gamma_{C}^t$,
respectively, as illustrated in Figure \ref{fig:contactfig_I}.
\begin{enumerate} 
    \item On $\Gamma_I^t$ where there is no collision, only the
    fluid-structure interface conditions \eqref{interface} are applied,
    i.e.,
    \begin{align}\label{struct-bc-1_I}
        \bv_f=\bv_s, \quad \bS_f\bn=\bS_s\bn,\quad \text{on }
        \Gamma_I^t\times[0,T].
    \end{align}
    \item On $ \Gamma_{C}^t$, the collision occurs by holding the following
    complementary conditions that are similar with, e.g., \cite{CFP2013},
    \begin{align}
        -\sigma_{s,\bn}\ge 0,\quad v_{s,\bn}\le 0,\quad \sigma_{s,\bn}v_{s,\bn} =&~ 0,\quad \text{on }
        \Gamma_C^t\times[0,T], 
        \label{struct-bc-2a}\\
        \bS_{s,\btau}= &~0,\quad \text{on }
        \Gamma_C^t\times[0,T], 
        \label{struct-bc-2b}
    \end{align}
    where $\sigma_{s,\bn}=\bn\cdot \bS_s \bn$, $v_{s,\bn}=\bv_s\cdot\bn$,
    $\bv_{s,\btau}=\bv_s-v_{s,\bn}\bn$ and
    $\bS_{s,\btau} = \bS_s\bn -\sigma_{s,\bn}\bn$. 
\end{enumerate}

Figure \ref{fig:contactfig_I}(b) illustrates a typical occurrence of
the collision between the immersed structure and the wall of fluidic channel,
All boundary/interface conditions with respect to the involved boundaries/interface
illustrated in Figure \ref{fig:contactfig_I} are further redescribed below.
\begin{enumerate} 
    \item On the fluidic channel wall $\partial\Omega$,
    the Dirichlet- and Neumann boundary conditions \eqref{fluid-BICs} are defined.
    \item On the fluid-structure interface $\Gamma^t$,
    the interface conditions \eqref{interface} are defined.
    \item On the contacting surface $\Gamma_C^t$ where the structural
    collision occurs, the contact conditions
    \eqref{struct-bc-2a} and \eqref{struct-bc-2b} are applied.
    Note that in this case $\Gamma_C^t$ is viewed as a part of structural
    boundary $\partial\Omega_s^t$ that is neither a part of the fluidic channel wall
    $\partial\Omega$ nor a part of the fluid-structure interface $\Gamma^t$.
    However, when no collision occurs, $\Gamma_C^t$ vanishes and retreats back to
    $\Gamma_I^t$ that is a part of $\Gamma^t$ essentially, as illustrated in Figure \ref{fig:contactfig_I}(a).
\end{enumerate}

By introducing a convex cone of admissible structural velocities which satisfy
the non-interpenetration condition on the contact surface $\Gamma_C^t$,
we might be able to obtain a weak form of the
aforementioned FSCI model \eqref{fluid-equation}-\eqref{interface},
\eqref{struct-bc-1_I}-\eqref{struct-bc-2b},
which is actually in the form of variational inequality
by following an analogous derivation in \cite{CFP2013,Han2022,Han2023}.
However at this stage, we do not intend to prove the well-posedness of
either this FSCI model or its relevant variational inequality that
involve the complemental contact condition \eqref{struct-bc-2a}.
In fact, the repulsive contact force on the contacting surface of
structure $\Gamma_{C}^t$ as well as $\Gamma_{C}^t$ itself are not
explicitly given in \eqref{struct-bc-2a}, which are actually determined by the
entire nonlinear FSCI system in an implicit fashion.
Inspired by similar arguments in \cite{Alart;Curnier1991,LAURSEN1992},
we make the following assumption in this paper in order to
process a successful numerical implementation later in Section \ref{sec:numericalexperiment}
and a possible numerical analysis in the future.
\begin{assumption}\label{assumption_1}
    The FSCI system  \eqref{fluid-equation}-\eqref{interface} and
    \eqref{struct-bc-1_I}-\eqref{struct-bc-2b} exists an unique solution, and, at any time when such a collision occurs,
    the contacting surface $\Gamma_C^t$ and the repulsive contact force $f_{s,C_W}:=-\sigma_{s,\bn}$ exist, uniquely.
\end{assumption}
Thus instead of  \eqref{struct-bc-2a}, we can apply the following structural boundary condition on $\Gamma_{C}^t$,
\begin{align}
    -\sigma_{s,\bn}= f_{s,C_W},\quad \text{on }
    \Gamma_C^t\times[0,T], 
    \label{struct-bc-2aNew}
\end{align}
which belongs to a Neumann-type of boundary condition. Therefore,
for the FSCI system defined by \eqref{fluid-equation}-\eqref{interface},
\eqref{struct-bc-1_I}, \eqref{struct-bc-2b} and \eqref{struct-bc-2aNew},
Assumption \ref{assumption_1} implies that if we utilize the uniquely existing
repulsive contact force $f_{s,C_W}$ and the contact surface $\Gamma_C^t$ in
\eqref{struct-bc-2b} and \eqref{struct-bc-2aNew}, then the solution of FSCI system
\eqref{fluid-equation}-\eqref{interface},
\eqref{struct-bc-1_I}, \eqref{struct-bc-2b} and \eqref{struct-bc-2aNew}
is also a solution of FSCI system  \eqref{fluid-equation}-\eqref{interface} and
\eqref{struct-bc-1_I}-\eqref{struct-bc-2b}, and  thus still satisfies the
complemental contact condition \eqref{struct-bc-2a}.
Thus, we gain $f_{s,C_W}\in [0,\infty)$,
$v_{s,\bn} \le 0$ and $\sigma_{s,\bn}v_{s,\bn} =0$,
which will be used to design a nonlinear iteration process
to determine $f_{s,C_W}$ and $\Gamma_C^t$
in the numerical implementation (see Algorithm \ref{algorithm3}).

\section{Weak formulations}

\subsection{Weak form of the fluid-structure-contact interaction model}\label{weakform-elastic}
Introduce the following notations of Sobolev spaces,
\begin{equation}\label{SobolevSpaces0}
\begin{array}{rcl}
\widetilde{\bV}&:=&\{(\bv_f,\bv_s)\big|\bv_f\in H^1(\Omega_f^t)^d,\,\bv_s\in H^1(\Omega_s^t)^d;
\,\bv_f=\bv_s\text{ on }{\Gamma^t},\bv_f=\bv_{f,D} \text{ on
}\partial\Omega_{f,D}^t,\bv_s=\bv_{s,D} \text{ on
}\partial\Omega_{s}^t\backslash\Gamma^t\},\\
\widetilde{\bV}_0&:=&\{(\bv_f,\bv_s)\big|\bv_f\in H^1(\Omega_f^t)^d,\,\bv_s\in H^1(\Omega_s^t)^d;
\,\bv_f=\bv_s\text{ on }{\Gamma^t},\bv_f=0 \text{ on
}\partial\Omega_{f,D}^t,\bv_s=0 \text{ on
}\partial\Omega_{s}^t\backslash\Gamma^t\},\\
\widetilde{W}&:=&\{(p_f,p_s)\big|p_f\in L^2(\Omega_f^t),\,p_s\in L^2(\Omega_s^t)\},
\end{array}
\end{equation}
based on which we can derive a natural weak formulation of the dimensional FSCI model
\eqref{fluid-equation}-\eqref{interface},
\eqref{struct-bc-1_I}, \eqref{struct-bc-2b} and \eqref{struct-bc-2aNew}
as follows: for any $
[(\tilde\bv_f,\tilde\bv_s),(\tilde p_f,\tilde p_s)]\in \widetilde{\bV}_{0} \times \widetilde{W}$, find $[(\bv_f,\bv_s),(p_f,p_s)]\in
H^1(0,T;\widetilde{\bV}) \times L^2(0,T;\widetilde{W})$ such that
\begin{eqnarray}
\left(\rho_f\frac{D{{\bv_f}}}{D
t},\tilde\bv_f\right)_{\Omega_f^t}+2\left(\mu_f\bD({\bv_f}),\bD(\tilde\bv_f)\right)_{\Omega_f^t}
-\left({p_f},\nabla
\cdot\tilde\bv_f\right)_{\Omega_f^t}+\left(\rho_s\frac{D{{\bv_s}}}{D
t},\tilde\bv_s\right)_{\Omega_s^t}+\left(\mu_s(\bF\bF^T-\bI),\nabla\tilde\bv_s\right)_{\Omega_s^t}\notag\\
-\left(p_s,\nabla \cdot\tilde\bv_s\right)_{\Omega_s^t}
=(\rho_f\bg,\tilde\bv_f)_{\Omega_f^t}
+(\rho_s\bg,\tilde\bv_s)_{\Omega_s^t}+\langle\bff_{f,N},\tilde\bv_f\rangle_{\partial\Omega_{f,N}^t}
-\langle f_{s,C_W},\tilde\bv_s\cdot \bn\rangle_{\Gamma_C^t},\label{weakform1-1}\\
\left(\nabla \cdot{{\bv_f}},\tilde p_f\right)_{\Omega_f^t}
+\left(\nabla\cdot{{\bv_s}},\tilde p_s\right)_{\Omega_s^t}=0,
\label{weakform1-2}
\end{eqnarray}
where the interface integral term vanishes from (\ref{weakform1-1}) due to
the first interface condition (\ref{interface})$_1$ that is built into the
space $\widetilde{\bV}_0$ as well as the second interface condition (\ref{interface})$_2$
when the integration by parts is applied to inner products of the fluidic and structural stress terms.
Here the repulsive contact force $f_{s,C_W}$ and the contact surface
$\Gamma_C^t$ are implicitly determined by a nonlinear iteration algorithm
as shown in Algorithm \ref{algorithm3}.

Because $\Omega=\Omega_f^t\cup\Omega_s^t$ for all $t\in[0,T]$, by
smoothly extending all fluid parts defined in $\Omega_f^t$ into the
structural domain $\Omega_s^t$ as a kind of fictitious fluid while
forcing the velocity to satisfy the interface condition (\ref{interface})$_1$
not only on the interface $\Gamma^t$ but also inside the entire structural
domain, $\mathring{\Omega}_s^t$, we can define the following Sobolev spaces of the unified
velocity $\bv$ and unified pressure $p$ in the entire domain $\Omega$, respectively,
\begin{equation}\label{SobolevSpaces}
\begin{array}{rcl}
\bV&:=&\{\bv\in H^1(\Omega)^d\big|\bv_f=\bv|_{\Omega_f^t}\in H^1(\Omega_f^t)^d,
\bv_s=\bv|_{\Omega_s^t}\in H^1(\Omega_s^t)^d;
\,
\bv=\bv_{f,D} \text{ on }\partial\Omega_{f,D}^t,\bv=\bv_{s,D} \text{ on
}\partial\Omega_{s}^t\backslash\Gamma^t\},\\
\bV_0&:=&\{\bv\in H^1(\Omega)^d\big|\bv_f=\bv|_{\Omega_f^t}\in H^1(\Omega_f^t)^d,
\bv_s=\bv|_{\Omega_s^t}\in H^1(\Omega_s^t)^d;
\,
\bv=0 \text{ on }\partial\Omega_{f,D}^t\cup\partial\Omega_{s}^t\backslash\Gamma^t\},\\
W&:=&\{p\in L^2(\Omega)\big|p_f=p|_{\Omega_f^t}\in L^2(\Omega_f^t),\,p_s=p|_{\Omega_s^t}\in L^2(\Omega_s^t)\},
\end{array}
\end{equation}
utilizing which we are able to rewrite (\ref{weakform1-1}) and (\ref{weakform1-2}) as the
following equivalent weak formulation in the fictitious domain frame,
i.e., for any $(\tilde\bv,\tilde p)\in {\bV}_{0} \times W$, find
$(\bv, p)\in H^1(0,T;{\bV}) \times L^2(0,T;W)$ such that
\begin{eqnarray}
\left(\rho_f\frac{D{{\bv}}}{D
t},\tilde\bv\right)_{\Omega}+2\left(\mu_f\bD({\bv}),\bD(\tilde\bv)\right)_{\Omega}
-\left({p},\nabla
\cdot\tilde\bv\right)_{\Omega}+\left((\rho_s-\rho_f)\frac{D{{\bv}}}{D
t},\tilde\bv\right)_{\Omega_s^t}+\left(\mu_s(\bF\bF^T-\bI),\nabla\tilde\bv\right)_{\Omega_s^t}\notag\\
-2\left(\mu_f\bD({\bv}),\bD(\tilde\bv)\right)_{\Omega_s^t}
=(\rho_f\bg,\tilde\bv)_{\Omega}
+\left((\rho_s-\rho_f)\bg,\tilde\bv\right)_{\Omega_s^t}+\langle\bff_{f,N},\tilde\bv\rangle_{\partial\Omega_{f,N}^t}
-\langle f_{s,C_W},\tilde\bv\cdot \bn\rangle_{\Gamma_C^t},\label{weakform2-1}\\
\left(\nabla \cdot{{\bv}},\tilde
p\right)_{\Omega}=0,\label{weakform2-2}
\end{eqnarray}
where all unified trial and test functions, $\bv,\,p,\,\tilde\bv$ and $\tilde p$
which are defined in the entire domain $\Omega$,
actually represent $\bv|_{\Omega_s^t},\,p|_{\Omega_s^t},\,\tilde\bv|_{\Omega_s^t}$ and
$\tilde p|_{\Omega_s^t}$ once they are involved in the inner products of
structural terms defined in $\Omega_s^t$.
Note that (\ref{weakform2-1}) is described in the Eulerian-updated Lagrangian frame since
the entire domain $\Omega$ is stationary thus the Eulerian
description is used therein to define the real and fictitious fluid equations,
whereas the updated Lagrangian description is applied to the structural domain
$\Omega_s^t$ that is kept being updated via its position variable $\bX_s(\hat\bx,t)$
along the time $t\in (0,T]$ determined by
\begin{equation}\label{position-strong}
\begin{array}{rcll}
\frac{\partial\bX_s}{\partial
t}(\hat\bx,t)&=&\bv(\bX_s(\hat\bx,t),t),& \text{in
}\Omega_s^0\times[0,T],\\
\bX_s(\hat\bx,0)&=&\bX_s^0,& \text{in }\Omega_s^0.
\end{array}
\end{equation}

\subsection{Derivation of an IBM-based strong form}
If we transfer all terms associated with $\Omega_s^t$ in (\ref{weakform2-1})
back to the initial domain of structure, $\Omega_s^0$,
then we have the following alternative weak formulation
in Eulerian-total Lagrangian description,
\begin{eqnarray}
\left(\rho_f\frac{D{{\bv}}}{D
t},\tilde\bv\right)_{\Omega}+2\left(\mu_f\bD({\bv}),\bD(\tilde\bv)\right)_{\Omega}
-\left({p},\nabla
\cdot\tilde\bv\right)_{\Omega}+\left((\rho_s-\rho_f)J\frac{\partial^2\bX_s}{\partial
t^2},\tilde\bv(\bX_s(\hat\bx,t))\right)_{\Omega_s^0}\notag\\
+\left(\mu_sJ(\bF-\bF^{-T}),\nabla_{\hat\bx}\tilde\bv(\bX_s(\hat\bx,t))\right)_{\Omega_s^0}
-2\big(\mu_fJ[\nabla_{\hat\bx}\bv(\bX_s(\hat\bx,t))\bF^{-1}+\bF^{-T}\nabla_{\hat\bx}\bv(\bX_s(\hat\bx,t))^T],\nabla_{\hat\bx}\tilde\bv(\bX_s(\hat\bx,t))\big)_{\Omega_s^0}\notag\\
=(\rho_f\bg,\tilde\bv)_{\Omega}
+\left((\rho_s-\rho_f)J\bg,\tilde\bv(\bX_s(\hat\bx,t))\right)_{\Omega_s^0}+\langle\bff_{f,N},\tilde\bv\rangle_{\partial\Omega_{f,N}^t}
-\langle Jf_{s,C_W}(\bX_s(\hat\bx,t)),\tilde\bv(\bX_s(\hat\bx,t))\cdot\hat\bn\rangle_{\Gamma_C^0},\label{weakform3-1}\\
\left(\nabla \cdot{{\bv}},\tilde
p\right)_{\Omega}=0,\quad \forall (\tilde\bv,\tilde p)\in \bV_{0} \times W,\label{weakform3-2}
\end{eqnarray}
where the structure domain remains as the initial one, $\Omega_s^0$,
and the one-to-one Lagrangian mapping, $\bX_s(t):\
\Omega_s^0\rightarrow\Omega_s^t$, is used to update the current
position of the structure, $\bx\in\Omega_s^t$, as
$\bx=\bX_s(\hat\bx,t)=\hat\bx+\hat\bu_s$ for any
$\hat\bx\in\Omega_s^0$. Thus the total Lagrangian description is
used in $\Omega_s^0$ while the entire domain $\Omega$ remains Eulerian.

Applying the defining property of the $d$-dimensional Dirac delta
distribution to any $\tilde\bv\in\bV_0$ as follows,
\begin{eqnarray}
\tilde\bv(\bX_s(\hat\bx,t))
=\int_{\Omega}\tilde\bv(\bx)\delta(\bx-\bX_s(\hat\bx,t))d\bx,\quad\forall
\hat\bx\in\Omega_s^0,\label{Dirac}
\end{eqnarray}
we can reformulate, for example, the fourth term on the left hand
side of (\ref{weakform3-1}) as follows,
\begin{eqnarray}
\left((\rho_s-\rho_f)J\frac{\partial^2\bX_s}{\partial
t^2},\tilde\bv(\bX_s(\hat\bx,t))\right)_{\Omega_s^0} &=&
\int_{\Omega_s^0}(\rho_s-\rho_f)J \frac{\partial^2\bX_s}{\partial
t^2}(\hat\bx,t)\left(\int_{\Omega}\tilde\bv(\bx)\delta(\bx-\bX_s(\hat\bx,t))d\bx\right)d\hat\bx\notag\\
&=& \int_{\Omega}\left(\int_{\Omega_s^0}(\rho_s-\rho_f)J
\frac{\partial^2\bX_s}{\partial
t^2}(\hat\bx,t)\delta(\bx-\bX_s(\hat\bx,t))d\hat\bx\right)\tilde\bv(\bx)d\bx\notag\\
&=&\left(\int_{\Omega_s^0}(\rho_s-\rho_f)J
\frac{\partial^2\bX_s}{\partial
t^2}\delta(\bx-\bX_s(\hat\bx,t))d\hat\bx,\tilde\bv\right)_{\Omega},
\end{eqnarray}
or, reformulate the more complicated fifth term on the left hand
side of (\ref{weakform3-1}) by integrating by parts first then using
(\ref{Dirac}) as follows,
\begin{eqnarray}
\left(\mu_sJ(\bF-\bF^{-T}),\nabla_{\hat\bx}\tilde\bv(\bX_s(\hat\bx,t))\right)_{\Omega_s^0}
&=& -\int_{\Omega_s^0}\nabla_{\hat\bx}\cdot\left(\mu_sJ(\bF-\bF^{-T})\right)\left(\int_{\Omega}\tilde\bv(\bx)\delta(\bx-\bX_s(\hat\bx,t))d\bx\right)d\hat\bx\notag\\
&&+ \int_{\partial\Omega_s^0}\left(\mu_sJ(\bF-\bF^{-T})\right)\hat\bn\cdot\left(\int_{\Omega}\tilde\bv(\bx)\delta(\bx-\bX_s(\hat\bx,t))d\bx\right)d\hat s \notag\\
&=&-\left(\int_{\Omega_s^0}\nabla_{\hat\bx}\cdot\left(\mu_sJ(\bF-\bF^{-T})\right)\delta(\bx-\bX_s(\hat\bx,t))d\hat\bx,\tilde\bv\right)_{\Omega}\notag\\
&&
+\left(\int_{\partial\Omega_s^0}\left(\mu_sJ(\bF-\bF^{-T})\right)\hat\bn\delta(\bx-\bX_s(\hat\bx,t))d\hat
s,\tilde\bv\right)_{\Omega}.
\end{eqnarray}
Similarly, we can also reformulate the rest structural terms in
(\ref{weakform3-1}). Then, conducting the integration by parts as
well for the second term on the left hand side of
(\ref{weakform3-1}), and considering that $\tilde\bv$ is arbitrary
in $\bV_0$, we can obtain the following strong form of the dimensional
momentum equation,
\begin{eqnarray}
&&\rho_f\frac{D{{\bv}}}{D
t}-\nabla\cdot\left(2\mu_f\bD({\bv})\right) +\nabla p
+\int_{\Omega_s^0}(\rho_s-\rho_f)J \frac{\partial^2\bX_s}{\partial
t^2}\delta(\bx-\bX_s(\hat\bx,t))d\hat\bx\notag\\
&&-\int_{\Omega_s^0}\nabla_{\hat\bx}\cdot\left(\mu_sJ(\bF-\bF^{-T})\right)\delta(\bx-\bX_s(\hat\bx,t))d\hat\bx
+\int_{\partial\Omega_s^0}\left(\mu_sJ(\bF-\bF^{-T})\right)\hat\bn\delta(\bx-\bX_s(\hat\bx,t))d\hat s\notag\\
&&+\int_{\Omega_s^0}\nabla_{\hat\bx}\cdot
\left(\mu_fJ[\nabla_{\hat\bx}\bv(\bX_s(\hat\bx,t))\bF^{-1}+\bF^{-T}\nabla_{\hat\bx}\bv(\bX_s(\hat\bx,t))^T]\right)\delta(\bx-\bX_s(\hat\bx,t))d\hat\bx\notag\\
&&-\int_{\partial\Omega_s^0}
\left(\mu_fJ[\nabla_{\hat\bx}\bv(\bX_s(\hat\bx,t))\bF^{-1}
+\bF^{-T}\nabla_{\hat\bx}\bv(\bX_s(\hat\bx,t))^T]\right)\hat\bn
\delta(\bx-\bX_s(\hat\bx,t))d\hat s\notag\\
&&=\rho_f\bg
+\int_{\Omega_s^0}(\rho_s-\rho_f)J\bg\delta(\bx-\bX_s(\hat\bx,t))d\hat\bx
-\int_{\Gamma_C^0}
Jf_{s,C_W}(\bX_s(\hat\bx,t))\delta(\bx-\bX_s(\hat\bx,t))d\hat
s,\quad \text{in }\Omega\times(0,T],\label{strongform1}
\end{eqnarray}
and the strong form of the dimensional mass equation due to the
arbitrary $\tilde p\in W$,
\begin{eqnarray}\label{strongform2}
\nabla\cdot\bv = 0,\quad \text{in }\Omega\times(0,T],
\end{eqnarray}
where $\bF=\nabla_{\hat\bx}\bX_s(\hat\bx,t)=\nabla_{\hat\bx}\bx(\hat\bx,t)$,
and $\bX_s(\hat\bx,t)$ satisfies (\ref{position-strong}). The strong form
(\ref{strongform1}) and (\ref{strongform2}) subject to the following boundary and
initial conditions,
\begin{equation}\label{BICs-strong}
\begin{array}{rcll}
\bv&=&\left\{
\begin{array}{l}
\bv_{f,D}, \\
\bv_{s,D},
\end{array} \right.
&\left.
\begin{array}{l}
\text{on } \partial\Omega_{f,D}^t\times[0,T],\\
\text{on } \partial\Omega_s^t\backslash\Gamma^t\times[0,T],
\end{array} \right.\\
(-p\bI+2\mu_f\bD({\bv}))\bn_{f,N}&=&\bff_{f,N}, & \hspace{0.2cm}\text{on } \partial\Omega_{f,N}^t\times[0,T],\\
\bv(\bx,0)&=&\left\{
\begin{array}{l}
\bv_{f}^0, \\
\bv_{s}^0,
\end{array} \right.
&\left.
\begin{array}{l}
\text{in } \Omega_f^0,\\
\text{in } \Omega_s^0.
\end{array} \right.
\end{array}
\end{equation}

The above derivations are completely reversible, i.e., we can trace
all the way back to regain (\ref{weakform3-1}) and
(\ref{weakform3-2}), which means (\ref{weakform3-1}) and
(\ref{weakform3-2}), or further back, (\ref{weakform2-1}) and
(\ref{weakform2-2}) are also the weak formulation of the strong
forms (\ref{strongform1})-(\ref{BICs-strong}). In
fact, (\ref{strongform1})-(\ref{BICs-strong}) generally
demonstrate the mathematical modeling of FSCI based on IBM. 
It is well known that the original IBM \cite{peskin2002immersed,peskin1993improved} was
carried out by means of the finite difference scheme to approximate the
strong form that involves the Dirac delta function everywhere, such as
(\ref{strongform1})-(\ref{BICs-strong}) here that requires the
construction of suitable approximations to the Dirac delta
distribution in order to discretize the terms in (\ref{strongform1}) that
contains the Dirac delta function. Our goal in this paper is to
propose a continuous Galerkin (CG) finite element approximation to
the fluidic and structure equations, as well as to the Lagrangian mapping $\bX_s$ that
describes the position of immersed structure along the time, within
the frame of unified-field monolithic fictitious domain method without
introducing the Dirac delta function, the Jacobian matrix $\bF$ and
the Jacobian $J$ into the structural terms.
To that end, we will develop our fictitious domain-finite element approximation
based on (\ref{weakform2-1}) and (\ref{weakform2-2}) in the next section,
thus there is no need to discretize the Dirac delta function and
no need to deal with the complicated nonlinear coefficients that
involve $\bF$ and $J$ in the structural terms.
All of these can be considered as significant improvements over the original IBM.

\section{Unified-field monolithic fictitious domain-finite element method}
In this section, we describe how the unified-field, monolithic fictitious domain-finite element
method (UFMFD-FEM) is developed on the basis of
the weak formulation (\ref{weakform2-1}) and (\ref{weakform2-2}).
First of all, we introduce a fixed Eulerian
mesh, ${\mathcal{T}}_{h}=\bigcup\limits_{k=1}^M e_k\ (0<h<1)$, to
triangulate the entire (real and fictitious) fluid domain $\Omega$ for up to $M$ simplicial
fluid elements $e_k\ (1\leq k\leq M)$, and adopt an updated
Lagrangian mesh, ${\mathcal{T}}_{h_s}^t=\bigcup\limits_{k=1}^{M_s}
e_{s,k}^t\ (0<h_s<1)$, to triangulate the structural domain
$\Omega_s^t$ along the time $t\in[0,T]$ for up to $M_s$ simplicial structural elements
$e_{s,k}^t\ (1\leq k\leq M_s)$ and $N_s$ structural nodes, where
$h_s$ can be different from $h$, and, ${\mathcal{T}}_{h}$ and
${\mathcal{T}}_{h_s}^t$ are nonconforming through the interface
$\Gamma^t$ for all $t\in[0,T]$, as illustrated in the right of Figure
\ref{fig:domain}. Therefore, material derivatives in these two
different domains have different expressions, as introduced before
and restated below,
\begin{equation}\label{Total-Derivatives}
\frac{D\bv}{Dt}=\left\{
\begin{array}{ll}
\frac{\partial\bv}{\partial
t}+\bv\cdot\nabla\bv, &\text{in }\Omega\times[0,T],\\
\frac{\partial\bv}{\partial t}=\frac{\partial^2\bX_s}{\partial
t^2}\circ\bX_s^{-1},& \text{in }\Omega_s^t\times[0,T].
\end{array} \right.
\end{equation}
To discretize the temporal derivatives in (\ref{Total-Derivatives}),
we introduce a uniform partition
$0=t_0<t_1<\cdot\cdot\cdot<t_N=T$ with the time-step size $\Delta
t=T/N$, and set $t^n=n\Delta t$, $\varphi^n=\varphi(\bx,t^n)$ for
$n=1,\cdots,N$.

\subsection{Reformulation of the deviatoric stress in update Lagrangian frame}
Since we now adopt the updated Lagrangian frame to
describe the structure motion, we shall update the deviatoric stress
term of the structure in (\ref{weakform2-1}) at the current
$(n+1)$-th time step, i.e.,
$\btau_s^{n+1}:=\mu_s(\bF^{n+1}(\bF^{n+1})^T-\bI)$, by utilizing the spatial gradient
$\nabla_{\bx^n}\bx^{n+1}$ that is computed on the known coordinate at
the previous $n$-th time step,
$\bx^n=\bX_{s}(\hat\bx,t^n)=\bX_s^n\in\Omega_s^n$, instead of using
$\nabla_{\hat\bx}\bx^{n+1}$ to directly compute $\bF^{n+1}$
which however belongs to the total Lagrangian description. To do so,
we employ the chain rule to reformulate $\btau_s^{n+1}$ through
$\bx^n$, leading to
\begin{eqnarray}
\btau_s^{n+1}&=&
\mu_s\left(\nabla_{\hat\bx}\bx^{n+1}(\nabla_{\hat\bx}\bx^{n+1})^T-\bI\right)\notag\\
&=&\mu_s\left(\nabla_{\bx^n}\bx^{n+1}(\nabla_{\bx^n}\bx^{n+1})^T-\bI\right)
+\mu_s\nabla_{\bx^n}\bx^{n+1}\left(\nabla_{\hat\bx}\bx^{n}(\nabla_{\hat\bx}\bx^{n})^T
-\bI\right)(\nabla_{\bx^n}\bx^{n+1})^T\notag\\
&=&\mu_s\left(\nabla_{\bx^n}\bx^{n+1}(\nabla_{\bx^n}\bx^{n+1})^T-\bI\right)
+ \nabla_{\bx^n}\bx^{n+1}\btau_s^n
(\nabla_{\bx^n}\bx^{n+1})^T,\label{tau_semi}
\end{eqnarray}
where $\btau_{s}^n$ is defined in $\Omega_s^n$ in the sense that
$\btau_{s}^n=\btau_{s}(\bx^n,t^n)$. It is easy to see $\btau_{s}^0=\bm{0}$.

Integrate (\ref{position-strong}) in time from $t^n$ to $t^{n+1}$,
yields
\begin{eqnarray}
\bX_{s}^{n+1}=\bX_{s}^{n}+\int_{t^n}^{t^{n+1}}\bv(\bX_s(\hat\bx,t),t)dt,
\quad \text{ or, }\quad
\bx^{n+1}
=\bx^{n}
+\int_{t^n}^{t^{n+1}}\bv(\bx(\hat\bx,t),t)dt.\label{position-update0}
\end{eqnarray}
Applying (\ref{position-update0}) to (\ref{tau_semi}), we can
further reformulate $\btau_s^{n+1}$ as follows,
\begin{equation}\label{tau_semi1}
\begin{array}{rl}
\btau_s^{n+1}
=&\mu_s\left(\left(\bI+\int_{t^n}^{t^{n+1}}\nabla_{\bx^n}\bv
dt\right) \left(\bI+\int_{t^n}^{t^{n+1}}\nabla_{\bx^n}\bv
dt\right)^T-\bI\right)+
\left(\bI+\int_{t^n}^{t^{n+1}}\nabla_{\bx^n}\bv dt\right)\btau_s^n
\left(\bI+\int_{t^n}^{t^{n+1}}\nabla_{\bx^n}\bv dt\right)^T\\
=&\mu_s\left(\int_{t^n}^{t^{n+1}}\nabla_{\bx^n}\bv dt+
\left(\int_{t^n}^{t^{n+1}}\nabla_{\bx^n}\bv dt\right)^T
+\int_{t^n}^{t^{n+1}}\nabla_{\bx^n}\bv dt
 \left(\int_{t^n}^{t^{n+1}}\nabla_{\bx^n}\bv dt\right)^T\right)
+ \btau_s^n +\\
&\int_{t^n}^{t^{n+1}}\nabla_{\bx^n}\bv dt\btau_s^n +
\btau_s^n\left(\int_{t^n}^{t^{n+1}}\nabla_{\bx^n}\bv dt\right)^T +
\int_{t^n}^{t^{n+1}}\nabla_{\bx^n}\bv dt\btau_s^n
\left(\int_{t^n}^{t^{n+1}}\nabla_{\bx^n}\bv dt\right)^T.
\end{array}
\end{equation}
Utilizing a proper quadrature rule to discretize the temporal
integrals in (\ref{tau_semi1}), e.g., simply picking up the right-endpoint
rule, we can obtain
\begin{equation}\label{tau_semi2}
\begin{array}{rcl}
\btau_s^{n+1}&=&\mu_s\Delta t\left(\nabla_{\bx^n}\bv^{n+1}+
(\nabla_{\bx^n}\bv^{n+1})^T +\Delta t\nabla_{\bx^n}\bv^{n+1}
 (\nabla_{\bx^n}\bv^{n+1})^T\right)
+ \btau_s^n +\\
&&\Delta t\nabla_{\bx^n}\bv^{n+1}\btau_s^n + \Delta
t\btau_s^n(\nabla_{\bx^n}\bv^{n+1})^T + (\Delta
t)^2\nabla_{\bx^n}\bv^{n+1}\btau_s^n (\nabla_{\bx^n}\bv^{n+1}
)^T+O((\Delta t)^2).
\end{array}
\end{equation}

\subsection{Full discretization of the UFMFD-FEM}
First, we introduce the following finite element spaces in which the velocity and pressure
are discretized, respectively,
\begin{equation}\label{FESpaces}
\begin{array}{rcl}
\bV_{h}&:=&\{\tilde\bv\in \bV \big| \tilde\bv|_{e_k}\in
P_m({e_k})^d,\forall {e_k}\in{\mathcal{T}}_{h}\},\quad
\bV_{0,h}\ :=\ \{\tilde\bv\in \bV_0 \big| \tilde\bv|_{e_k}\in
P_m({e_k})^d,\forall {e_k}\in{\mathcal{T}}_{h}\},\\
W_{h}&:=&\big\{\tilde p\in W \big| \tilde p|_{e_k}\in
P_1({e_k}),\forall {e_k}\in{\mathcal{T}}_{h}\big\},
\end{array}
\end{equation}
where $P_k$ denotes the $k$-th degree piecewise polynomial space. Here we
can either employ the Taylor-Hood ($P_2/P_1$) mixed finite element by letting $m=2$
in (\ref{FESpaces}), or the lowest equal-order mixed $P_1/P_1$ element by letting $m=1$
in (\ref{FESpaces}) with
the pressure stabilization scheme
\cite{Hughes.T;Franca.L;Balestra.M1986,Tezduyar.T1992,braack2007stabilized,Liu.X;Li.S2006a,Nemer;Larcher;Coupez;Hachem2021},
to approximate the saddle-point problem arising from the weak form
(\ref{weakform2-1}) and (\ref{weakform2-2}) in the finite element
spaces $\bV_{h}\times W_h\subset\bV\times W$.
The different choice of $m$ in (\ref{FESpaces}) will lead to
a large difference on the computational cost due to the quadratic element involved
in $P_2/P_1$ mixed element with at least twice as many degrees of freedom (DOFs)
of the velocity variable as the linear element in $P_1/P_1$ mixed element.
On the other hand, the $P_1/P_1$ element with pressure stabilization needs an
empirically well tuned stabilization parameter to make its approximation
stable and accurate. Their use on real world problems sometimes depends on the
practical interest. In fact, with a reasonably acceptable accuracy,
if the computational efficiency is more preferred in practice,
then the $P_1/P_1$ element with pressure stabilization scheme
is superior to the $P_2/P_1$ element.

Moreover, we choose the backward Euler scheme to
approximate temporal derivatives in (\ref{Total-Derivatives}), resulting in
the following fully discrete UFMFD-FEM for the FSI model presented in Section \ref{sec:FSImodel}.
We first interpolate the initial value function $\bv(\bx,0)$,
which is defined in (\ref{BICs-strong}), into $\bV_h$ to get the
discrete initial value $\bv_h^0=\prod_h \bv(\bx,0)\in\bV_h$, where
the interpolation operator, $\prod_h$, is described by the
definition of $\bV_h$ in (\ref{FESpaces}). Then, we define
the following fully discrete UFMFD-FEM: for
$n=0,1,\cdots,N-1$ and $\forall (\tilde\bv_h,\tilde p_h)\in \bV_{0,h} \times W_h$,
find $(\bv_{h}^{n+1},\ p_{h}^{n+1})\in
\bV_h\times W_h$ such that
\begin{eqnarray}
\left(\rho_f\frac{\bv_h^{n+1}-\bv_h^{n}}{\Delta
t},\tilde\bv_h\right)_{\Omega}+(\rho_f\bv_{h}^{n+1}\cdot\nabla\bv_{h}^{n+1},\tilde\bv_{h})_\Omega
+2\left(\mu_f\bD({\bv_h^{n+1}}),\bD(\tilde\bv_h)\right)_{\Omega}
-\left({p_h^{n+1}},\nabla
\cdot\tilde\bv_h\right)_{\Omega}-\left(\nabla \cdot{{\bv_h^{n+1}}},\tilde
p_h\right)_{\Omega}\notag\\
+\left((\rho_s-\rho_f)\frac{\bv_h^{n+1}\big|_{\Omega_s^{n}}-\bv_h^{n}\big|_{\Omega_s^{n}}}{\Delta
t},\tilde\bv_h\big|_{\Omega_s^{n}}\right)_{\Omega_s^{n}}
+\left(\btau_{s,h}^{n+1}
,\nabla_{\bx^n}\tilde\bv_h\big|_{\Omega_s^{n}}\right)_{\Omega_s^{n}}
-2\left(\mu_f\bD\big({\bv_h^{n+1}}\big|_{\Omega_s^{n}}\big),\bD\big(\tilde\bv_h\big|_{\Omega_s^{n}}\big)\right)_{\Omega_s^{n}}\notag\\
-\sum\limits_{k=1}^{M}\xi_{k}\lb\bR_h^{n+1},\nabla\tilde
p_h\rb_{e_k}
=(\rho_f\bg,\tilde\bv_h)_{\Omega}
+\left((\rho_s-\rho_f)\bg,\tilde\bv_h\big|_{\Omega_s^{n}}\right)_{\Omega_s^{n}}+\langle\bff_{f,N}^{n+1},\tilde\bv_h\rangle_{\partial\Omega_{f,N}^{n}}
-\langle f_{s,C_W}^{n+1},\tilde\bv_h\big|_{\Gamma_C^{n}}\cdot\bn\rangle_{\Gamma_C^{n}},\label{FDFEM}
\end{eqnarray}
where the discrete deviatoric stress
$\btau_{s,h}^{n+1}
$, which approximates to
$\btau_s^{n+1}$ defined in (\ref{tau_semi2}), is defined as
\begin{eqnarray}
\btau_{s,h}^{n+1}
&=&\mu_s\Delta
t\left(\nabla_{\bx^n}\bv_h^{n+1}\big|_{\Omega_s^{n}}
+\left(\nabla_{\bx^n}\bv_h^{n+1}\big|_{\Omega_s^{n}}\right)^T
+\Delta
t\nabla_{\bx^n}\bv_h^{n+1}\big|_{\Omega_s^{n}}\left(\nabla_{\bx^n}\bv_h^{n+1}\big|_{\Omega_s^{n}}\right)^T\right)
+ \btau_{s,h}^n +\notag\\
&&\Delta t\nabla_{\bx^n}\bv_h^{n+1}\big|_{\Omega_s^{n}}\btau_{s,h}^n
+ \Delta
t\btau_{s,h}^n\left(\nabla_{\bx^n}\bv_h^{n+1}\big|_{\Omega_s^{n}}\right)^T
+ (\Delta
t)^2\nabla_{\bx^n}\bv_h^{n+1}\big|_{\Omega_s^{n}}\btau_{s,h}^n\left(\nabla_{\bx^n}\bv_h^{n+1}\big|_{\Omega_s^{n}}\right)^T.\label{tau-dis}
\end{eqnarray}
Note that in (\ref{FDFEM}) all structural integrals are processed in
the known structural domain at the $n$-th time step, $\Omega_s^n$,
which is essential, thus the discrete deviatoric stress,
$\btau_{s,h}^{n+1}
$, is defined in $\Omega_s^{n}$ that makes an accurate sense because
all involved spatial gradients in (\ref{tau-dis}) are for the
restriction of $\bv_h^{n+1}$ in $\Omega_s^n$ with respect to
$\bx^n\in\Omega_s^n$, and, the involved $\btau_{s,h}^n$ belongs to
$\Omega_s^n$ as well in the sense that
$\btau_{s,h}^n=\btau_{s,h}(\bx^n,t^n)$. 
We will discuss in Section \ref{sec:postprocessing} about how to compute $\btau_{s,h}^n$ by means
of the finite element discretization in $\Omega_s^n$, for
$n=0,1,\cdots,N-1$.

If $m=2$ is chosen in (\ref{FESpaces}), then the last term on the left hand side of
(\ref{FDFEM}) vanishes, otherwise it remains with $m=1$ as
the pressure stabilization term defined in all (real and fictitious) fluidic elements,
$e_{k}\in{\mathcal{T}}_{h}\
(k=1,\cdots,M)$, elementwisely, where $\bR_h^{n+1}$ is
the discrete residual of the momentum equation of both the real
and fictitious fluid at the $(n+1)$-th time
step, which is defined below in terms of the unified finite element
solution $(\bv_h,p_h)$ that is only associated with
${\mathcal{T}}_{h}$,
\begin{eqnarray}
\bR_h^{n+1}&=&\rho_f\frac{\bv_h^{n+1}-\bv_h^{n}}{\Delta
t}+\rho_f\bv_{h}^{n+1}\cdot\nabla\bv_{h}^{n+1}
+\nabla p_h^{n+1}-\rho_f\bg,\label{residuals}
\end{eqnarray}
in which the viscous term (the second-order spatial derivative term)
vanishes from $\bR_h^{n+1}$ since $\bv_h\in\bV_h$ that is
discretized by piecewise linear ($P_1$) polynomials. $\xi_{k}\ (1\leq k\leq M)$ is
the pressure stabilization parameter that may affect
the computational results dearly. According to the usual choices
of the pressure stabilization parameter
\cite{Nemer;Larcher;Coupez;Hachem2021,Tezduyar.T1992,Franca;Hughes1993,Onate;Garcia1997},
we define the following elementwise stabilization
parameter $\xi_{k}$,
\begin{eqnarray}
\xi_{k}&=&\left[\left(\zeta_{0}\frac{\rho_f}{\Delta t}\right)^2+
\left(\zeta_{1}\frac{\mu_f}{h_k^2}\right)^2+
\left(\zeta_{2}\frac{\rho_f\|\bv_h^{n+1}\|}{h_k}\right)^2
\right]^{-\frac{1}{2}},\label{stab_parameter}
\end{eqnarray}
where $\zeta_{i}\ (i=0,1,2)$ are the
tunable, dimensionless parameters that are
independent of the element length, $h_k\ (1\leq k\leq M)$, and,
their exact values can only be determined empirically.

In addition, since both the trial function $\bv_h\in\bV_h$ and the
test function $\tilde\bv_h\in\bV_{0,h}$ in (\ref{FDFEM}) are
associated with the fixed Eulerian mesh ${\mathcal{T}}_{h}$ only,
to compute all structural terms in (\ref{FDFEM}),
we need to interpolate both $\bv_h$ and
$\tilde\bv_h$ into the updated Lagrangian mesh
${\mathcal{T}}_{h_s}^n$ in $\Omega_s^n$, which are represented by
$\bphi_h|_{\Omega_s^n}$ in (\ref{FDFEM}) and (\ref{tau-dis}),
generally. A detailed interpolation procedure for how to implement
$\bphi_h|_{\Omega_s^n}$ from ${\mathcal{T}}_{h}$ to ${\mathcal{T}}_{h_s}^n$
is described in \ref{appendix:interpolation}.

\begin{remark}\label{rmk:theory}
Although a rigorous error analysis for the proposed UFMFD-FEM (\ref{FDFEM})-(\ref{stab_parameter})
is not yet studied in this paper (which will be considered in our next work),
we suppose that the proposed method may own a similar convergence property
with that of the distributed Lagrange multiplier/fictitious domain (DLM/FD)-finite element method
for various interface problems \cite{Boffi2016,Wang2016,Lundberg;Sun2016,Sun2016,Sun2017},
i.e., when $h_s/h \approx 1$, the spatial convergence rate of our developed UFMFD-FEM
with a pressure-stabilized $P_1/P_1$ mixed element for the velocity
and pressure on a quasi-uniform grid supposes to be
\begin{eqnarray}
\|\bv-\bv_h\|_{L^2(\Omega)^d}=O(h^r),\ \|p-p_h\|_{L^2(\Omega)}=O(h^q),
\hbox{ for } \bv\in H^r(\Omega)^d\ (1< r \leq 3/2) \hbox{ and }
p\in H^q(\Omega)\ (0 < q \leq 1/2),
\label{errorestimate}
\end{eqnarray}
while its temporal convergence rate is of the first order, i.e., $O(\Delta t)$
due to the employment of backward Euler scheme on the temporal discretization.
We will attempt to prove (\ref{errorestimate}) in our next paper
based upon the Babu\v{s}ka--Brezzi's theory
\cite{Babuska1971,Brezzi_intr} by proving a corresponding
\emph{inf-sup} condition of a properly defined total bilinear form
associated with the saddle-point problem that is involved in the developed
UFMFD-FEM, as we do in \cite{Wang2016,Lundberg;Sun2016,Sun2016,Sun2017}.
In this paper, we focus more on the development of advanced numerical
discretizations and algorithms, as well as their numerical validations
through benchmark problems, therefore instead of proving (\ref{errorestimate})
theoretically, we conduct a series of numerical
convergence tests to validate (\ref{errorestimate}), as shown in Section \ref{sec:numericalexperiment}.
\end{remark}

\subsection{Post-processing the structural information}\label{sec:postprocessing}
After obtaining $\bv_h^{n+1}$ from (\ref{FDFEM}), we then use,
\begin{eqnarray}
\bX_{s,h_s}^{n+1}&=&\bX_{s,h_s}^{n}+\Delta
t\bv_h^{n+1}\big|_{\Omega_s^{n}},\label{position-update}
\end{eqnarray}
at the $(n+1)$-th time step to attain the position of the immersed
structure domain, $\Omega_s^{n+1}$, or, the updated Lagrangian mesh
of the structure, ${\mathcal{T}}_{h_s}^{n+1}$, for
$n=0,1,\cdots,N-1$, where $\bX_{s,h_s}^{0}$ denotes the initial
structural position that is also discretized by ${\mathcal{T}}_{h_s}^{0}$,
on which we set $\btau_{s,h}^0=\bm{0}$.

Instead of a direct interpolation, we can employ the
$L^2$-projection method to compute $\bv_h^{n+1}\big|_{\Omega_s^{n}}$
shown in (\ref{position-update}) within the structure domain
$\Omega_s^n$ in a more accurate fashion. To do so, we first define
the following finite element spaces in $\Omega_s^t\times[0,T]$,
\begin{equation}\label{FESpaces-structure}
\begin{array}{rcl}
\bU_{h_s}^0&:=&\{\bphi\in L^2(\Omega_s^0)^d \big|
\bphi|_{e_{s,k}^{0}}\in P_1({e_{s,k}^{0}})^d,\forall
{e_{s,k}^{0}}\in{\mathcal{T}}_{h_s}^0\},\\
\bU_{h_s}^t&:=&\{\bphi:\Omega_s^t\times[0,T]\rightarrow
\mathbb{R}^d\big| \bphi=\hat{\bphi}\circ\bX^{-1}_{s,h_s}(t),\forall
\hat{\bphi}\in\bU_{h_s}^0\},
\end{array}
\end{equation}
where $\bX_{s,h_s}(t):\Omega_s^0\rightarrow\Omega_s^t$ is a discrete
Lagrangian mapping approximated by $P_1^d$ Lagrange-type finite
elements such that $\bx(\hat\bx,t)=\bX_{s,h_s}(\hat\bx,t)$ for any
$\bx\in\mathcal{T}_{h_s}^t$ that corresponds to its initial position
$\hat\bx\in\mathcal{T}_{h_s}^0$. And, due to the following long-standing
Proposition \ref{prop1} (see, e.g., \cite{Grisvard1985}), we
know $\bU_{h_s}^t\subset L^2(\Omega_s^t)^d$.
\begin{prop}\label{prop1}
Consider $\Omega_1$ and $\Omega_2$ are two bounded open subsets of
$\mathbb{R}^d$, and assume $\bX\in W^{1,\infty}(\Omega_1)$.
Suppose also that $\bX: \Omega_1\rightarrow\Omega_2$ is invertible and
such that $\bX^{-1}\in W^{1,\infty}(\Omega_2)$. Then for any $u\in
H^1(\Omega_2)$ we have $u\circ \bX\in H^1(\Omega_1)$.
\end{prop}

Now we can define the following $L^2$-projection for
$\bv_h^{n+1}\big|_{\Omega_s^{n}}$ by the finite element discretization:
find $\bphi_{h_s}^n\in\bU_{h_s}^n$ such that
\begin{eqnarray}
\left(\bphi_{h_s}^n,\tilde
\bphi\right)_{\Omega_s^{n}}&=&\sum\limits_{k=1}^{M_s}
\left(\bv_h^{n+1}\big|_{\Omega_s^{n}},
\tilde\bphi\right)_{e_{s,k}^n},\quad \forall \tilde\bphi\in
\bU^n_{h_s},\qquad n=0,1,\cdots,N-1.\label{L2-projection}
\end{eqnarray}




Then, (\ref{position-update}) can be changed to
\begin{eqnarray}
\bX_{s,h_s}^{n+1}&=&\bX_{s,h_s}^{n}+\Delta t\bphi_{h_s}^n, \quad
n=0,1,\cdots,N-1.\label{position-update1}
\end{eqnarray}
The new position of $\Omega_s^{n+1}$ and its triangulation,
$\mathcal{T}_{h_s}^{n+1}$, are thus updated.

Next, we can also adopt the $L^2$ projection to update
$\btau_{s,h}^{n+1}$ in ${\Omega_s^{n+1}}$ by using the newly
obtained $\bv_h^{n+1}$, $\bX_{s,h_s}^{n+1}$ and
$\mathcal{T}_{h_s}^{n+1}$. Define the corresponding finite element
spaces first,
\begin{equation}\label{FESpaces-tau-structure}
\begin{array}{rcl}
\bm{\Upsilon}_{h_s}^0&:=&\{\bphi\in
L^2(\Omega_s^0)^{\frac{d(d+1)}{2}} \big| \bphi|_{e_{s,k}^{0}}\in
P_1({e_{s,k}^{0}})^{\frac{d(d+1)}{2}},\forall
{e_{s,k}^{0}}\in{\mathcal{T}}_{h_s}^0\},\\
\bm{\Upsilon}_{h_s}^t&:=&\{\bphi:\Omega_s^t\times[0,T]\rightarrow
\mathbb{R}^{\frac{d(d+1)}{2}}\big|
\bphi=\hat{\bphi}\circ\bX^{-1}_{s,h_s}(t),\forall
\hat{\bphi}\in\bm{\Upsilon}_{h_s}^0\},
\end{array}
\end{equation}
which leads to $\bm{\Upsilon}_{h_s}^t\subset
L^2(\Omega_s^t)^{\frac{d(d+1)}{2}}$.
On the other hand, from (\ref{tau_semi2}) we see
$\nabla_{\bx^n}\bv^{n+1}$ needs to be approximated using the
obtained $\bv_h^{n+1}$ in the updated ${\Omega_s^{n+1}}$, i.e., we
need to frequently compute its discrete version,
$\nabla_{\bx^n}\bv_h^{n+1}\big|_{\Omega_s^{n+1}}$, in order to
finally update $\btau_{s,h}^{n+1}$ in a discrete manner. However, it
is not a direct computation because
$\bv_h^{n+1}\big|_{\Omega_s^{n+1}}$, which is defined in
${\Omega_s^{n+1}}$, takes a spatial gradient with respect to
$\bx^n\in \Omega_s^{n}$. Actually, by the chain rule we can rewrite
$\nabla_{\bx^n}\bv_h^{n+1}\big|_{\Omega_s^{n+1}}=\bG^{n+1}\bF^{n+1,n}$,
where
$\bG^{n+1}:=\nabla_{\bx^{n+1}}\bv_h^{n+1}\big|_{\Omega_s^{n+1}}$
that is defined in $\Omega_s^{n+1}$, and
$\bF^{n+1,n}:=\nabla_{\bx^{n}}\bx^{n+1}$ that can be reformulated as
$\bF^{n+1,n}=\bI+\Delta t\bQ^{n}$, where
$\bQ^{n}:=\nabla_{\bx^{n}}\bv_h^{n+1}\big|_{\Omega_s^{n}}$ due to
(\ref{position-update}). Thus $\bF^{n+1,n}$ is defined in
$\Omega_s^{n}$, essentially.
Before defining the $L^2$-projection for
$\btau_{s,h}^{n+1}\in\bm{\Upsilon}_{h_s}^{n+1}$, we first adopt the
$L^2$-projection method to project $\bG^{n+1}$ onto
$\bm{\Upsilon}_{h_s}^{n+1}$, i.e., find
$\bG_{P}^{n+1}\in\bm{\Upsilon}_{h_s}^{n+1}$ such that
\begin{eqnarray}
\left(\bG_{P}^{n+1},\tilde
\bphi\right)_{\Omega_s^{n+1}}&=&\sum\limits_{k=1}^{M_s}
\left(\bG^{n+1}, \tilde\bphi\right)_{e_{s,k}^{n+1}},\quad \forall
\tilde\bphi\in \bm{\Upsilon}_{h_s}^{n+1},\qquad
n=0,1,\cdots,N-1.\label{L2-projection-G}
\end{eqnarray}
Now based on (\ref{tau_semi2}), we are able to define the following
$L^2$-projection for $\btau_{s,h}^{n+1}$ in the finite element approximation:
find $\btau_{s,h}^{n+1}\in\bm{\Upsilon}^{n+1}_{h_s}$ such that
\begin{eqnarray}
&&\left(\btau_{s,h}^{n+1},\tilde \bphi\right)_{\Omega_s^{n+1}}
=\mu_s\Delta
t\sum\limits_{k=1}^{M_s}\left(\bG_{P}^{n+1}\bF^{n+1,n}+(\bG_{P}^{n+1}\bF^{n+1,n})^T
+\Delta t\bG_{P}^{n+1}\bF^{n+1,n}(\bG_{P}^{n+1}\bF^{n+1,n})^T,\tilde\bphi\right)_{e_{s,k}^{n+1}}\notag\\
&&+\Delta
t\sum\limits_{k=1}^{M_s}\left(\bG_{P}^{n+1}\bF^{n+1,n}\btau_{s,h}^n+\btau_{s,h}^n(\bG_{P}^{n+1}\bF^{n+1,n})^T
+\Delta t\bG_{P}^{n+1}\bF^{n+1,n}\btau_{s,h}^n(\bG_{P}^{n+1}\bF^{n+1,n})^T,\tilde\bphi\right)_{e_{s,k}^{n+1}}\notag\\
&&+\sum\limits_{k=1}^{M_s}\left(\btau_{s,h}^n,\tilde\bphi\right)_{e_{s,k}^{n+1}}=\sum\limits_{i=1}^7R_i,
\qquad\qquad \forall \tilde\bphi\in \bm{\Upsilon}^{n+1}_{h_s},\qquad
n=0,1,\cdots,N-1,\label{L2-projection-tau}
\end{eqnarray}
where all seven terms on the right hand side of
(\ref{L2-projection-tau}), which are to be integrated in
$\Omega_s^{n+1}$, involve integrant functions defined in
$\Omega_s^n$, such as $\btau_{s,h}^n$ and $\bQ^{n}$ included in $\bF^{n+1,n}$
which are actually represented by a composite function
with respect to $\bx^{n+1}$ through the discrete Lagrangian mapping
$\bX_{s,h_s}(t)$, i.e.,
$\btau_{s,h}^n\circ\bX_{s,h_s}^{n}(\bX_{s,h_s}^{n+1})^{-1}(\bx^{n+1})$ and
$\bQ^{n}\circ\bX_{s,h_s}^{n}(\bX_{s,h_s}^{n+1})^{-1}(\bx^{n+1})$.
Thus, in order to carry out all integrals on the right hand side of
(\ref{L2-projection-tau}), it is necessary to change variables from
$\Omega_s^{n+1}$ to $\Omega_s^n$ through the following discrete
mapping, $\bX_{s,h_s}^{m,n}:\Omega_s^m\rightarrow\Omega_s^n$ such
that
\begin{eqnarray}\label{transfer_mn}
\bX_{s,h_s}^{m,n}:=\bX_{s,h_s}^{n}(\bX_{s,h_s}^{m})^{-1},
\end{eqnarray}
by introducing the Jacobian
$J^{n+1,n}=\operatorname{det}(\bF^{n+1,n})$. For example, the most
complicated term on the right hand side of
(\ref{L2-projection-tau}), i.e., $R_6$, can actually be rewritten as
follows in an exact fashion,
\begin{eqnarray}
R_6=(\Delta
t)^2\sum\limits_{k=1}^{M_s}\left(J^{n+1,n}\bG_{P}^{n+1}\circ\bX_{s,h_s}^{n,n+1}\big(\bI+\Delta
t\bQ^{n}\big)\btau_{s,h}^n\big(\bI+\Delta
t(\bQ^{n})^T\big)(\bG_{P}^{n+1})^T\circ\bX_{s,h_s}^{n,n+1},\tilde\bphi\circ\bX_{s,h_s}^{n,n+1}\right)_{e_{s,k}^{n}},
\end{eqnarray}
where $\bG_{P}^{n+1}\circ\bX_{s,h_s}^{n,n+1}$ and
$\tilde\bphi\circ\bX_{s,h_s}^{n,n+1}$ are defined in $\Omega_s^n$ through
the discrete mapping $\bX_{s,h_s}^{n,n+1}$, although both $\bG_{P}^{n+1}$
and $\tilde\bphi$ belong to
$\bm{\Upsilon}^{n+1}_{h_s}\ (n=0,1,\cdots,N-1)$ which is constructed
using the same finite elements as for $\bm{\Upsilon}^{0}_{h_s}$
based on the Lagrangian mapping $\bX_{s,h_s}(t)$ all the
time. Other less complicated terms on the right hand side of
(\ref{L2-projection-tau}) can also be similarly transformed to the integral
form in $\Omega_s^n$ as well.

\section{Numerical algorithms}
It is obvious that (\ref{FDFEM}) is highly nonlinear due to the
nonlinear convection term of the fluid, the updated
Lagrangian-based deviatoric stress of the structure
$\btau_{s,h}^{n+1}$ as shown in (\ref{tau-dis}), the unknown contact force
$f_{s,C_w}$ that depends on the normal structural stress
over the unknown contacting surface $\Gamma_C^t$, as well as the updating
structural domain that solely depends on the structural velocity all the time.
In the following, we describe in Algorithm \ref{algorithm1} how the
developed UFMFD-FEM is implemented for the presented FSCI model by means of
the Newton's linearization for handling nonlinear fluid convection and
structural stress terms, and by means of the fixed-point iteration for updating
the structural domain and computing the contact force when exists.
\begin{alg}\label{algorithm1} The overall algorithm of implementing UFMFD-FEM for FSCI problems.
\begin{enumerate}
\item Initialization of the time marching. Set the time step $n=0$ and let
$\bv_h^0=\prod_h \bv(\bx,0)\in\bV_h$, $p_h^0=0$, $\btau_{s,h}^0=\bm{0}$. $\Omega$ and $\Omega_s^0$
are triangulated as $\mathcal{T}_{h}$ and $\mathcal{T}_{h_s}^{0}$, respectively.
\item Call Algorithm \ref{algorithm3} to carry out the fixed-point iteration for
the structural collision determination and Newton's linearization iteration for
solving the entire FSCI system at the $(n+1)$-th time step $(n\geq 0)$.
\item  Update the discrete deviatoric stress of the structure,
$\btau_{s,h}^{n+1}$, by means of (\ref{L2-projection-G}) and (\ref{L2-projection-tau}).
\item  Time marching. If $n+1<N$, then set $n+1$ to $n$, go back to Step 2 and continue the time marching. Otherwise, stop the entire computation.
\end{enumerate}
\end{alg}

Next, we carry out the following structural collision algorithm to determine the
repulsive contact force $f_{s,C_W}^{n+1}$ and the contacting surface $\Gamma_C^{n+1}$
between the structure and fluidic channel wall at the $(n+1)$-th time step.
\begin{alg}\label{algorithm3} The nonlinear iteration algorithm for
structural collision determination at the $(n+1)$-th time step.
    \begin{enumerate}
        \item Initialization of the fixed-point iteration
        for collision determination.
        Set $k=0$, let $(\bv_h^{n+1,0},p_h^{n+1,0})=(\bv_h^{n},p_h^{n})$ and $f_{s,C_W}^{n+1,0}=0$ be the initial
        guess of the fixed-point iteration at the $(n+1)$-th time step $(n\geq 0)$.
        Let $\varepsilon$ be a given tolerance.
        \item Call Algorithm \ref{algorithm1-1} to carry out the Newton's linearization
        for solving the FSCI system at the $(n+1)$-th time step and the $(k+1)$-th
        fixed-point iteration step, and obtain the desired numerical solution
        $(\bv_h^{n+1,k+1},p_h^{n+1,k+1})$,
        including the structural velocity
        $\bv_{h_s}^{n+1,k+1}\circ\bX_{s,h_s}^{n,n+1}=\bv_{h}^{n+1,k+1}\big|_{\Omega_s^n}$ by (\ref{L2-projection})
        and the structural mesh $\mathcal{T}_{h_s}^{n+1,k+1}$ by (\ref{position-update1}).
        \item To determine $\Gamma_C^{n+1,k+1}$, 
we first introduce a fluidic boundary layer whose boundary is $\partial \Omega\cup \Gamma_h^{BL}$
and whose triangulation is $\mathcal{T}_h^{BL}\subset \mathcal{T}_h$
along the fluidic channel wall $\partial \Omega$, as illustrated in Figure \ref{fig:boundarylayer},
where $\mathcal{T}_h^{BL}$ might be
the first two-layer fluidic elements within the boundary layer
that is attached to $\partial \Omega$ and
bounded by $\Gamma_h^{BL}$.
Then, let $\mathcal{T}_{h_s}^{n+1,k+1}(\partial \Omega_s)$ be a partition of $\partial \Omega_s^{n+1,k+1}$ associated with $\mathcal{T}_{h_s}^{n+1,k+1}$,
and define
\begin{align}\label{eq_61}
    \Gamma_{C}^{n+1,k+1} &= \left\{e\in\mathcal{T}_{h_s}^{n+1,k+1}(\partial \Omega_s): e\cap \mathcal{T}_h^{BL}\neq \emptyset \right\},
\end{align}
as sketched by red line segments in Figure \ref{fig:boundarylayer}
over the structural surface.
\begin{figure}[hbt]
    \centering
    \includegraphics[height=5.5cm]{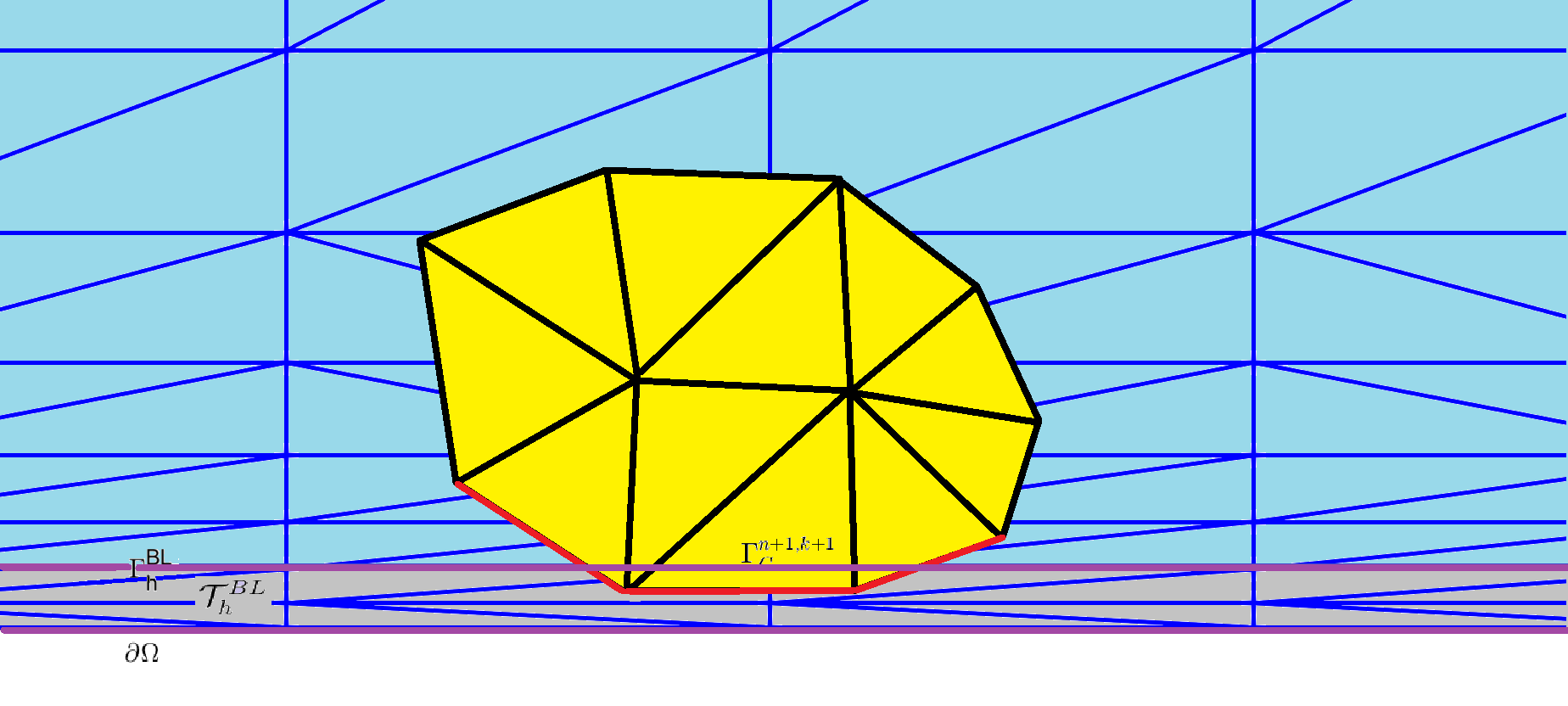}
    \caption{A schematic background fluidic and foreground structural mesh portion near the fluidic
    channel wall (the bottom purple line) and the approaching structure (the yellow area),
    which indicates the boundary-layer mesh $\mathcal{T}_h^{BL}$
        (the grey area) that is overlapped by the structural mesh, and
        the contacting structural surface elements $\Gamma_C^{n+1,k+1}$ that are
        marked in red line segments.}\label{fig:contactfig4}\label{fig:boundarylayer}
\end{figure}

        \item If  $\Gamma_{C}^{n+1,k+1} =\emptyset$, or
        $\bv_{h_s}^{n+1,k+1}\cdot \bn\big|_{\Gamma_{C}^{n+1,k+1}}<\varepsilon$,
        then stop the fixed-point iteration, take all obtained numerical solutions
        as the desired ones at the $(n+1)$-th time step, and exit the collision algorithm.
        Otherwise, go to Step 5.

        \item Determine $f_{s,C_W}^{n+1,k+1}$ by
        \begin{align}\label{contact_force_num}
            f_{s,C_W}^{n+1,k+1} =
            \max\bigg\{\zeta \frac{m_s}{\Delta t\text{A}(\Gamma_{C}^{n+1,k+1})}\bv_{h_s}^{n+1,k+1}\cdot \bn
            +f_{s,C_W}^{n+1,k},0\bigg\},
            & \text{ on } \Gamma_{C}^{n+1,k+1},
        \end{align}
        where $f_{s,C_W}^{n+1,k+1}$ is defined as a piecewise constant
        function on $\Gamma_{C}^{n+1,k+1}$,
        $\zeta$ is a tunable positive dimensionless constant, $m_s$ denotes the mass of structure,
        and A($\Gamma_{C}^{n+1,k+1}$) the measure of contacting surface $\Gamma_{C}^{n+1,k+1}$.

        \item Set $k+1$ to $k$, go back to Step 2 with the updated $\Gamma_{C}^{n+1,k}$
        and $f_{s,C_W}^{n+1,k}$ in (\ref{FDFEM-linear})
        to continue the fixed-point iteration.

    \end{enumerate}
\end{alg}
\begin{remark}
In Algorithm \ref{algorithm3}, when the structural boundary node
enters $\mathcal{T}_h^{BL}$,
a nonlinear iteration is triggered to prevent the structural boundary node
from penetrating the fluidic boundary layer as well as the fluidic channel wall.
In the process of this nonlinear iteration, the magnitude of
repulsive contact force increases until it is large enough to
make the velocity at this boundary node in the outward normal direction
is very small or negative. We point out that \eqref{contact_force_num} and the stopping criteria in Step 4 of this algorithm
can be viewed as a numerical implementation of the complemental contact condition \eqref{struct-bc-2a}.
\end{remark}

When calling Algorithm \ref{algorithm3} in Step 2 of Algorithm \ref{algorithm1},
we need to iteratively call Algorithm \ref{algorithm1-1} in Step 2 of Algorithm \ref{algorithm3}
to solve (\ref{FDFEM}) by virtue of Newton's linearization, as shown below.
\begin{alg}\label{algorithm1-1} Newton's linearization algorithm for
UFMFD-FEM at the $(n+1)$-th time step
and the $(k+1)$-th fixed-point iteration step.
\begin{enumerate}
\item Initialization of the nonlinear iteration.
Let $(\bv_h^{n+1,k+1,0},p_h^{n+1,k+1,0})=(\bv_h^{n+1,k},p_h^{n+1,k})$ be the initial
guess at the $(n+1)$-th time step $(n\geq 0)$ and the $(k+1)$-th fixed-point iteration step $(k\geq 0)$.
If $k=0$, let $(\bv_h^{n+1,0},p_h^{n+1,0})=(\bv_h^{n},p_h^{n})$.
\item Solve the linearized UFMFD-FEM at the $(m+1)$-th iteration step for
$(\bv_h^{n+1,k+1,m+1},p_h^{n+1,k+1,m+1})\in \bV_h\times W_h$ such that
\begin{numcases}{}
\left(\rho_f\frac{\bv_h^{n+1,k+1,m+1}-\bv_h^{n}}{\Delta
t},\tilde\bv_h\right)_{\Omega}+(\rho_f\bv_{h}^{n+1,k+1,m}\cdot\nabla\bv_{h}^{n+1,k+1,m+1},\tilde\bv_{h})_\Omega
+(\rho_f\bv_{h}^{n+1,k+1,m+1}\cdot\nabla\bv_{h}^{n+1,k+1,m},\tilde\bv_{h})_\Omega\notag\\
+2\left(\mu_f\bD({\bv_h^{n+1,k+1,m+1}}),\bD(\tilde\bv_h)\right)_{\Omega}
-\left({p_h^{n+1,k+1,m+1}},\nabla
\cdot\tilde\bv_h\right)_{\Omega}-\left(\nabla
\cdot{{\bv_h^{n+1,k+1,m+1}}},\tilde p_h\right)_{\Omega}\notag\\
+\left((\rho_s-\rho_f)\frac{\bv_h^{n+1,k+1,m+1}\big|_{\Omega_s^{n}}-\bv_h^{n}\big|_{\Omega_s^{n}}}{\Delta
t},\tilde\bv_h\big|_{\Omega_s^{n}}\right)_{\Omega_s^{n}}
+\left(\bar\btau_{s,h}^{n+1,k+1,m+1}
,\nabla_{\bx^n}\tilde\bv_h\big|_{\Omega_s^{n}}\right)_{\Omega_s^{n}}\notag\\
-2\left(\mu_f\bD\big({\bv_h^{n+1,k+1,m+1}}\big|_{\Omega_s^{n}}\big),\bD\big(\tilde\bv_h\big|_{\Omega_s^{n}}\big)\right)_{\Omega_s^{n}}
-\sum\limits_{j=1}^{M}\xi_{j}\lb\bR_h^{n+1,k+1,m+1},\nabla\tilde
p_h\rb_{e_j}\notag\\
=(\rho_f\bv_{h}^{n+1,k+1,m}\cdot\nabla\bv_{h}^{n+1,k+1,m},\tilde\bv_{h})_\Omega
+\mu_s(\Delta
t)^2\left(\nabla_{\bx^n}\bv_h^{n+1,k+1,m}\big|_{\Omega_s^{n}}\left(\nabla_{\bx^n}\bv_h^{n+1,k+1,m}\big|_{\Omega_s^{n}}\right)^T
,\nabla_{\bx^n}\tilde\bv_h\big|_{\Omega_s^{n}}\right)_{\Omega_s^{n}}\notag\\
\quad-\left(\btau_{s,h}^{n}
,\nabla_{\bx^n}\tilde\bv_h\big|_{\Omega_s^{n}}\right)_{\Omega_s^{n}}
+(\Delta
t)^2\left(\nabla_{\bx^n}\bv_h^{n+1,k+1,m}\big|_{\Omega_s^{n}}\btau_{s,h}^n\left(\nabla_{\bx^n}\bv_h^{n+1,k+1,m}\big|_{\Omega_s^{n}}\right)^T
,\nabla_{\bx^n}\tilde\bv_h\big|_{\Omega_s^{n}}\right)_{\Omega_s^{n}}\notag\\
\quad+(\rho_f\bg,\tilde\bv_h)_{\Omega}
+\left((\rho_s-\rho_f)\bg,\tilde\bv_h\big|_{\Omega_s^{n}}\right)_{\Omega_s^{n}}+\langle\bff_{f,N}^{n+1},\tilde\bv_h\rangle_{\partial\Omega_{f,N}^{n+1}}
-\langle f_{s,C_W}^{n+1,k},\tilde\bv_h\big|_{\Gamma_C^{n+1,k}}\cdot\bn\rangle_{\Gamma_C^{n+1,k}},\label{FDFEM-linear}\\
\hspace{10cm}\forall (\tilde\bv_h,\tilde p_h)\in \bV_{0,h} \times W_h,\quad m=0,1,2, \cdots,\notag
\end{numcases}
where
\begin{eqnarray}
&&\bar\btau_{s,h}^{n+1,k+1,m+1}
=\mu_s\Delta
t\left(\nabla_{\bx^n}\bv_h^{n+1,k+1,m+1}\big|_{\Omega_s^{n}}
+\left(\nabla_{\bx^n}\bv_h^{n+1,k+1,m+1}\big|_{\Omega_s^{n}}\right)^T\right.\notag\\
&&\left.+\Delta
t\nabla_{\bx^n}\bv_h^{n+1,k+1,m}\big|_{\Omega_s^{n}}\left(\nabla_{\bx^n}\bv_h^{n+1,k+1,m+1}\big|_{\Omega_s^{n}}\right)^T+\Delta
t\nabla_{\bx^n}\bv_h^{n+1,k+1,m+1}\big|_{\Omega_s^{n}}\left(\nabla_{\bx^n}\bv_h^{n+1,k+1,m}\big|_{\Omega_s^{n}}\right)^T\right)\notag\\
&&+ \Delta
t\nabla_{\bx^n}\bv_h^{n+1,k+1,m+1}\big|_{\Omega_s^{n}}\btau_{s,h}^n +
\Delta
t\btau_{s,h}^n\left(\nabla_{\bx^n}\bv_h^{n+1,k+1,m+1}\big|_{\Omega_s^{n}}\right)^T
\notag\\
&&+(\Delta
t)^2\nabla_{\bx^n}\bv_h^{n+1,k+1,m}\big|_{\Omega_s^{n}}\btau_{s,h}^n\left(\nabla_{\bx^n}\bv_h^{n+1,k+1,m+1}\big|_{\Omega_s^{n}}\right)^T\notag\\
&&+(\Delta
t)^2\nabla_{\bx^n}\bv_h^{n+1,k+1,m+1}\big|_{\Omega_s^{n}}\btau_{s,h}^n\left(\nabla_{\bx^n}\bv_h^{n+1,k+1,m}\big|_{\Omega_s^{n}}\right)^T,\label{tau-dis-iterate}\\
&&\bR_h^{n+1,k+1,m+1}=\rho_f\frac{\bv_h^{n+1,k+1,m+1}-\bv_h^{n}}{\Delta
t}+\rho_f\bv_{h}^{n+1,k+1,m}\cdot\nabla\bv_{h}^{n+1,k+1,m+1} +\nabla
p_h^{n+1,k+1,m+1}-\rho_f\bg,\label{residuals1}\\
&&\xi_{j}=\left[\left(\zeta_{0}\frac{\rho_f}{\Delta t}\right)^2+
\left(\zeta_{1}\frac{\mu_f}{h_j^2}\right)^2+
\left(\zeta_{2}\frac{\rho_f\|\bv_h^{n+1,k+1,m}\|}{h_j}\right)^2
\right]^{-\frac{1}{2}},\quad 1\leq j\leq M.\label{stab_parameter1}
\end{eqnarray}
\item Check the stopping criteria for the nonlinear iteration. Stop the iteration if \be
\|\bv_h^{n+1,k+1,m+1}-\bv_h^{n+1,k+1,m}\|_{0}+\|p_h^{n+1,k+1,m+1}-p_h^{n+1,k+1,m}\|_{0}
\leq\varepsilon\nonumber, \ee where $\varepsilon$ is a given
tolerance, set
$(\bv_h^{n+1,k+1},p_h^{n+1,k+1})=(\bv_h^{n+1,k+1,m+1},p_h^{n+1,k+1,m+1})$,
update the structural domain $\Omega_s^{n+1,k+1}$ and
its updated Lagrangian mesh $\mathcal{T}_{h_s}^{n+1,k+1}$ by
(\ref{L2-projection}) and (\ref{position-update1}),
and stop the nonlinear iteration. Otherwise, set $m+1$ to $m$, go back to Step 2 and continue
the nonlinear iteration.
\end{enumerate}
\end{alg}

\section{Numerical experiments}\label{sec:numericalexperiment}
In this section, we validate numerical performances of the developed UFMFD-FEM via two
benchmark FSI problems first: the particulate fluid of (1) the simple shear flow type;
and (2) the plane Poiseuille flow type, then via a self-defined small-scale
DLD problem second, finally via
a realistic large-scale DLD problem by comparing with its physical experiments,
where the particle collision with tilted pillar arrays inside the microfluidic channel are engaged
in the fluid flow.

\subsection{Example 1: The simple shear flow-particle interaction benchmark problem}\label{sec:benchmark1}
First, we consider the case of a neutrally buoyant particle of circular
shape locating at the middle of the fluidic channel between two walls.
This example was considered by Pan et al. in \cite{Pan;Glowinski2015},
where the distributed Lagrange multiplier/fictitious domain (DLM/FD) method was employed
to solve the problem. Here we use the same sets of parameters as in \cite{Pan;Glowinski2015}
to validate our proposed UFMFD-FEM that produces very similar results with those in \cite{Pan;Glowinski2015},
as shown below.

As illustrated in the left part of Figure \ref{fig:benchmark1}, the computational domain
$\Omega: =[0,2]\times[0,2]$, and the mass center of the circular particle
locates at $(1,1)$. The periodic velocity boundary condition is defined on
the left and right boundary, $\Gamma_P^f$, and the shear flow condition
is assigned to the top boundary as $\bv_f|_{\Gamma_{D,1}^f}=(\frac{\gamma H}{2},0)=(1,0)^T$
and to the bottom boundary as $\bv_f|_{\Gamma_{D,2}^f}=(-\frac{\gamma H}{2},0)=(-1,0)^T$, respectively,
where $\gamma$ denotes the shear rate of the fluid with the fixed value $\gamma=1$ (s$^{-1}$),
and $H=2$ is the distance between the top and bottom walls.
\begin{figure}[htb]
    \centering
    \includegraphics[height=4cm]{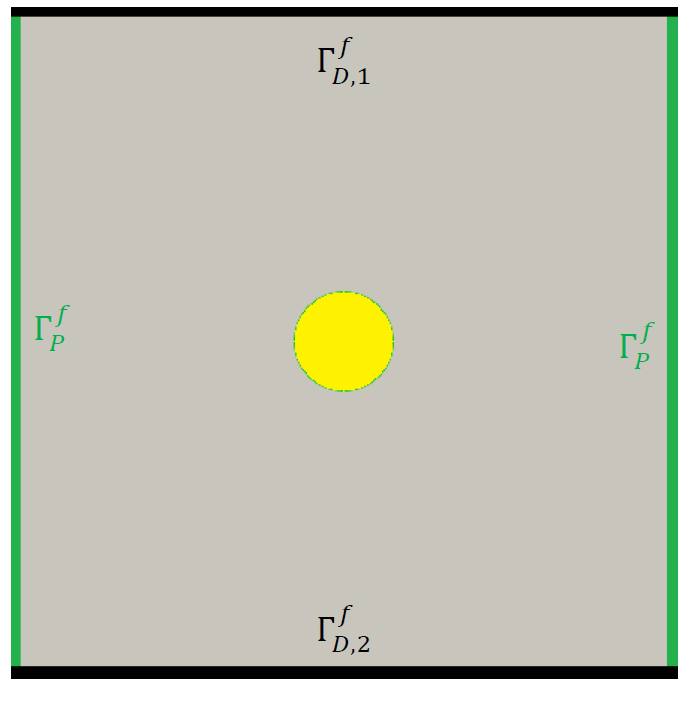}
    \hspace{0.2in}
    \includegraphics[height=4cm]{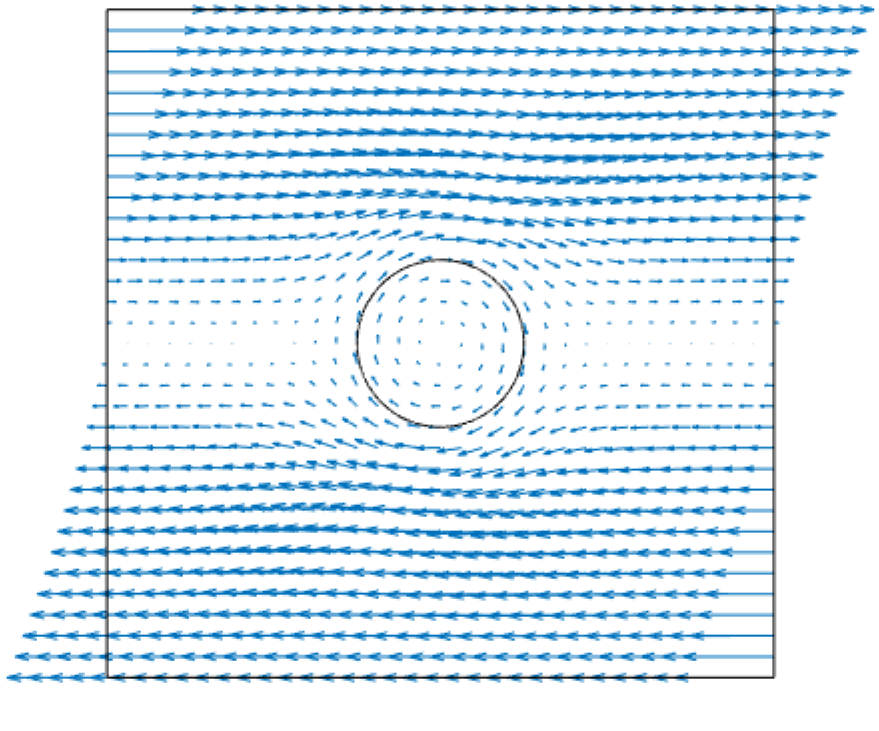}
    \hspace{0.0in}
    \includegraphics[height=4cm]{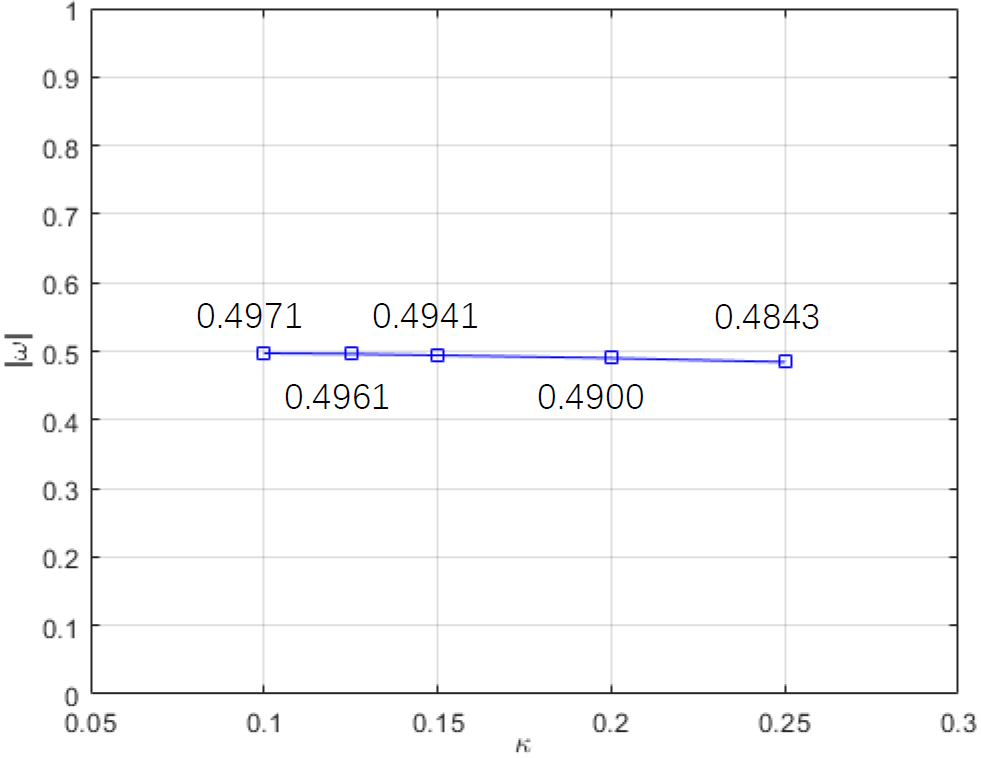}
    \caption{The simple shear flow-particle interaction problem:
    the computational domain (left);
    the velocity field of the case $\kappa=0.2$ (middle);
    the angular speed $|\omega|$ versus the confined ratio $\kappa$ (right).}\label{fig:benchmark1}
\end{figure}
We introduce the confined ratio defined as
$\kappa=2r/H$, where $r$ is the radius of circular particle.
Thus $\kappa=r$ in this example. Further, we take the fluid viscosity as $\mu_f=1$,
and densities of both the fluid and particle as $\rho_f=\rho_s=1$.
Finally, different from \cite{Pan;Glowinski2015} where the initial condition
of fluid velocity is defined by a simple shear flow associated with
the given shear rate without particle, here we treat the initial velocity field
of both the fluid and particle at rest.

Under the shear flow condition, it is known that the relationship
between the angular velocity and rotational angle of an elliptic shape particle
in an unbounded shear flow is given as follows according to Jeffery's solution \cite{jeffery1922motion}
\begin{equation}\label{ellp}
    \omega=-\gamma\frac{r_a^2\sin^2\theta+r_b^2\cos^2\theta}{r_a^2+r_b^2},
\end{equation}
which shows that the angular velocity is
related with the ratio of $r_a$ and $r_b$, theoretically, and is
independent of their specific values.
Nevertheless, the values of $r_a$ and $r_b$ used in the
numerical experiment should be as small as possible as Jeffery's solution
is derived under the assumption that the range of the domain is much
larger than the particle size so that to minimize the fluid-particle
interaction effect. Hence in this example, the Jeffery's exact solution
of the angular speed of a circular particle is $|\omega|=|-\gamma/2|=0.5$
in an unbounded shear flow field. In our numerical experiments below,
we let the confined ratio $\kappa$ (thus the radius of particle $r$)
vary from 0.1 to 0.25 to see its influence on the numerical angular speed.

In contrast to the rigid particle that is used in this benchmark problem,
our structural equation (\ref{solid-equation}) defined in the presented
FSCI model however belongs to an incompressible neo-Hookean material.
To make an elastic structure approximate to a rigid body at utmost, we can choose
a sufficiently large shear modulus, $\mu_s$, in the constitutive relation
of the elastic structure. Now we are able to carry out the developed UFMFD-FEM
through Algorithm \ref{algorithm1-1} without
considering the structure collision effect to solve this benchmark problem
for its numerical angular speed $|\omega_h|$, by taking the time step size $\Delta t=0.001$,
$m=1$ in (\ref{FESpaces}), 
and $\zeta_0=2,\ \zeta_1=12$ and $\zeta_2=2$ as pressure stabilization
parameters in (\ref{stab_parameter}).

We first investigate the influences
of increasing shear modulus $\mu_s$ and decreasing mesh size $h$ on the angular
speed for a circular particle with a fixed radius of $0.2$, and the obtained
numerical results are shown in Table \ref{tab:benchmark1}, where we can see that
for each chosen $\mu_s$, the corresponding numerical angular speed,
$|\omega_h|$, becomes stably more accurate towards $0.5$
along with the decreasing $h$, showing a stable convergence with
respect to the mesh size $h$ for our developed method.
In addition, amongst three chosen shear moduli,
$\mu_s=10^8$ (Pa) seems large enough to deliver an accurate angular speed
especially for a smaller mesh size $h$. On the other hand, we also see from
numerical computations that a much larger $\mu_s$ may possibly make Newton's nonlinear
iteration described in Algorithm \ref{algorithm1-1} harder to converge.
Therefore under the consideration of both accuracy and efficiency,
a reasonably large shear modulus such as $\mu_s=10^8$ (Pa) shall be
a good choice for the simulation of fluid-rigid particle interactions
by means of our proposed FSI model and UFMFD-FEM, in practice. A snapshot of the velocity field
in the case of $\kappa=r=0.2$ is also illustrated in the middle part of Figure \ref{fig:benchmark1}.
\begin{table}[H]
    \caption{Numerical angular speeds of a circular particle with a radius of 0.2
    on refining meshes under different shear moduli.\label{tab:benchmark1}}
\centering
\begin{tabular}[t]{ccc}
\toprule
    Mesh Size $h$ & Shear Modulus $\mu_s$ & Angular Speed $|\omega_h|$ \\
\midrule
     &   $10^6$ &      0.488623 \\
$\frac{1}{32}$ & $10^8$   &   0.487293 \\
    & $10^{10}$   &      0.487285 \\
\midrule
     &   $10^6$   &      0.489294 \\
$\frac{1}{64}$    & $10^8$   &  0.489108 \\
     &   $10^{10}$   &      0.489107 \\
\midrule
     & $10^6$   &  0.490046 \\
$\frac{1}{128}$     &   $10^8$ &  0.490028\\
     &   $10^{10}$ &  0.490038 \\
\bottomrule
\end{tabular}
\end{table}

Hence next, we take $h =1/128$, $\Delta t=0.001$ and $\mu_s=10^8$ to investigate
the influence of the confined ratio $\kappa$ on the numerical angular speed.
The right part of Figure \ref{fig:benchmark1} illustrates the relationship
of angular speeds versus the confined ratios, where we can see that
numerical angular speeds turn to well agree with the Jeffery's solution along with
decreasing particle sizes.
Note that relative errors on the numerical angular speed for the cases of
particle's radii $r =0.1, 0.125$ and $0.15$ are less than or equal to $1\%$,
while they become larger with increasing radius of the circular particle,
which can be ascribed to the increasingly subtle influence of
fluid-particle interaction effect due to the increasing particle size.

Finally, to investigate the numerical convergence rate between the obtained
numerical solutions $(\bv_h,p_h)$ and real solutions $(\bv,p)$
that are impossibly prescribed for the studied benchmark problem,
we conduct an error estimate based upon numerical solutions on a sequence of nested meshes with the mesh sizes
$2^jh,\ j=0,1,2,\cdots$ as follows,
$\|\phi_{2^{j}h}-\phi_{2^{j+1}h}\|_{L^2(\Omega)}\leq
\|\phi-\phi_{2^{j}h}\|_{L^2(\Omega)}+\|\phi-\phi_{2^{j+1}h}\|_{L^2(\Omega)}$,
which is under the circumstance that our proposed
numerical scheme converges to the real solution of the benchmark problem,
and here $\phi$ can be $\bv$ or $p$.
Then, while $\Delta t$ is sufficiently small, apply (\ref{errorestimate})
to two adjacent mesh levels with the mesh size $2^{j}h$ and $2^{j+1}h$, respectively, yields
$$
\|\bv_{2^{j}h}-\bv_{2^{j+1}h}\|_{L^2(\Omega)^d}
=O\big(2^{(j+1)r}h^r\big),\quad
\|p_{2^{j}h}-p_{2^{j+1}h}\|_{L^2(\Omega)}
=O\big(2^{(j+1)q}h^{q}\big),\quad \text{for }1< r \leq 3/2,\ 0 < q \leq 1/2.
$$
Thus, we have the following numerical convergence rate indicators of the velocity
and pressure attained on every three adjacent mesh levels,
\begin{equation}\label{estimator}
\frac{\|\bv_{2^{j+1}h}-\bv_{2^{j+2}h}\|_{L^2(\Omega)^d}}{\|\bv_{2^{j}h}-\bv_{2^{j+1}h}\|_{L^2(\Omega)^d}}=2^r,
\quad
\frac{\|p_{2^{j+1}h}-p_{2^{j+2}h}\|_{L^2(\Omega)}}{\|p_{2^{j}h}-p_{2^{j+1}h}\|_{L^2(\Omega)}}=2^q,
\end{equation}
which demonstrate the $r$th-order convergence rate of velocity in $L^2$ norm and
the $q$th-order convergence rate of pressure in $L^2$ norm as well, where
$1< r \leq 3/2$ and $0 < q \leq 1/2$.

Applying (\ref{estimator}) to numerical results of the developed UFMFD-FEM obtained
for the benchmark problem on every three adjacent meshes of the aforementioned
mesh doubling, we eventually gain numerical convergence rates of the velocity and
pressure as shown in Table \ref{tab:convergence-rate0}, where the convergence errors
$\|\bv_{2^{j-1}h}-\bv_{2^{j}h}\|_{L^2(\Omega)^d}$ and $\|p_{2^{j-1}h}-p_{2^{j}h}\|_{L^2(\Omega)}$
with $h=1/256$ and $j=3,2,1$, reversely, are indeed reduced by factors of
at least $2^{3/2}$ and of at least $2^{1/2}$ (indicated by values in columns labeled ``Ratio''),
respectively, for each doubling of resolution.
Thus, we can conclude that the obtained numerical solutions
are accurate by holding the convergence
rates of order $\frac{3}{2}$ and of order $\frac{1}{2}$ in $L^2$ norm for
velocity and pressure, respectively,
which are consistent with the supposed numerical convergence property
(\ref{errorestimate}) and the associated convergence rate indicators (\ref{estimator}).
\begin{table}[!htb]
\caption{Convergence errors \& rates of velocity and pressure
on successively nested grids for Example 1}\label{tab:convergence-rate0} \noindent
\begin{center}
\begin{tabular}{|c|c|c|c|c|c|}
\hline \multicolumn{2}{|c|}{Mesh Sizes} & \multicolumn{2}{|c|}{Velocity} & \multicolumn{2}{|c|}{Pressure} \\
\cline{1-6}
$j$ & \multicolumn{1}{|c|}{$2^jh$} &  $\|\bv_{2^{j-1}h}-\bv_{2^{j}h}\|_{L^2(\Omega)^d}$ & \multicolumn{1}{|c|}{Ratio} & $\|p_{2^{j-1}h}-p_{2^{j}h}\|_{L^2(\Omega)}$ &
\multicolumn{1}{|c|}{Ratio } \\ \hline
\multicolumn{1}{|c|} {$3$}& $1/32$ &
\multicolumn{1}{|c|}{$2.4000\times10^{-2}$} & 2.62 &
\multicolumn{1}{|c|}{$2.1471\times10^{-1}$} & 1.67 \\
\hline
\multicolumn{1}{|c|} {$2$} & $1/64$ &
\multicolumn{1}{|c|}{$9.1558\times10^{-3}$} & 3.50 &
\multicolumn{1}{|c|}{$1.2864\times10^{-1}$} & 1.56 \\
\hline
\multicolumn{1}{|c|}{$1$}& $1/128$ &
\multicolumn{1}{|c|}{$ 2.6157\times10^{-3}$} & --- &
\multicolumn{1}{|c|}{$ 8.2492\times10^{-2}$} & --- \\
\hline
\multicolumn{1}{|c|}{$0$} & $1/256$ &
\multicolumn{1}{|c|}{---} &--- &
\multicolumn{1}{|c|}{---} &---
\\ \hline
\end{tabular}
\end{center}
\end{table}

\subsection{Example 2: The plane Poiseuille flow-particle interaction benchmark problem}\label{sec:benchmark2}
\begin{figure}[htb]
\centerline{
    \includegraphics[height=4cm]{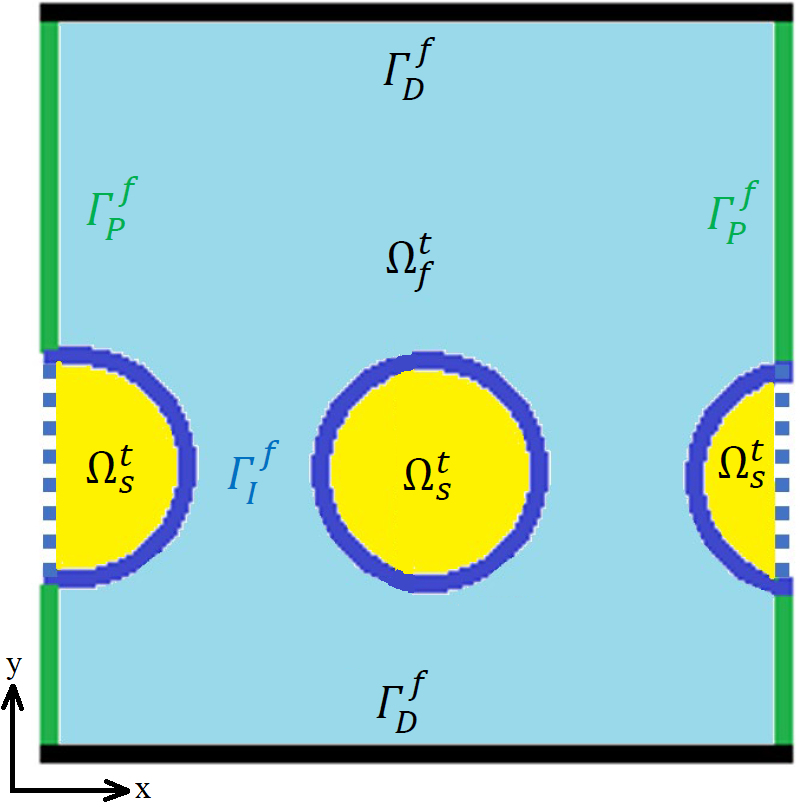}
    \hspace{0.2in}
    \includegraphics[height=4cm]{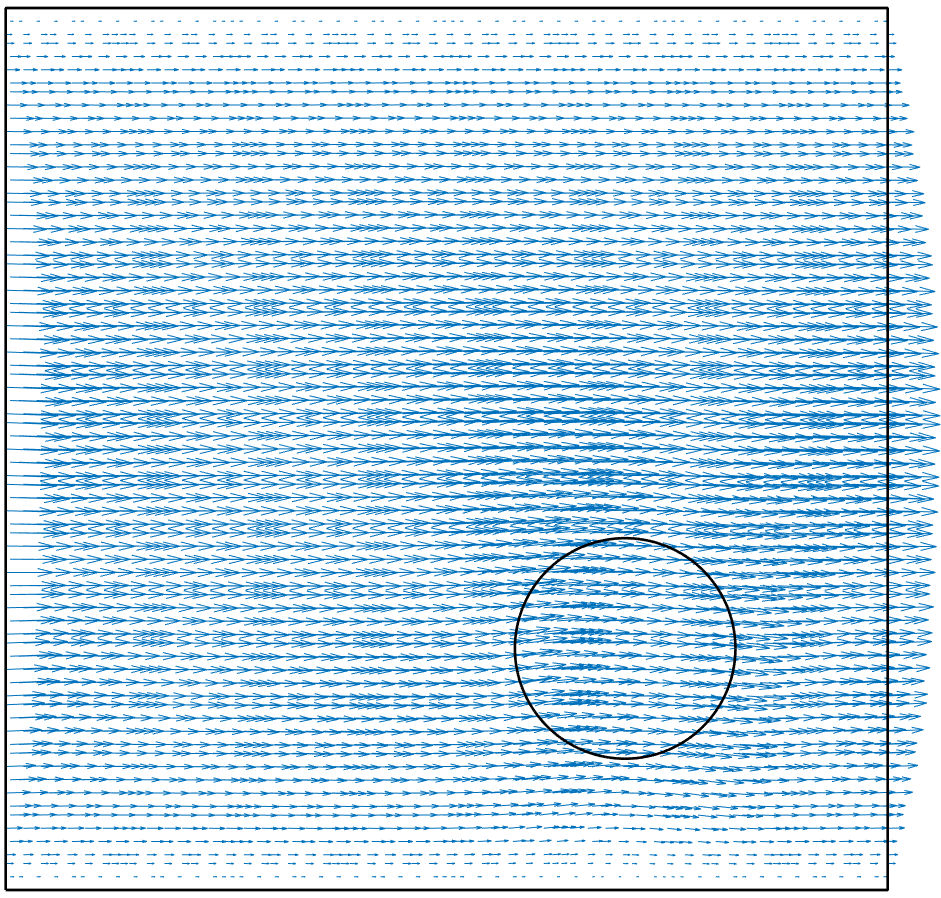}
    \hspace{0.2in}
    \includegraphics[height=4cm]{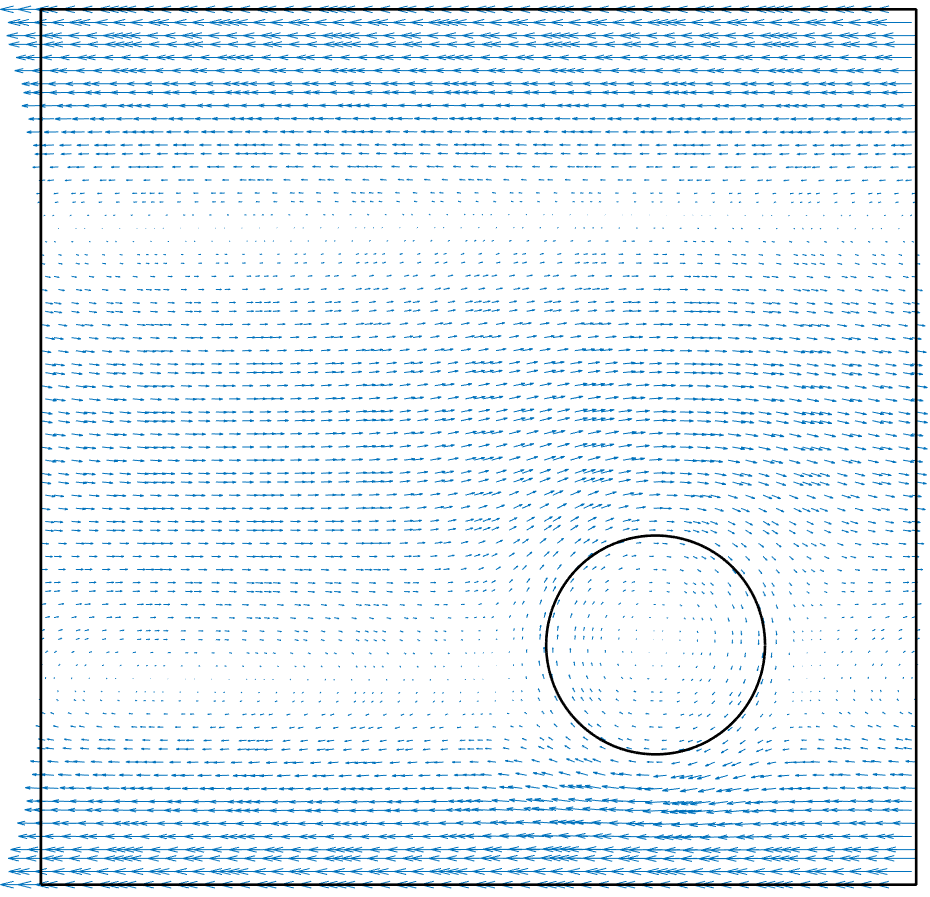}}
\centerline{(a)\hspace{4.5cm}(b)\hspace{4.5cm}(c)}
\vspace{0.2in}
\centerline{
    \includegraphics[height=4cm]{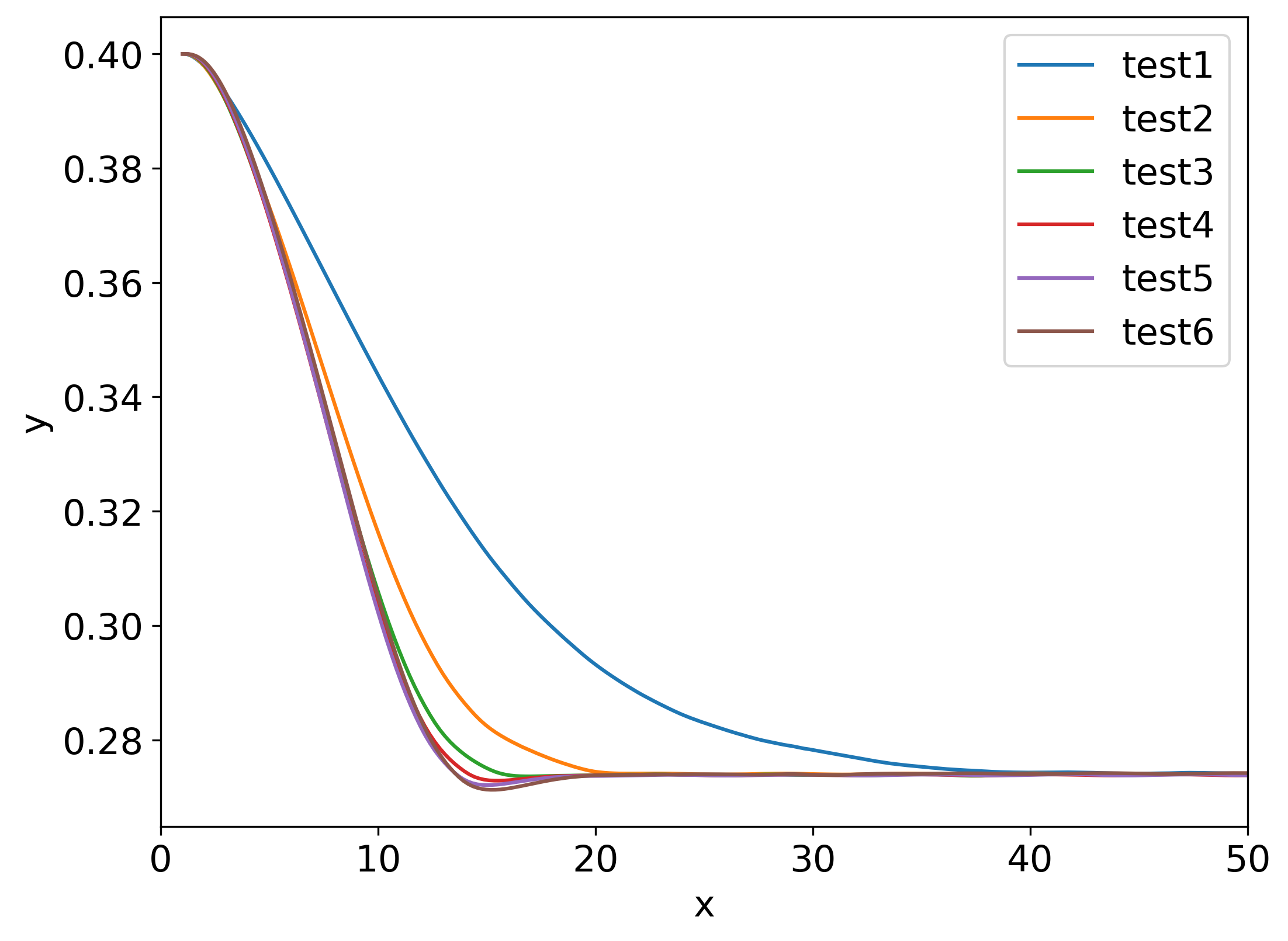}
    \hspace{0.2in}
    \hspace{0.0in}
    \includegraphics[height=4cm]{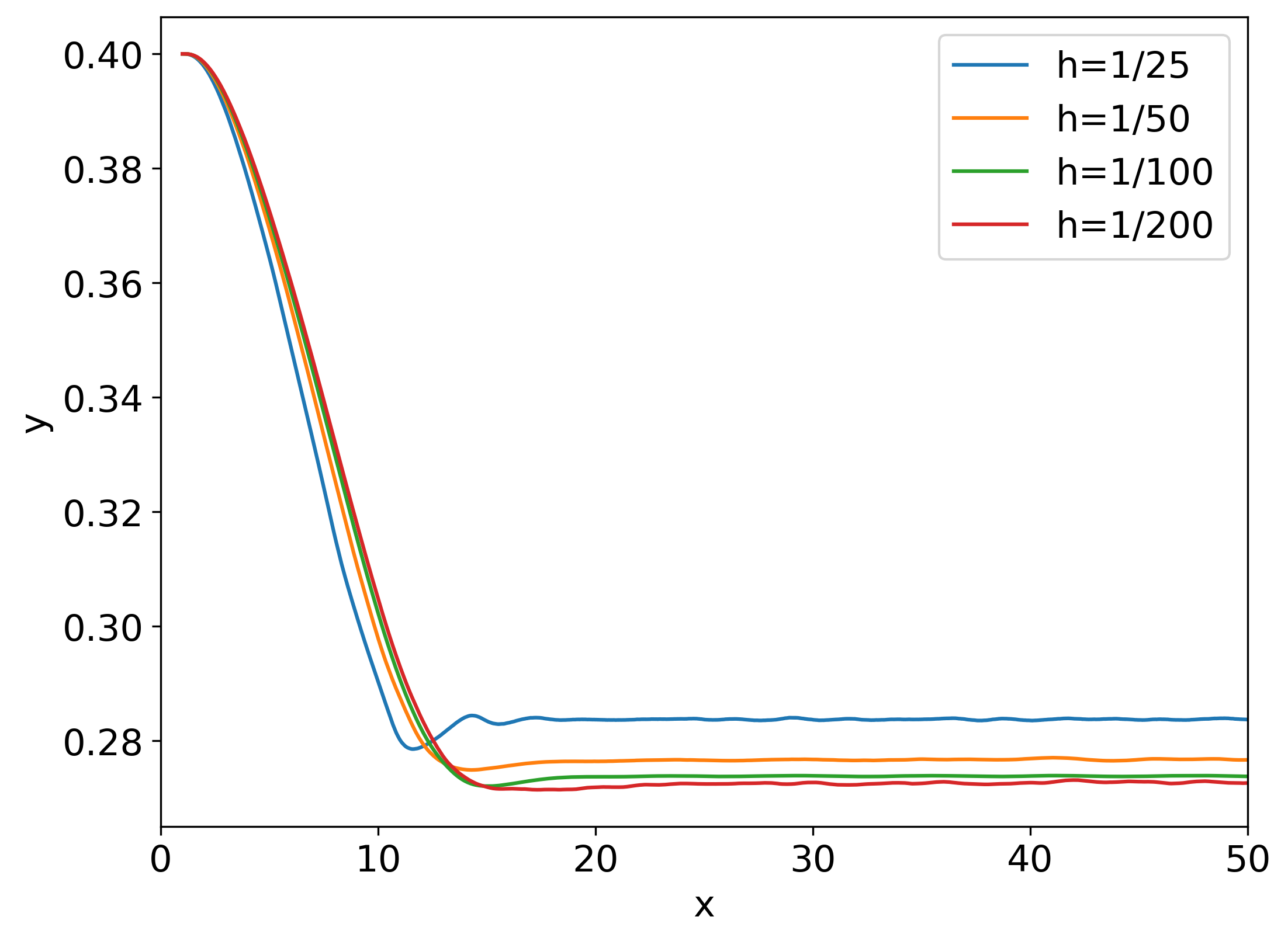}}
\centerline{\hspace{0.5cm}(d)\hspace{5.8cm}(e)}
    \caption{The plane Poiseuille flow-particle interaction problem:
    (a) the computational domain with the length $L=1$;
    (b) the total velocity field at the equilibrium position in Test \#5;
    (c) the rotating velocity field surrounding the particle at the equilibrium position
    in Test \#5 after an average (translational) velocity of the particle
    is deducted from the total velocity field;
    (d) lateral migrations of a circular particle in different tests;
    (e) lateral migrations of a circular particle with different mesh sizes in Test \#5.
    }\label{fig:benchmark2}
\end{figure}

In the second example, we consider to simulate the motion of
a neutrally buoyant rigid particle in a pressure-driven Poiseuille flow,
which was studied by Pan et al. in \cite{Pan;Glowinski2002} using the DLM/FD method as well,
and also by Inamuro et al. in \cite{Inamuro2000} via the lattice Boltzmann method.
Here we employ the same sets of parameters as in \cite{Pan;Glowinski2002} to
validate our developed UFMFD-FEM. As shown in Figure \ref{fig:benchmark2}(a),
the computational domain is $\Omega:= [0,1] \times [0,1]$ (i.e., $L=1$), where the periodic- and
no-slip velocity boundary conditions are assigned to the left and right boundary,
$\Gamma_P^f$, and to the top and bottom boundary, $\Gamma_D^f$, respectively.
In addition, to drive the fluid flow from the left to the right,
we introduce a pressure drop, $\Delta p$ (Pa), to this example which is equivalent with
a compound Neumann (stress) boundary condition on $\Gamma_P^f$, and,
we take $\Delta p$ as variable values in different testing cases (see Table \ref{tab:benchmark2}).
The fluidic and particle velocities are set as zero, initially, i.e., $\bv_f^0=\bv_s^0=0$,
the circular particle's diameter is taken as $0.25$, and its center's
vertical coordinate as $y_c=0.4$, initially.
Densities of the fluid and the particle are the same and both are set as $1$,
and, the fluid viscosity $\mu_f$ (kg/m/s) varies as well in different testing cases (see Table \ref{tab:benchmark2}).
Again, we take the structural shear modulus $\mu_s=10^8$ (Pa),
and employ Algorithm \ref{algorithm1-1} to solve this benchmark problem,
where the pressure-stabilized $P_1/P_1$ mixed finite element is adopted
with well tuned stabilization parameters, as chosen in Section \ref{sec:benchmark1}.
\begin{table}[H]
    \caption{Parameters and numerical results in six testing cases
    compared with \cite{Pan;Glowinski2002} denoted as [$\,\bullet\,$].\label{tab:benchmark2}}
\centering
{\small\begin{tabular}[t]{ccccccccc}
\toprule
 Test\#  & $\mu_f$ & $\Delta p$ & [$\bar u$] & $\bar u$ & [$y_c$]&$y_c$
    & [$\omega$] & $\omega$\\
\midrule
  1  &  $3.2498036\times10^{-3}$ & $1.763\times10^{-3}$ & 0.04155 &0.04131 & 0.2732 &0.2735 & -0.05345 &-0.05315 \\
  2  &  $1.5000000\times10^{-3}$ & $8.167\times10^{-4}$ & 0.04159 &0.04141 & 0.2725 &0.2731 & -0.05343 &-0.05304 \\
  3  &  $9.4984908\times10^{-4}$ & $5.133\times10^{-4}$ & 0.04122 &0.04104 & 0.2723 &0.2727 & -0.05264 &-0.05215 \\
  4  &  $7.5000000\times10^{-4}$ & $4.100\times10^{-4}$ & 0.04166 &0.04146 & 0.2720 &0.2724 & -0.05283 &-0.05222 \\
  5  &  $6.0000000\times10^{-4}$ & $3.270\times10^{-4}$ & 0.04148 &0.04114 & 0.2719 &0.2723 & -0.05206 &-0.05135 \\
  6  &  $4.2834760\times10^{-4}$ & $2.337\times10^{-4}$ & 0.04144 &0.04116 & 0.2722 &0.2728 & -0.05052 &-0.04962 \\
\bottomrule
\end{tabular}}
\end{table}

First of all, we fix $h=1/200$ and let $\Delta t=0.05$ before the particle reaches the equilibrium position,
thereafter $\Delta t=0.002$, as chosen in \cite{Pan;Glowinski2002}.
The obtained numerical results of six testing cases with different values of pressure drop and
of fluid viscosity as used in \cite{Pan;Glowinski2002} are shown in Table \ref{tab:benchmark2},
where the space-averaged inlet velocity $\bar u$ (m/s), the equilibrium position and
angular velocity of the circular particle, $y_c$ (m) and $\omega$ (rad/s), are listed and
hold the relative errors from corresponding results in \cite{Pan;Glowinski2002}
(denoted as [$\bar u$], [$y_c$] and [$\omega$] in Table \ref{tab:benchmark2})
between $0.43\%\sim0.82\%$, $0.11\%\sim0.22\%$, and $0.56\%\sim1.78\%$, respectively,
illustrating that the particle's angular velocity $\omega$ is influenced
by Reynolds numbers of the fluid (Re = $\rho_f\bar uL/\mu_f$) the most, and, its relative errors
even increase along with Reynolds numbers, whereas the particle's equilibrium
position $y_c$ possesses the smallest error range and is thus affected
by Reynolds numbers the least. The same fact is also displayed
in Figure \ref{fig:benchmark2}(d), where we can see all lateral migration
curves of the evolving particle, which start from the same initial position
in six testing cases, eventually approach to the same equilibrium position,
resulting in the following phenomenon that the equilibrium state of the
evolving particle is independent of Reynolds numbers of the plane Poiseuille flow.
Figure \ref{fig:benchmark2}(b) shows a snapshot of velocity field
at the equilibrium position in Test \#5, illustrating that the particle mainly
conducts a translational motion due to the horizontal fluid flow
driven by the pressure drop from the left to the right. To show the particle
also does a slightly rotational motion at the same time in Test \#5, we first
calculate an average velocity over totally $N_s$ mesh nodes of the particle
with an equal weight $1/N_s$, which can be considered
as an approximation to the translational velocity of the rigid particle.
Then, we subtract such an average velocity of the particle
from the velocity field shown in Figure \ref{fig:benchmark2}(b) to obtain
a clockwise rotational velocity field (the angular velocity is thus negative)
surrounding the center of particle as shown in Figure \ref{fig:benchmark2}(c),
and, the rotational velocity field fades away as approaching to the center of particle,
showing a decreasing rotational velocity field along the radius direction
that is nearly anti-proportional to the distance from the center of particle.

Next, to investigate the convergence property of our proposed method
on this example, without loss of generality, we pick Test \#5 in
particular, and produce a series of
successively nested meshes with $h=\frac{1}{25},\ \frac{1}{50},\ \frac{1}{100},
\ \frac{1}{200}$, on which we carry out Algorithm \ref{algorithm1-1} for Test \#5
with the time step size $\Delta t = 0.002$ to
gain numerical solutions, respectively. Then,
as done in Section \ref{sec:benchmark1}, we investigate the numerical
convergence rate between the obtained numerical solutions $(\bv_h,p_h)$
and the real solution $(\bv,p)$
by checking if (\ref{estimator}) still holds amongst every three adjacent mesh levels
for this benchmark problem. As Table \ref{tab:convergence-rate} shows,
we observe that the convergence errors
$\|\bv_{2^{j-1}h}-\bv_{2^{j}h}\|_{L^2(\Omega)^d}$ and $\|p_{2^{j-1}h}-p_{2^{j}h}\|_{L^2(\Omega)}$
with $h=1/200$ and $j=3,2,1$, reversely, are still roughly reduced by factors of
$2^{3/2}$ and of $2^{1/2}$ or so, respectively, for each doubling of resolution, illustrating that
numerical solutions of the benchmark problem still roughly hold the convergence rates
of order $\frac{3}{2}$ and of order $\frac{1}{2}$ in $L^2$ norm for
velocity and pressure, respectively, which agree with the supposed
numerical convergence property (\ref{errorestimate}) and the associated
convergence rate indicators (\ref{estimator}), again.
\begin{table}[!htb]
\caption{Convergence errors \& rates of velocity and pressure
on successively nested grids for Example 2}\label{tab:convergence-rate} \noindent
\begin{center}
\begin{tabular}{|c|c|c|c|c|c|}
\hline \multicolumn{2}{|c|}{Mesh Sizes} & \multicolumn{2}{|c|}{Velocity} & \multicolumn{2}{|c|}{Pressure} \\
\cline{1-6}
$j$ & \multicolumn{1}{|c|}{$2^jh$} &  $\|\bv_{2^{j-1}h}-\bv_{2^{j}h}\|_{L^2(\Omega)^d}$ & \multicolumn{1}{|c|}{Ratio} & $\|p_{2^{j-1}h}-p_{2^{j}h}\|_{L^2(\Omega)}$ &
\multicolumn{1}{|c|}{Ratio } \\ \hline
\multicolumn{1}{|c|} {$3$}& $1/25$ &
\multicolumn{1}{|c|}{$1.7655\times10^{-3}$} & 2.57 &
\multicolumn{1}{|c|}{$4.4926\times10^{-5}$} & 2.36 \\
\hline
\multicolumn{1}{|c|} {$2$} & $1/50$ &
\multicolumn{1}{|c|}{$6.8578\times10^{-4}$} & 2.14 &
\multicolumn{1}{|c|}{$1.9036\times10^{-5}$} & 1.54 \\
\hline
\multicolumn{1}{|c|}{$1$}& $1/100$ &
\multicolumn{1}{|c|}{$3.2107\times10^{-4}$} & --- &
\multicolumn{1}{|c|}{$1.2399\times10^{-5}$} & --- \\
\hline
\multicolumn{1}{|c|}{$0$} & $1/200$ &
\multicolumn{1}{|c|}{---} &--- &
\multicolumn{1}{|c|}{---} &---
\\ \hline
\end{tabular}
\end{center}
\end{table}
On the other hand, Figure \ref{fig:benchmark2}(e) illustrates lateral migration curves
of the circular particle from the same initial position
over four different mesh sizes, where we can see that along with the decreasing mesh size
the equilibrium states of four curves are farther away from the center axis
of the fluidic channel (i.e., the horizontal line $y=0.5$) one after the other,
and, the distance between every two adjacent curves' equilibrium state
also keeps decreasing while the mesh size decreases, showing also an approximation
process with respect to the mesh size $h$ as well, qualitatively.

\subsection{Example 3: A self-defined FSCI problem with particle collisions
over the fluidic channel wall}\label{sec:selfdefinedFSCI}
\begin{figure}[htb]
    \centering
    \includegraphics[height=6cm]{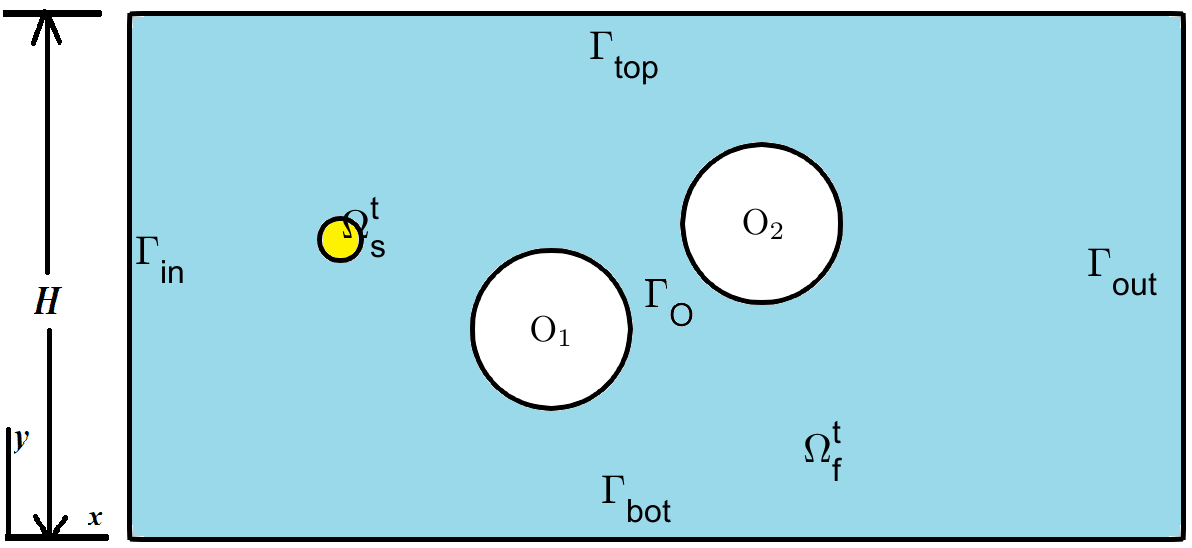}
    \caption{Computational domain of a self-defined FSCI problem.}\label{fig:fig7}
\end{figure}

In this example, we numerically test a FSCI model problem with contact effects between
the structure and fluidic channel wall
occurring in $\Omega=[0,2]\times[0,1]\backslash \left(\text{O}_1\cup \text{O}_2\right)$ (i.e., $L=2,\ H=1$),
where O$_1$ and O$_2$ are two circular pillar obstacles with the same radius $0.15$, and
centered at $(0.8,0.4)$ and $(1.2,0.6)$, respectively, and,
$\Gamma_{\text{O}}=\partial\text{O}_1\cup\partial\text{O}_2$. The structure is a
circular particle with radius $0.04$ that is initially positioned at $(0.4,0.57)$,
as shown in Figure \ref{fig:fig7}. Note that in this example
units of length and time are micrometers ($\mu$m) and microsecond (ms), respectively.
All involved physical parameters are listed in Table \ref{tab:parameters}.
\begin{table}[H]
    \caption{{Symbol descriptions and values in unit for physical parameters}}\label{tab:parameters}
    \newcolumntype{C}{>{\centering\arraybackslash}}
    \begin{tabularx}{\textwidth}{@{} Cp{0.1\textwidth} Cp{0.4\textwidth} Cp{0.2\textwidth} Cp{0.2\textwidth} @{}}
    \toprule
      Symbol & Description & Value &Unit \\
    \midrule
      $\rho_f$ & Density of fluid & $1.005584\times{10}^{-12}$  & g$\cdot\mu$m$^{-3}$ \\
      $\mu_f$ & Viscosity of fluid & $1.0219\times 10^{-9}$ & g$\cdot\mu$m$^{-1}\cdot$ms$^{-1}$ \\
      $\rho_s$ & Density of particle & $1.05\times{10}^{-12}$ &  g$\cdot \mu$m$^{-3}$ \\
      $\mu_s$ & Shear modulus of particle & $1$ & g$\cdot\mu$m$^{-1}\cdot$ms$^{-2}$ \\
    \bottomrule
    \end{tabularx}
\end{table}

We set zero for initial conditions of all variables in the interior domain,
while subjecting to their own boundary conditions on the boundary as shown below.
And, the initial structural displacement is also at rest, i.e.,
the particle stays still initially. As for the setup of boundary conditions,
we define a parabolic profile for the incoming
flow on the inlet $\Gamma_{\text{in}}$ as follows
\begin{align}
\bv_f^{\text{in}}=\left(v_{f,x}^{\text{in}},0\right)^T, \text{ and, }
    v_{f,x}^{\text{in}}&=\left\{
    \begin{array}{ll}
        v_{\text{max}}\frac{2y(H-y)}{H^2}\left(1-\cos\left(\frac{\pi t}{t_0}\right)\right),& t\le t_0\label{inflow}\\
        v_{\text{max}}\frac{4y(H-y)}{H^2},& otherwise,
    \end{array}
    \right.
\end{align}
where $v_{\text{max}}=1 ~\mu \text{m}/\text{ms}$ denotes the maximum value of parabolic incoming
velocity profile, $t_0$ is a prescribed buffer time period within which the
fluid flow can be developed from its zero initial state
to a stable state throughout $\Omega_f^{t_0}$ with an increasing incoming
velocity boundary condition from zero to $\bv_f^{\text{in}}(\bx,t_0)$. Additionally, the do-nothing
boundary condition is applied to the outlet $\Gamma_{\text{out}}$, i.e.,
$\bff_{f,N}=0$ in (\ref{fluid-BICs})$_2$, and the no-slip boundary condition
(\ref{fluid-BICs})$_1$ with $\bv_{f,D}=0$ is defined on $\Gamma_{\text{top}}\cup
\Gamma_{\text{bot}}\cup \Gamma_{\text{O}}$.

\begin{figure}[htb]
    \centering
    \includegraphics[height=6cm]{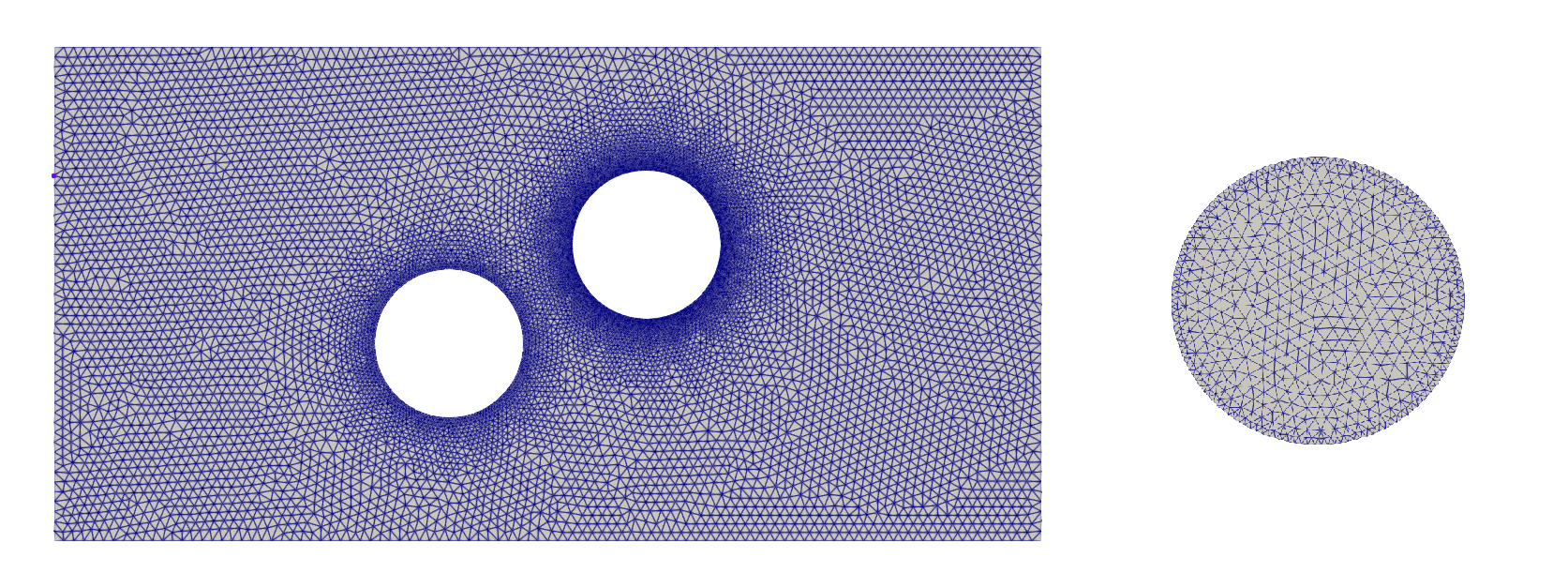}
    \caption{Meshes for $\Omega$ (left) and $\Omega_s^0$ (right).}\label{fig:fig8}
\end{figure}

As shown in Figure \ref{fig:fig8}, we triangulate the entire domain $\Omega$
with a fixed mesh $\mathcal{T}_h$, and the initial structure domain $\Omega_s^0$
with a mesh $\mathcal{T}_{h_s}^0$, where $\mathcal{T}_h$ consists of 9895 nodes
and 19042 elements, while $\mathcal{T}_{h_s}^0$ contains 660 nodes and 1180 elements.
The time step size is taken as $\Delta t = 0.01$ ms. Further, we
distribute a finer fluidic boundary-layer mesh along $\Gamma_{\text{O}}$ with a size
$h_{\text{C}} =0.003~\mu$m. The parameter $\zeta$ in \eqref{contact_force_num}
is set to be $10^4$, 
and the tolerance $\varepsilon$ in Step 4 of Algorithm \ref{algorithm3}
is chosen as $0.001 h_{\text{C}}{\Delta t}^{-1}$.

\begin{figure}[htb]
    \centering
    \includegraphics[height=6cm]{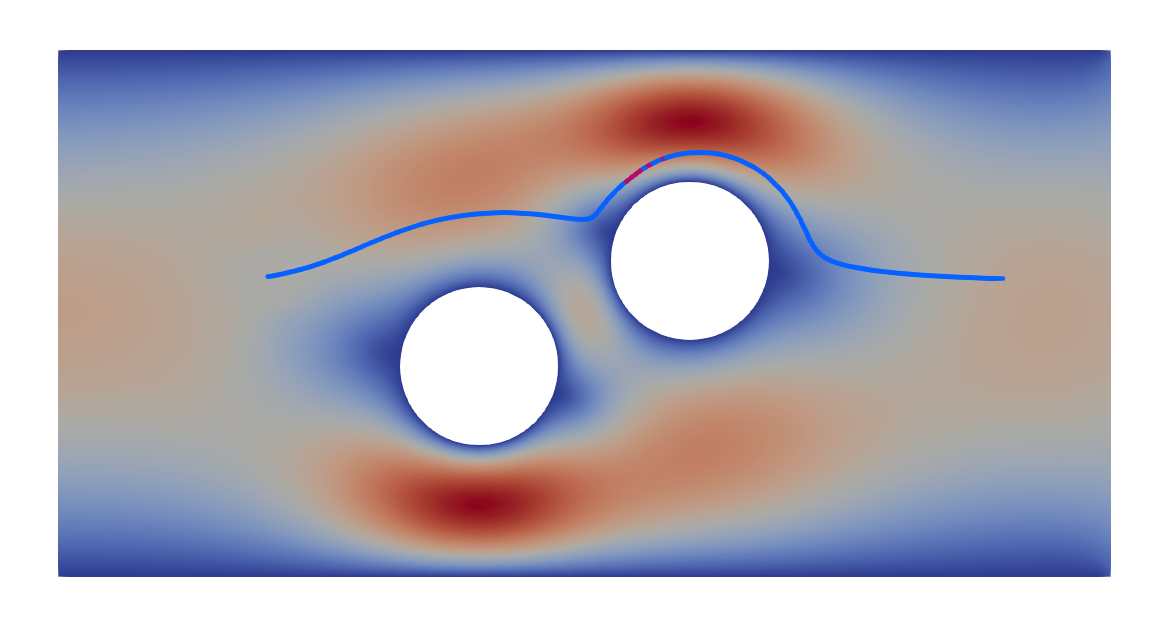}
    \caption{The lateral migration of particle with the field
    of velocity magnitude as the background.}\label{fig:fig9}
\end{figure}

\begin{figure}[htb]
    \centering
    \includegraphics[height=6cm]{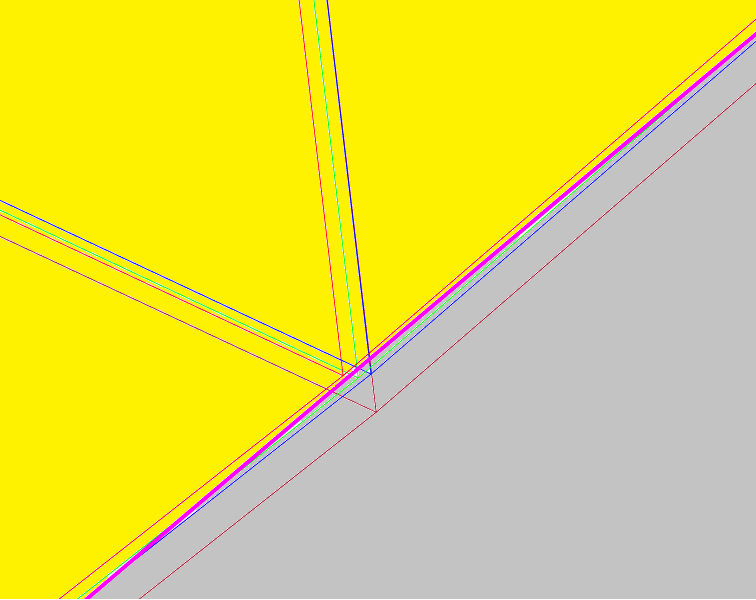}
\caption{The nonlinear iteration process of determining the repulsive contact force:
the respective location of structural elements at the beginning of the nonlinear iteration
(brown line), after 10 times iteration (blue line),
after 20 times iteration (green line), and the final location (red line),
where the grey area is a part of the fluidic boundary-layer mesh, the yellow area is a part
of the structural mesh, the thick purple line is a part of $\Gamma_h^{BL}$ that is the
boundary of first two-layer meshes from the fluidic channel wall within the fluidic boundary layer.}\label{fig:fig10}
\end{figure}

The lateral migration of particle is shown in Figure \ref{fig:fig9}
with the field of velocity magnitude as the background,
where some red parts on the trajectory of the particle's center
near the right pillar illustrate that the contacting phenomena
occur between the particle and the right pillar obstacle at those places.
A typical contacting process at one time step obtained by
Algorithm \ref{algorithm3} is shown in Figure \ref{fig:fig10}, where
the thin lines in brown, blue, green and red colors represent
the location of structural elements at the beginning of nonlinear iteration, after
10 and 20 times iteration, and at the end of (after 35 times) iteration,
respectively, until no structural element
falls into the fluidic boundary-layer mesh that is colored in grey and bounded by
a thick purple line. In Figure \ref{fig:fig10} one can observe that,
at the beginning of nonlinear iteration, some vertices of structural elements
enter the fluidic boundary layer, which indicates the contacting collision
between the structure and fluidic channel wall numerically occurs and thus
triggers Algorithm \ref{algorithm3}. The repulsive contact force $f_{s,C_W}$ and
the contacting surface $\Gamma_C^t$ are then updated and determined by the
nonlinear iteration process described in Algorithm \ref{algorithm3}.
According to \eqref{contact_force_num}, $f_{s,C_W}$ is gradually increased to
prevent the structure from penetrating into the fluidic boundary layer
when $\bv_{h_s}^{n+1,k+1}\cdot \bn\big|_{\Gamma_{C}^{n+1,k+1}}$ keeps positive.
The stopping criteria of Algorithm \ref{algorithm3}, $\Gamma_{C}^{n+1,k+1} =\emptyset$ or
$\bv_{h_s}^{n+1,k+1}\cdot \bn\big|_{\Gamma_{C}^{n+1,k+1}}<\varepsilon$, together with
\eqref{contact_force_num} that is to update the non-negative repulsive
contact force, imply that the original contact condition \eqref{struct-bc-2a}
numerically holds in an equivalent fashion.
Finally, Figure \ref{fig:fig10} also illustrates that the structure
is pushed away from the fluidic boundary layer as desired (shown
by the thin red line) at the end of nonlinear iteration process.

\subsection{Example 4: A realistic particulate FSCI problem in a deterministic lateral displacement microchip}
\begin{figure}[htb]
    \centering
    \includegraphics[height=6cm]{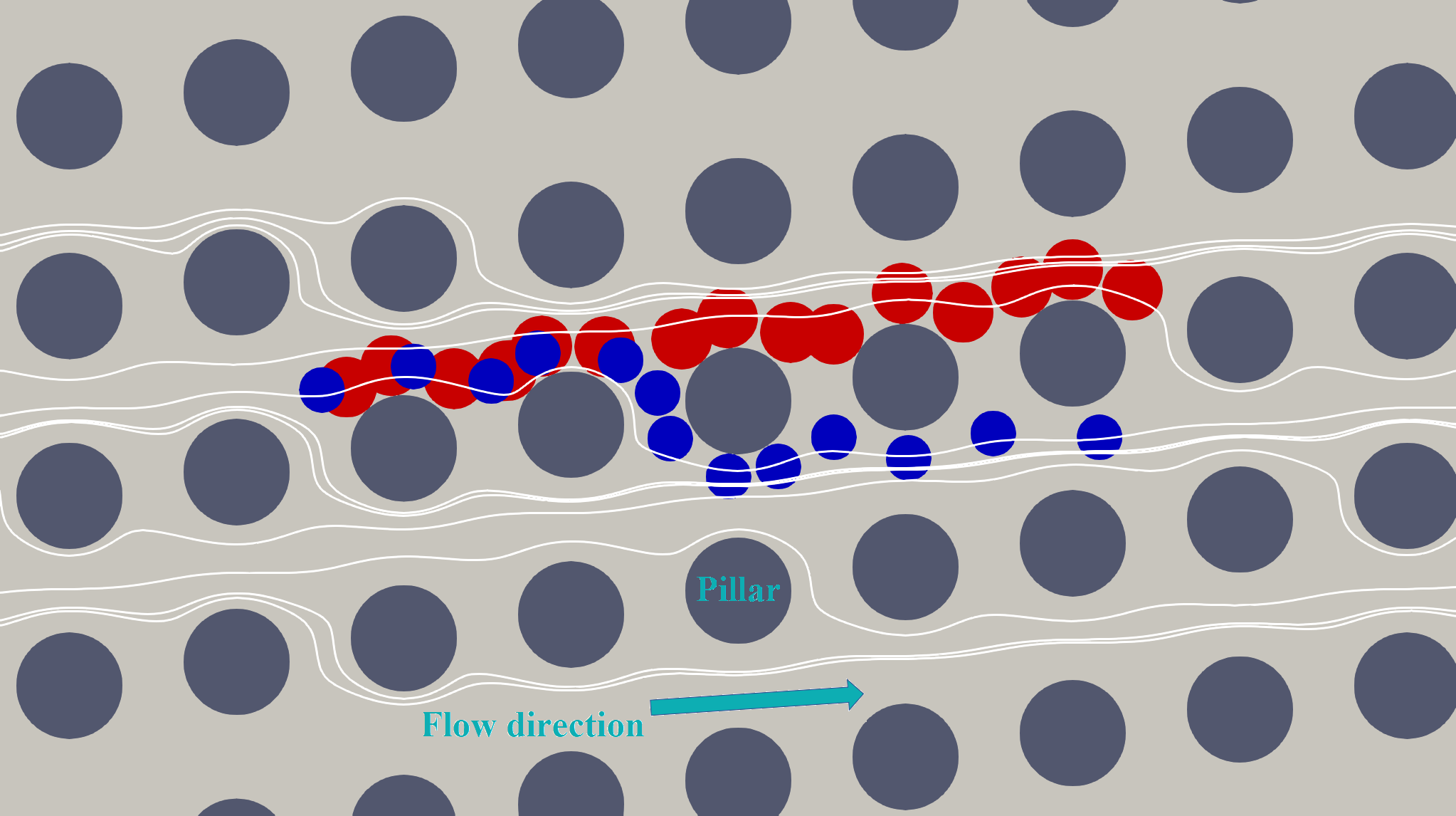}
    \caption{DLD principle:
    small particles whose sizes are less than the critical diameter follow
    the initial streamline and move
in a zigzag direction (blue), while large particles travel in bumping
mode (red) and move to the next streamline, 
where white streamlines of fluid flow are shown in the background.}\label{fig:DLD-criticaldiameter}
\end{figure}

In this section, we further validate the developed UFMFD-FEM by investigating its
numerical performances on a realistic particulate FSCI problem that arises
from a DLD microchip.
The DLD method is a robust passive microfluidic particle separation technique
established by Huang et al. \cite{Huang2004} for the first time to
sort particles based on their size with pillar arrays. During the last decade,
DLD has become popular and widely used for particle separation and detection
by holding a promise to provide a precise particle manipulation
with a high-resolution separation in a robust fashion and at low cost.
The mechanism of DLD is to use arrays of offset micro-pillars
within a flow channel to sort particles based on diameter in a high
throughput manner. DLD mostly operates at low Reynolds numbers and
provides high dynamic size separation, which ranges from millimeter
to micro- and nanometer sizes. Particle flows in the DLD array are
influenced by both the fluidic forces and the pillar obstacles effect.
A critical size for particle separation, which is so called
the DLD critical diameter denoted by D$_c$, is determined by
the gap between the pillars
and the angle of the pillar array in relation to the main direction of flow.
When the particle is located in the pillar gap, the particle with a diameter
smaller than D$_c$ will follow the initial streamline
and travel in the zigzag mode, while the particle larger than D$_c$
will be bumped to the pillar and displace laterally
to the next streamline, as shown in Figure \ref{fig:DLD-criticaldiameter}.

\begin{figure}[htb]
    \centering
    \includegraphics[height=3.5cm]{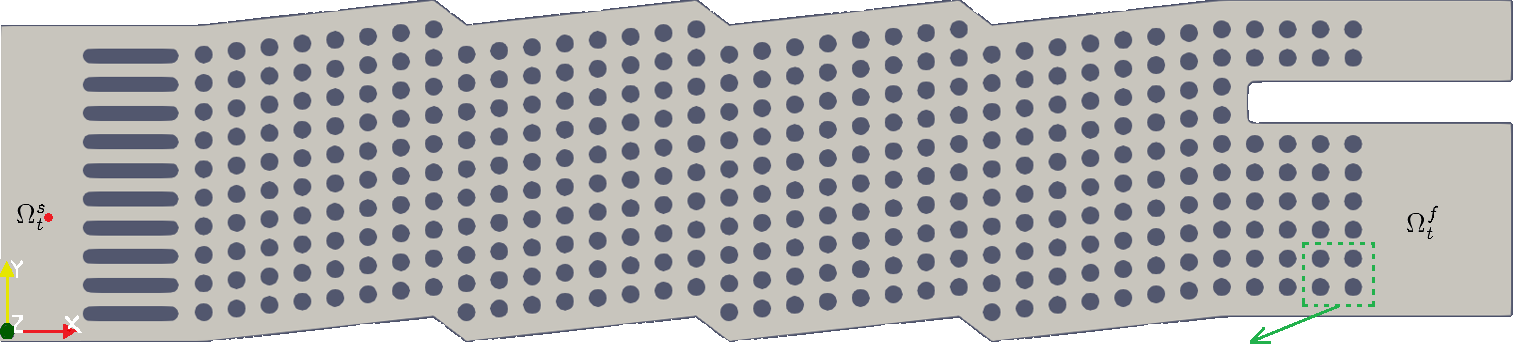}\\
    \includegraphics[height=5cm]{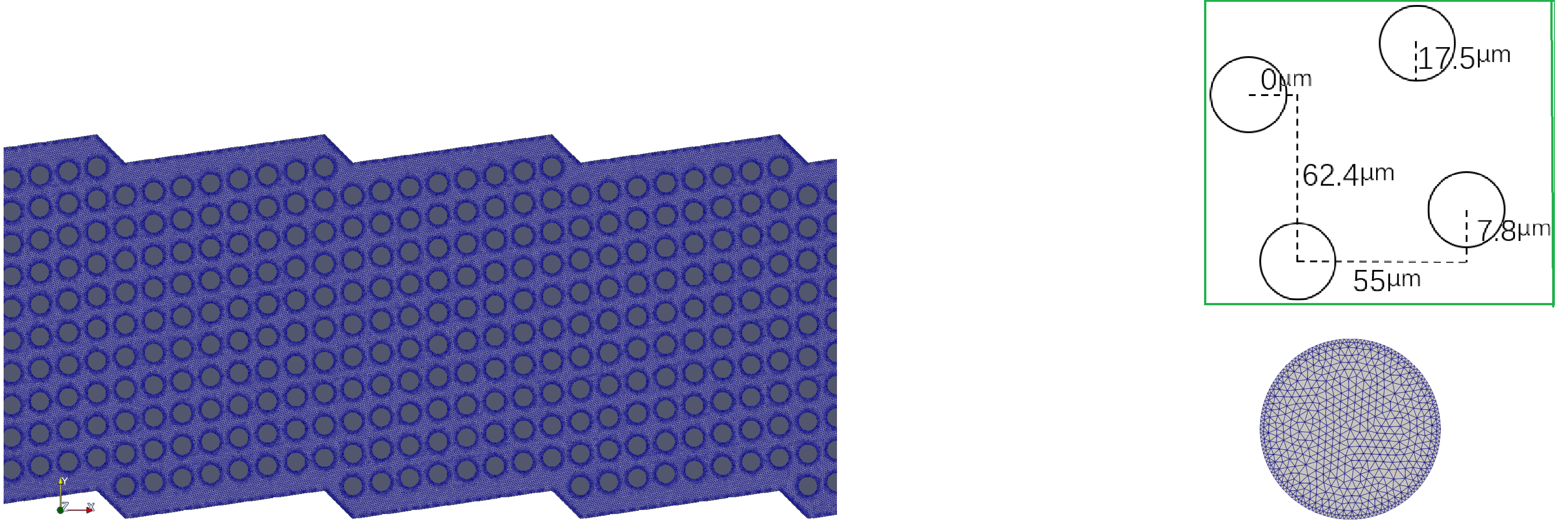}
    \caption{Computational domain of the FSI problem in a DLD microchip
    (top), a portion of its mesh triangulation (bottom left), and
    the structural mesh (bottom right).}\label{fig:DLD-mesh}
\end{figure}

A two-dimensional computational domain of DLD microchip as well as its mesh triangulation
are depicted in Figure \ref{fig:DLD-mesh},
where a matrix of pillar obstacles, i.e., the micro-posts array that are also
inner walls of the filter microfluidic channel $\Omega^f_t$, are built into the
microchip at an oblique angle to attempt to produce an effect of
particle isolation, and, such an oblique angle is particularly illustrated in
the green box shown in Figure \ref{fig:DLD-mesh}, where the shift distances
between adjacent pillar obstacles in $x$- and $y$-direction are $0~\mu$m and $7.8~\mu$m,
respectively, while the horizontal and vertical distances among adjacent
pillar obstacles are $55~\mu$m and $62.4~\mu$m, respectively.
The radius of each circular pillar obstacle is $17.5~\mu$m.
In addition, $\Omega^s_t$ denotes the particle that is plotted as a small red disk.
As for the mesh triangulation, we adopt $h=h_s=3~\mu$m to accommodate a
large-scale domain triangulation while making a fluidic boundary-layer
mesh along each pillar obstacle with a mesh size
$h_{\text{C}} =0.003~\mu$m, and, we take the time step size as $\Delta t=0.01$ms.
We further adopt the maximum incoming velocity on the inlet,
$v_{\text{max}}=50~\mu$m/ms.
Moreover, boundary conditions, initial conditions, involved physical coefficients
of this example, and control parameters for handling structural collisions
described in Algorithm \ref{algorithm3} are all set as the same with
those given in Section \ref{sec:selfdefinedFSCI}.

\begin{figure}[htb]
    \centering
    \includegraphics[height=4cm]{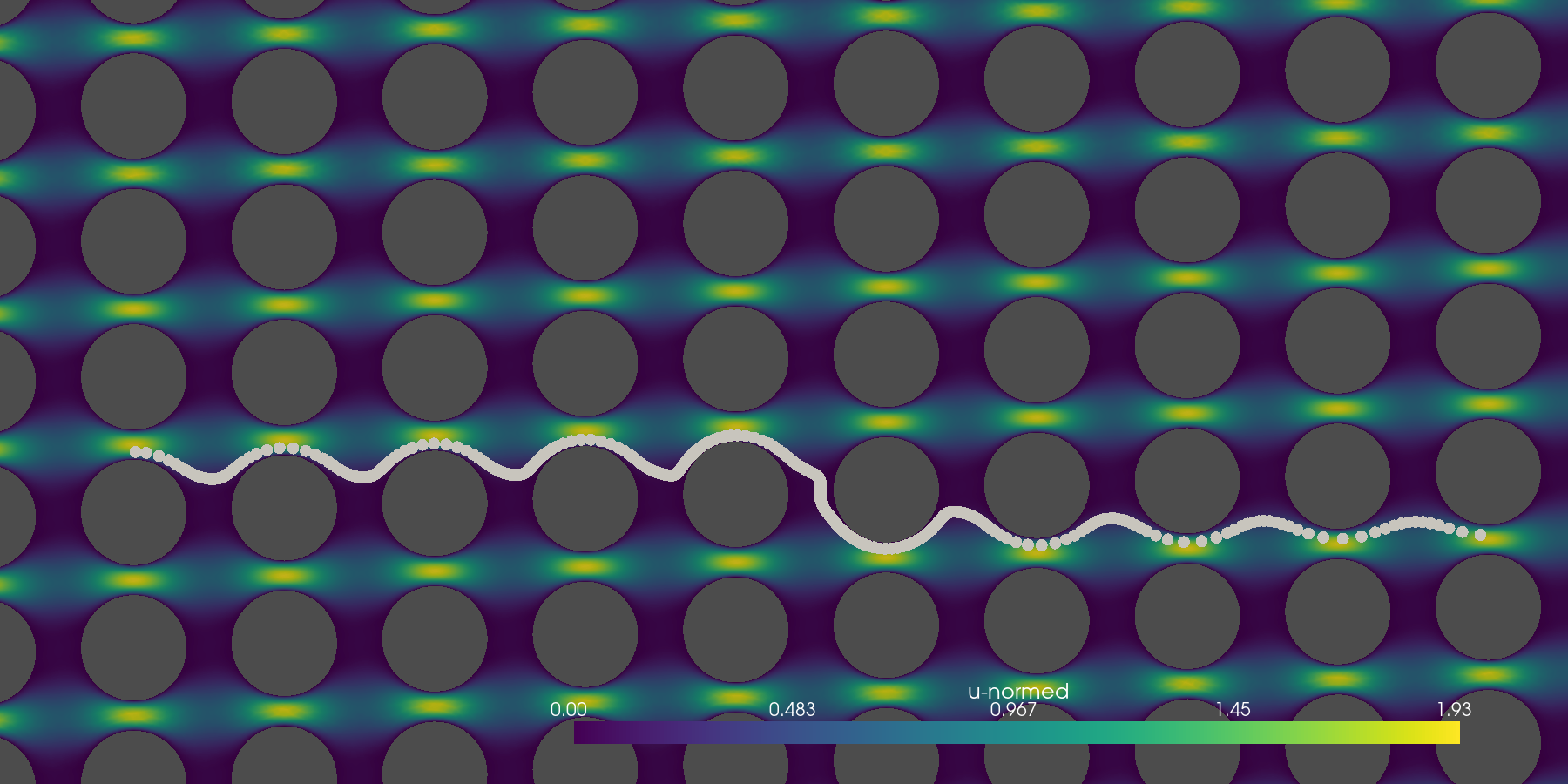}
    \hspace{0.5cm}
    \includegraphics[height=4cm]{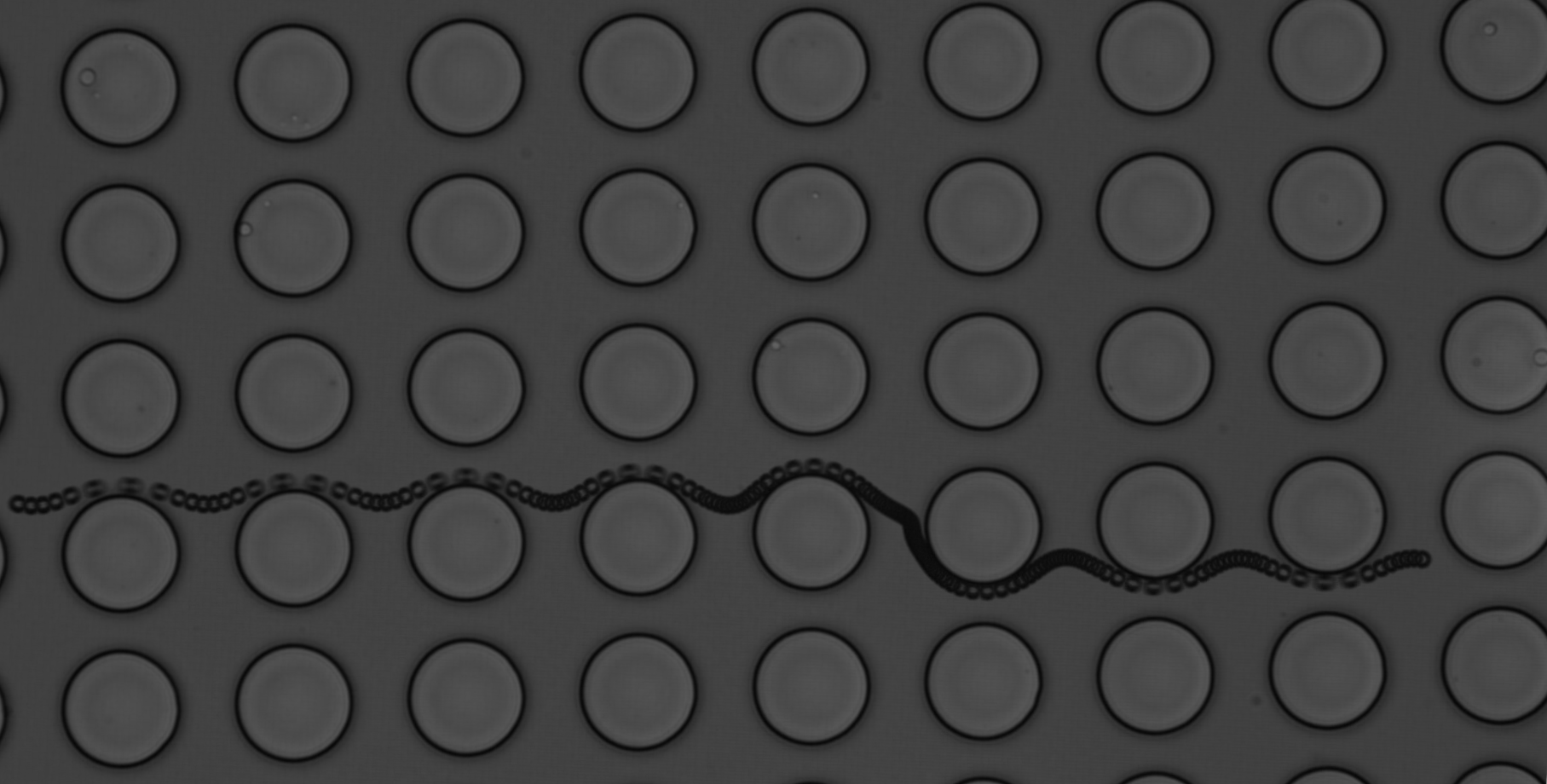}\\\vspace{0.5cm}
    \includegraphics[height=4cm]{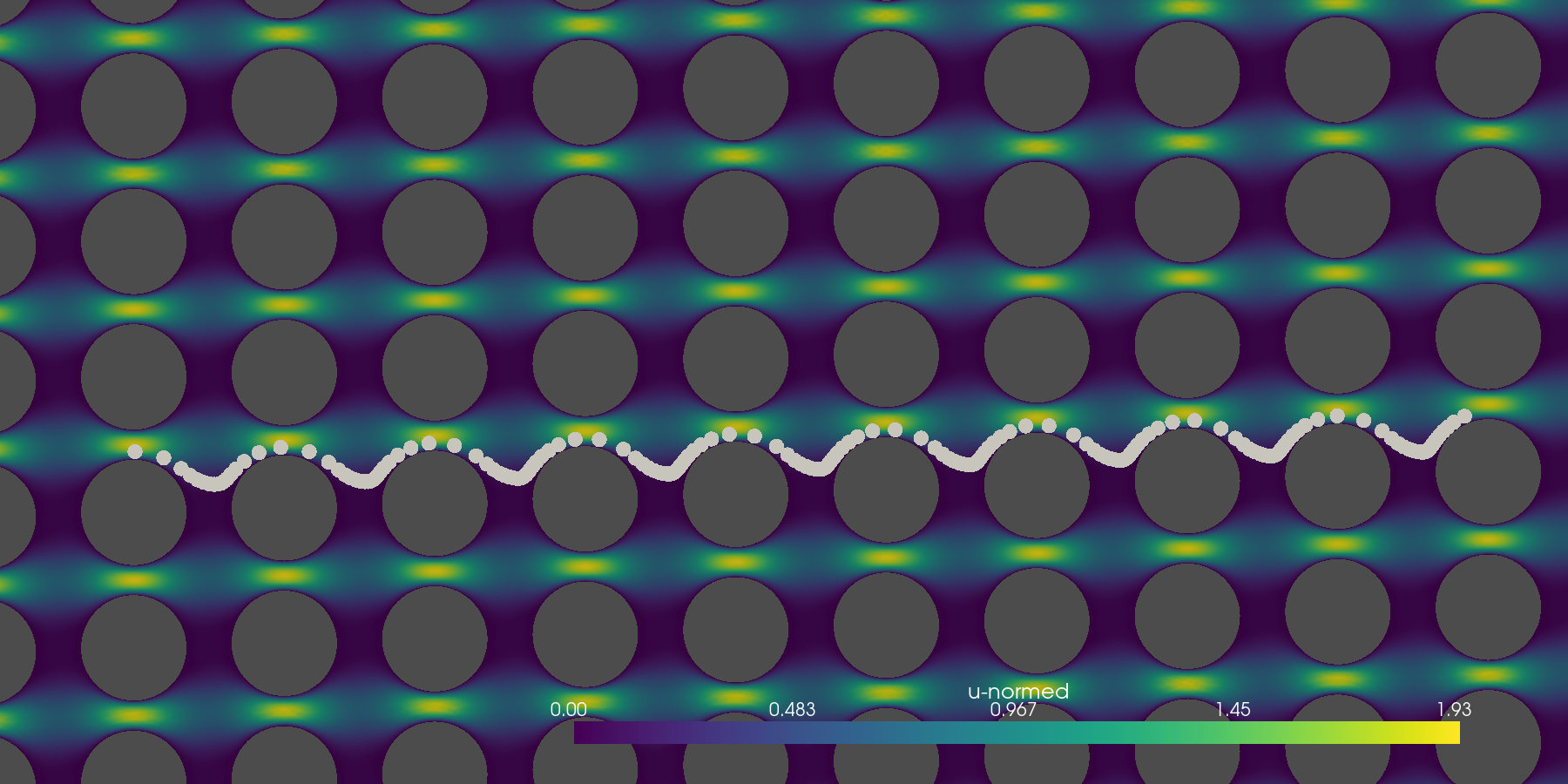}
    \hspace{0.5cm}
    \includegraphics[height=4cm]{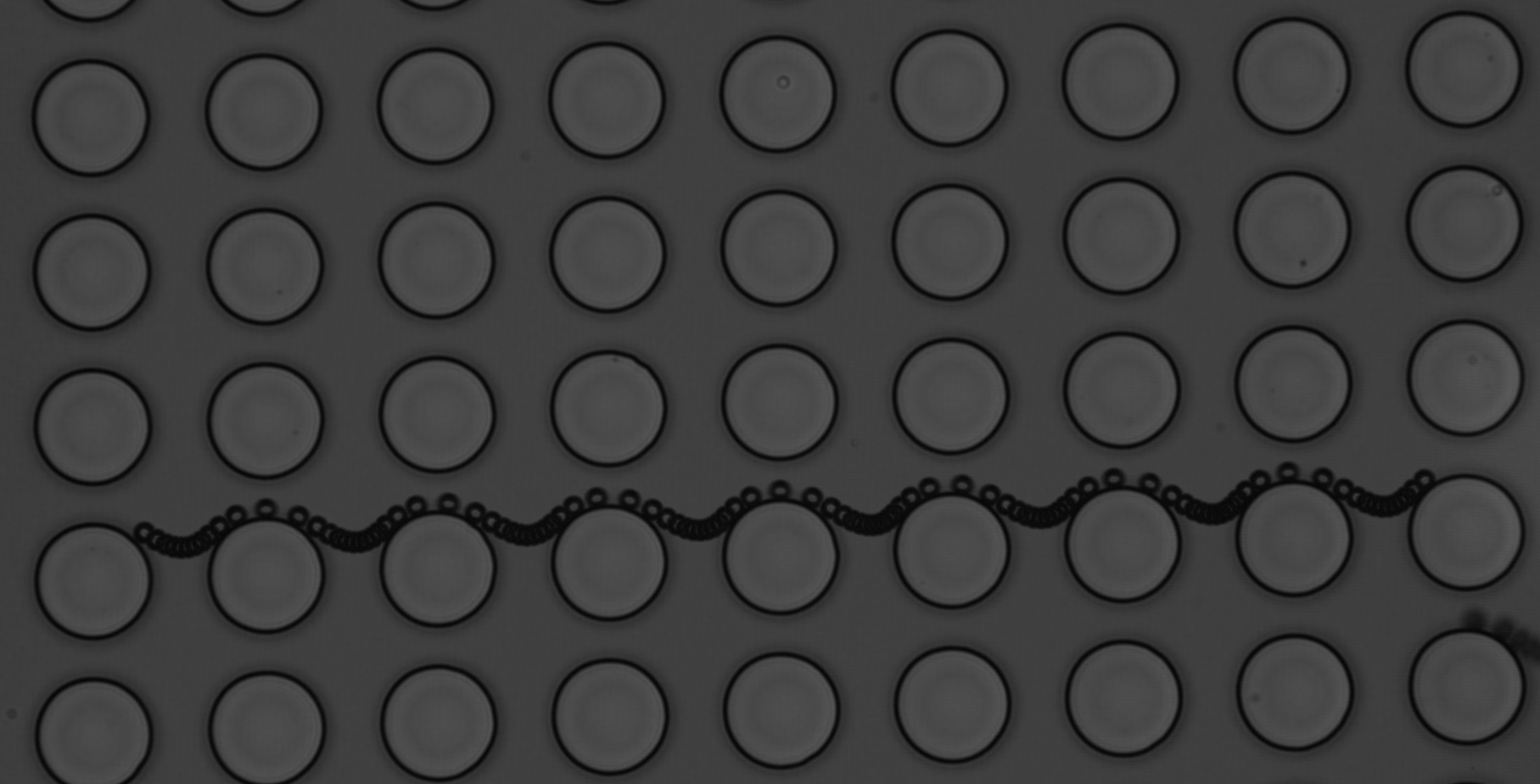}\\
    \caption{Comparisons of the lateral migration of particles
    between numerical (\textbf{left column}) and physical (\textbf{right column}) experiments
    for the diameter of particle D$=4$ (\textbf{top row}) and
    D$=5$ (\textbf{bottom row}), which illustrate that both numerical and
    physical results show the same critical diameter D$_c\in(4,5)$.}\label{fig:DLD-real}
\end{figure}

Then, we carry out Algorithm \ref{algorithm1} to solve the FSCI problem
occurring in the depicted DLD domain. Numerical results and their
comparisons with physical experiments 
are shown in Figure \ref{fig:DLD-real}, where the top row displays
the lateral migration of the particle of radius $r=2~\mu$m,
and the bottom row shows that of the particle of radius $r=2.5~\mu$m,
while the left column displays the numerical results versus
the physical experimental results that are shown in the right column.
We can see that both the numerical and physical
experimental results show a very similar lateral migration of the particle,
illustrating that the particle of radius $r=2~\mu$m travels in the zigzag mode
along the surface of pillar obstacles while the particle of radius $r=2.5~\mu$m
is bumped to the pillar and displaces laterally to the next streamline
without doing a zigzag move. Thus by the theory of DLD critical diameter,
we know that $4~\mu$m$~<~$ D$_c<5~\mu$m in this DLD microchip,
showing a good agreement between the numerical and physical experiments.

\section{Conclusion and future work}
In this paper, we develop a nonlinearly coupled system of
partial differential equations for one kind of fluid-structure-contact interaction (FSCI)
problem involving the interactional effect between the fluid and structure,
and the collision effect between the structure and fluidic channel wall,
where the structural motion may undergo large translational and rotational
displacements and/or deformations. To numerically tackle such a FSCI problem,
we develop a stabilized mixed finite element method within
the frame of monolithic fictitious domain/immersed boundary (FD/IB) approach to solve for
a unified variable pair of velocity and pressure defined only on the background
(Eulerian) fluidic mesh nodes in the entire domain. We further design advanced numerical implementation
algorithms to deal with the structural collision by means of a fixed-point nonlinear
iteration, to solve the developed nonlinear finite element approximation
system by Newton's linearization, and to handle the interpolation process
between the background Eulerian fluidic mesh and the foreground updated-Lagrangian
structural mesh through elemental computations and numerical quadrature points in an elementwise manner.
Our proposed numerical methods are validated through two benchmark problems,
self-designed and physically designed microfluidic DLD problems, respectively, by checking numerical convergence rates
and/or comparing with benchmark/physical experimental results,
which all return good agreements. Numerical simulations carried out by the
developed unified-field, monolithic fictitious domain-based
mixed finite element method (UFMFD-FEM) can help to optimize the design of cascaded
filter DLD microchips in an accurate and efficient fashion
for the sake of practically isolating circulating
tumor cells from blood cells in the blood fluid.

We will extend our numerical study to a rigorous theoretical analysis of the
proposed UFMFD-FEM in our next work as indicated in Remark \ref{rmk:theory}.
In addition, in our next paper we will also develop more advanced numerical
algorithms for the proposed UFMFD-FEM on the aspect of efficiency improvement,
such as the domain decomposition method (DDM), the adaptive time-step method,
the locally adaptive mesh refinement method, as well as the multi-timescale approach.
On the other hand, we will investigate a possible dimension error brought
by the two-dimensional FSCI simulation on an actually three-dimensional
DLD microchip whose thickness (about $30~\mu$m) is relatively thin though,
by doing 3D simulations for more realistic DLD problems
since our developed FSCI modeling approach can be seamlessly extended to high dimensions.
Moreover, we will conduct more in-depth studies on
the modeling of structural collisions and its numerical techniques
in more accurate and robust manner, and investigate more complicated cases
of collision effects amongst multiple particles that is beyond the
current contact phenomenon in this paper between a single particle and the
fluidic channel wall.

\section*{Acknowledgments}
C. Wang was supported by National Natural Science Foundation
of China grant (No. 12171366), P. Sun was partially supported by a
grant from the Simons Foundation (MPS-706640), J. Xu was
supported in part by KAUST Baseline Research Fund.

\appendix

\section{Description of interpolation process between fluidic and structural meshes}\label{appendix:interpolation}
In this appendix, we give a detailed implementation description for
the interpolation-related inner product terms arising from the developed UFMFD-FEM 
which occur in the structural domain $\Omega_s^t$
associated with its triangulation $\mathcal{T}_{h_s}^t\
(t\in(0,T])$ and whose integrand functions involve the
mesh-dependent trial function $\bv_h\in\bV_h$ and/or its test
function $\tilde\bv_h\in\bV_{0,h}$ that are associated with the
fixed background fluidic mesh $\mathcal{T}_h$. To conveniently describe the
interpolation procedure for inner products over two overlapped
meshes $\mathcal{T}_h$ and $\mathcal{T}_{h_s}^t$, without loss of
generality, we only consider the following two general terms in the
form of
\begin{eqnarray}
\left( u, v\right)_{\Omega_s^t} \text{ and } \left(\nabla u,\nabla v
\right)_{\Omega_s^t},\label{example-interpolation}
\end{eqnarray}
where $u$ and $v$ are two general finite element functions associated with
$\mathcal{T}_h=\bigcup\limits_{i=1}^M e_{f,i}$ and $\mathcal{T}_{h_s}^t
=\bigcup\limits_{i=1}^{M_s} e_{s,i}^t$, respectively.
Thus, the two terms in (\ref{example-interpolation}) can be rewritten as
\begin{eqnarray}
\sum_{i=1}^{M_s}\left( u, v\right)_{e_{s,i}^t} \text{ and }
\sum_{i=1}^{M_s}\left(\nabla u,\nabla v \right)_{e_{s,i}^t}.
\end{eqnarray}
Moreover, if we adopt the following quadrature scheme with $N_q$
quadrature points $\left\{\hat\bx_{j}\right\}_{j=1}^{N_q}$ defined
in an element $\hat e$,
\begin{eqnarray*}
\int_{\hat{e}}\hat{f}(\hat{\bx}) d\hat{\bx} \approx \sum_{j=1}^{N_q}
\omega_j\hat{f}(\hat{\bx}_j),
\end{eqnarray*}
then we can approximate two inner product terms in
(\ref{example-interpolation}) as follows,
\begin{eqnarray}
\left( u, v\right)_{\Omega_s^t}
&=&\sum\limits_{i=1}^{M_s}\sum\limits_{j=1}^{N_q} \omega_j{u}({\bx}_{ij}){v}({\bx}_{ij}),\notag\\
\left(\nabla u,\nabla v \right)_{\Omega_s^t}
&=&\sum\limits_{i=1}^{M_s}\sum\limits_{j=1}^{N_q} \omega_j {\nabla u}({\bx}_{ij}){\nabla v}({\bx}_{ij}),\notag
\end{eqnarray}
where $\bx_{ij}$ denotes the coordinate of the $j$-th quadrature
point $\hat{\bx}_j$ in the $i$-th structural element $e_{s,i}^t$.

Note that $v(\bx)$ is defined on $\mathcal{T}_{h_s}^t$, which
implies
that $v(\bx_{ij})$ and $\nabla v(\bx_{ij})$ can be computed by 
the following isoparametric transformations,
\begin{eqnarray}
    v(\bx_{ij})=\sum_{m=1}^{N_{s,p}}v_{i,m}\psi_{i,m}(\bx_{ij}),\quad
    \nabla v(\bx_{ij})=\sum_{m=1}^{N_{s,p}}v_{i,m}\nabla
    \psi_{i,m}(\bx_{ij}),
\end{eqnarray}
where $\{\psi_{i,m}\}_{m=1}^{N_{s,p}}$ are $N_{s,p}$ nodal basis
functions in the structural element $e_{s,i}^t\subset
\mathcal{T}_{h_s}^t$, and, the restriction of $v$ in $e_{s,i}^t$ is
denoted by $\sum\limits_{m=1}^{N_{s,p}}v_{i,m}\psi_{i,m}$. On the other
hand, we use the following three steps to compute $u(\bx_{ij})$ and
$\nabla u(\bx_{ij})$: (i) for each quadrature point $\bx_{ij}\
(i=1,\cdots,M_s,\ j=1,\cdots,N_q)$, we find the (first) fluidic
element $e_{f,k_{ij}}\subset\mathcal{T}_{h}$ that contains the
quadrature point $\bx_{ij}$. Note that the element in
$\mathcal{T}_{h}$ that contains $\bx_{ij}$ may be more than one in
the sense that $\bx_{ij}$ coincides with the vertices of
$\mathcal{T}_{h}$; (ii) Let $\{\phi_{k_{ij},l}\}_{l=1}^{N_{f,p}}$ be
$N_{f,p}$ nodal basis functions in the fluidic element
$e_{f,k_{ij}}\subset \mathcal{T}_h$.
The values of $\phi_{k_{ij},l}(\bx_{ij})$ and  $\nabla \phi_{k_{ij},l}(\bx_{ij})$ can also be computed by 
an isoparametric transformation, where the reference element
$\hat{e}_{f,k_{ij}}$ of $e_{f,k_{ij}}$ is used and the reference
quadrature point of $\bx_{ij}$ in $\hat{e}_{f,k_{ij}}$ is computed;
(iii) the desired values, $u(\bx_{ij})$ and $\nabla u(\bx_{ij})$, are
computed by
\begin{eqnarray}
    u(\bx_{ij}) = \sum_{l=1}^{N_{f,p}}u_{k_{ij},l}\phi_{k_{ij},l}(\bx_{ij}),\qquad
    \nabla u(\bx_{ij}) = \sum_{l=1}^{N_{f,p}}u_{k_{ij},l}\nabla \phi_{k_{ij},l}(\bx_{ij}),
\end{eqnarray}
where the restriction of $u$ in $e_{f,k_{ij}}\subset \mathcal{T}_h$
is denoted by  $\sum\limits_{l=1}^{N_{f,p}} u_{k_{ij},l}\phi_{k_{ij},l}$ .

Therefore, we have
\begin{eqnarray}
    \left( u, v\right)_{\Omega_s^t}
    =&\sum\limits_{i=1}^{M_s}\sum\limits_{j=1}^{N_q}\sum\limits_{l=1}^{N_{f,p}}\sum\limits_{m=1}^{N_{s,p}}\omega_j u_{k_{ij},l}v_{i,m}\phi_{k_{ij},l}(\bx_{ij}) {\psi_{i,m}}({\bx}_{ij}),\label{eq_54}\\
    \left(\nabla u,\nabla v \right)_{\Omega_s^t}
    =&\sum\limits_{i=1}^{M_s}\sum\limits_{j=1}^{N_q} \sum\limits_{l=1}^{N_{f,p}}\sum\limits_{m=1}^{N_{s,p}} \omega_j u_{k_{ij},l}v_{i,m}\nabla \phi_{k_{ij},l}(\bx_{ij}) {\nabla \psi_{i,m}}({\bx}_{ij}).\label{eq_55}
\end{eqnarray}
where the element index $k_{ij}$ of $\mathcal{T}_h$ depends on
the structural element index $i$ of $\mathcal{T}_{h_s}^t$ and on the quadrature
point index $j$. We remark that $k_{ij}$ may take different values
even for the cases that the element index $i$ is the same but the
quadrature point index $j$ is different.

Furthermore, to compute and assemble the elementwise finite element
matrix, we only need to take $v= \psi_{i,m}$ as the test function in
(\ref{example-interpolation}). Thus, \eqref{eq_54} and \eqref{eq_55}
can be further simplified as
\begin{eqnarray}
\left( u, v\right)_{\Omega_s^t}
=&\sum\limits_{i=1}^{M_s}\sum\limits_{j=1}^{N_q}\sum\limits_{l=1}^{N_{f,p}}\omega_j u_{k_{ij},l}\phi_{k_{ij},l}(\bx_{ij}) {\psi_{i,m}}({\bx}_{ij}),\label{eq_54B}\\
\left(\nabla u,\nabla v \right)_{\Omega_s^t}
=&\sum\limits_{i=1}^{M_s}\sum\limits_{j=1}^{N_q}
\sum\limits_{l=1}^{N_{f,p}} \omega_j u_{k_{ij},l}\nabla
\phi_{k_{ij},l}(\bx_{ij}) {\nabla
\psi_{i,m}}({\bx}_{ij}),\label{eq_55B}
\end{eqnarray}
for $m=1,\cdots,N_{s,p}$, where $u_{k_{ij},l}$ is called the degree
of freedoms (DOFs) defined on the vertices of $\mathcal{T}_h$.

\bibliographystyle{plain}
\bibliography{DLD,DLMFD-FSI,FSI,publicationlist4Sun}

\end{document}